\documentclass[12pt]{amsart}
\usepackage{amssymb}
\renewcommand{\a}{\alpha}
\renewcommand{\o}{\omega}
\renewcommand{\b}{\beta}
\renewcommand{\d}{\delta}
\newcommand{\D}{\Delta}
\newcommand{\e}{\varepsilon}
\newcommand{\g}{\gamma}
\renewcommand{\l}{\lambda}
\renewcommand{\L}{\Lambda}
\newcommand{\var}{\varphi}
\newcommand{\s}{\sigma}
\renewcommand{\th}{\theta}
\renewcommand{\O}{\Omega}
\newcommand{\z}{\zeta}
\newcommand{\Th}{\Theta}

\newcommand{\G}{\Gamma}

\renewcommand{\i}{\infty}
\newcommand{\p}{\partial}
\newcommand{\wt}{\widetilde}
\newcommand{\wh}{\widehat}
\newcommand{\x}{\times}
\newcommand{\<}{\langle}
\renewcommand{\>}{\rangle}

\newcommand{\Ad}{\mathop{\mathrm{Ad}}\nolimits}
\newcommand{\id}{\mathop{\mathrm{id}}\nolimits}
\renewcommand{\Im}{\mathop{\mathrm{Im}}\nolimits}
\newcommand{\Tr}{\mathop{\mathrm{Tr}}\nolimits}
\newcommand{\Gr}{\mathop{\mathrm{Gr}}\nolimits}

\newcommand{\cA}{{\mathcal A}}
\newcommand{\cB}{{\mathcal B}}
\newcommand{\cC}{{\mathcal C}}
\newcommand{\cD}{{\mathcal D}}
\newcommand{\cH}{{\mathcal H}}
\newcommand{\cL}{{\mathcal L}}
\newcommand{\cO}{{\mathcal O}}
\newcommand{\cQ}{{\mathcal Q}}
\newcommand{\cR}{{\mathcal R}}
\newcommand{\cT}{{\mathcal T}}
\newcommand{\cU}{{\mathcal U}}
\newcommand{\cV}{{\mathcal V}}
\newcommand{\cX}{{\mathcal X}}
\newcommand{\cZ}{{\mathcal Z}}

\newcommand{\bC}{{\mathbb C}}
\newcommand{\bH}{{\mathbb H}}
\newcommand{\bN}{{\mathbb N}}
\newcommand{\bP}{{\mathbb P}}
\newcommand{\bR}{{\mathbb R}}
\newcommand{\bT}{{\mathbb T}}
\newcommand{\bZ}{{\mathbb Z}}

\newcommand{\fM}{{\mathfrak M}}
\newcommand{\fS}{{\mathfrak S}}

\newcommand{\bpm}{\begin{pmatrix}}
\newcommand{\epm}{\end{pmatrix}}

\begin{document}

\title[Free Analysis Questions II]{Free Analysis Questions II: \\
the Grassmannian Completion and \\
the Series Expansions at the Origin}
\author{Dan-Virgil Voiculescu}

\begin{abstract}
The fully matricial generalization in part~I, of the difference quotient derivation on holomorphic functions, in which ${\mathbb C}$ is replaced by a Banach algebra $B$, is extended from the affine case to a Grassmannian completion.  The infinitesimal bialgebra duality, the duality transform generalizing the Stieltjes transform and the spectral theory with non-commuting scalars all extend to this completion.  The series expansions of fully matricial analytic functions are characterized, providing a new way to generate fully matricial functions.
\end{abstract}

\maketitle
\allowdisplaybreaks

\section{Introduction}
\label{sec1}

In our first paper in this series \cite{14}, we produced a fully matricial object which generalizes from the complex numbers $\bC$ to a unital Banach algebra $B$ the structure on $\bC$ in which on each open subset of $\bC$ we are given the bialgebra of scalar analytic functions with the comultiplication defined by the difference quotient derivation.  One main theme of the present paper can be described as the extension of the construction from the affine line to the Riemann sphere.  This means that with $\bC$ replaced by $B$, we will exhibit a fully matricial object generalizing the Riemann sphere and the bialgebras of analytic functions on open sets with respect to the difference quotient.  This will be achieved by tying together certain Grassmannians for the matrices $\fM_n(B)$.

Among the other main questions we study are the series expansions at the origin for fully matricial functions.  Using our results on the series expansions we are able to clarify the structure of fully matricial analytic functions.  For the bounded functions in the matricial analogue of the disk we show that the non-commutative $H^{\i}$-spaces which arise are connected to certain non-self-adjoint algebras in full free product $C^*$-algebras ($\fM_k *_{\bC} C(\bT)$ in case $B = \fM_k$ and the full free group $C^*$-algebra in case $B = \bC^n$).  A further feature of these matricial holomorphic functions is that they naturally have a large $N$ limit which lives in certain non-self-adjoint subalgebras of free group $II_1$-factors.

There are several reasons for passing to the Grassmannian completion.  One reason is the duality transform for the coalgebra of the free difference quotient (\cite{14} see also \cite{13}).  We used only a part of the corepresentations (i.e., group-like elements) which arise from generalized matricial resolvents.  Clearly, this means the matricial analytic functions of the duality transform have analytic extensions to a larger natural domain of definition and we should try to find these extensions.

Another motivation is provided by spectral theory.  Duality for the free difference quotient can be viewed as an extension of the spectral analysis of an operator from the case of complex scalars $\bC$ to the case of an algebra $B$ of scalars which may not commute with the operator.  To handle the behavior at infinity of the generalized resolvents it is natural to use the Grassmannian.  This is also convenient for dealing with unbounded self-adjoint operators.  Note that similar motives underlie the use of Grassmannians in operator theory problems of electrical engineering (see \cite{5}).  

The point at infinity of the Grassmannian completion is also well suited for dealing with the fully matricial version of the operator-valued $R$-transform (\cite{11}).  The result on the series expansion when applied to this fully matricial $R$-transform yields a noncommutative series $R$-transform, which in view of formal analogies coincides with the unsymmetrized $R$-transform of \cite{4a} (where also the relation to free comulants \cite{8a} is discussed).

A central result of our paper, the proof of which is scattered through several sections, can be stated roughly as follows.

\bigskip
\noindent
{\bf Fact.} If $B$ is a unital complex Banach algebra and $\O$ a fully matricial $B$-set of the Grassmannian, then the fully matricial scalar analytic functions $\cA(\O)$ form a ``topological'' infinitesimal bialgebra (\cite{7}, \cite{1}) (or GDQ ring in the terminology which we used in \cite{13}, \cite{14}) $(\cA(\O),\mu,\wt{\p})$ with multiplication $\mu: \cA(\O;\O) \to \cA(\O)$ and comultiplication-derivation $\wt{\p}: \cA(\O) \to \cA(\O;\O)$ and with coassociativity requirements stated in terms of $\cA(\O;\O;\O)$.  Moreover, there is a coderivation $\L: \cA(\O) \to \cA(\O)$ with respect to $\wt{\p}$ such that $\L - \id$ is a derivation of $\cA(\O)$.  If $B$ is a Banach algebra with involution and $\O = \O^*$, then there is an involution on $\cA(\O)$, compatible with $(\cA(\O),\mu,\wt{\p},\L)$ (that is, satisfying the ``topological'' version of Definition~{\em \ref{sec5.1}} of \cite{13}).

\bigskip
For the precise definitions of certain terms in the previous statement, it is necessary to read the paper beyond the introduction.  Note also, that the fully matricial affine space, to which the considerations in \cite{14} were restricted, is an open subset of the Grassmannian.  In case $B =\bC$, $\wt{\p}$ is the difference quotient $(z_1-z_2)^{-1}(f(z_1)-f(z_2))$ acting on analytic functions on an open subset of the Riemann sphere.  The difference quotient $\wt{\p}$ is not invariant under automorphisms of $\bP^1(\bC)$ and there is a certain vanishing of $\wt{\p}f$ at infinity.

Another main result is the extension of the duality transform to the Grassmannian context.  Part of this can also be seen as the extension of the related non-commutative spectral theory, with particular emphasis on the cases of unitary and self-adjoint operators.  While the properties of the duality transform relative to an operator, or graph of operator, and an ``algebra of scalars'' depend on whether certain technical conditions are satisfied, we show in the unitary and self-adjoint cases that one can always pass to the universal such operators which satisfy these requirements and use the corresponding universal duality transforms.

We have also included in this second paper new results in the affine context, like the sub-bialgebra of polynomial functions on the fully matricial affine space.  This clarifies the action of the generalized difference quotient and underlies the series expansions at the origin and the beginnings of a corresponding study of functions in the analogue of the disk.

It is not clear whether the non-commutative $H^{\i}$-spaces which we encounter here can benefit from the recent progress on the non-commutative $H^{\i}$-spaces arising from generalizations of the shift operator and of Fock spaces (see for instance \cite{8} and the references therein).

There is a certain peculiarity of the Grassmannian completion that should be pointed out.  In the affine framework of the preceding paper (\cite{14}), for most purposes it was sufficient to assume that $B$ is just a Banach space with a distinguished unit vector and sometimes to add to this a conjugate-linear involution.  The Grassmannian completion uses the fact that $B$ is an algebra, that is that there is a multiplication operation.  It is conceivable that while the multiplication may be necessary, its associativity could be relaxed.

The paper has seventeen sections, including the introduction and there are also two Appendices.  Some preliminaries from \cite{14} are briefly recalled in section~\ref{sec2}, but we could only cover a part of the instances where it is necessary to go back to \cite{14}.  Section~\ref{sec3} introduces the fully matricial Grassmannian and its fully matricial sets and fully matricial functions.  Grassmannian resolvents and resolvent sets are the subject of section~\ref{sec4}.  The long section~\ref{sec5} is devoted to the definition of the generalization of the difference quotient derivation to the fully matricial $B$-Grassmannian and to checking its basic properties.  Section~\ref{sec6} is also one of the longer sections of the paper and deals with the Grassmannian resolvent equation and duality transform.  The important algebras $\cL\cR(\pi;B)$ of coefficients of Grassmannian resolvents are also introduced in section~\ref{sec6}.  The polynomial bialgebra on the fully matricial affine space and some basic remarks about the fully matricial affine space are the subject of section~\ref{sec7}.  We then return to the Grassmannian in section~\ref{sec8} and we define the coderivation $\L$ and check its properties with respect to the bialgebra structure and duality transform.  In section~\ref{sec9} we turn to the involution on the Grassmannian and to its properties with respect to the bialgebra structure and duality transform.  (This is an instance when the arguments in the affine case are substantially shorter.)  We take up dual positivity and duality transforms of positive functionals in the Grassmannian context in section~\ref{sec10}.  The generalizations to our context of certain spectrally  important sets such as half-planes and the unit ball leads to stably matricial sets, which are more general than fully matricial sets, and to which we devote section~\ref{sec11}.  We also introduce in section~\ref{sec11} the mixed unit balls and the mixed half-planes, which are always present in the fully matricial resolvents of unitary and respectively of hermitian elements in addition to the usual unit ball and upper and lower half-planes.  We conclude section~\ref{sec11} with results about the use of stably matricial sets for obtaining fully matricial inverse function results.  Section~\ref{sec12} collects some further remarks on the duality transforms for unitary and hermitian elements and introduces the universal unitary and hermitian Grassmmannian elements.  We show that certain technical assumptions for the duality transform appearing in previous sections hold in the universal case.  Section~\ref{sec13} establishes the form of the series expansions at the origin for fully matricial analytic functions and their behavior under composition of such functions.  We should mention that using automorphisms of the fully matricial Grassmannian this also can be used for instance for the series expansions for the point at infinity or for ``$B$-points of the affine space'' that is for fully matricial sets $(\{b \otimes I_n\})_{n \in \bN}$.  In section~\ref{sec14} we find the analogue in case $B = \fM_k$ for the formulae connecting Fourier coefficients of boundary values on the circle to coefficients of the Taylor series at the origin.  The analogue of the integrations on the circle is here the limit of the integrations over the unitary groups $U(\fM_N(\fM_k))$.  We also show that the analogue of the $H^{\i}$-algebra in this case is related to a subalgebra of the full free product $\fM_k *_{\bC} C(\bT)$.  In section~\ref{sec15} we use the results of the preceding section to show in the case $B = \fM_k$ the existence of a large $N$ limit of the functions which are in the analogue of $H^{\i}$ of the unit disk.  The large $N$ limit is in a non-self-adjoint subalgebra of a $II_1$-factor.  Section~\ref{sec16} deals with the analogues in case $B = \bC^k$ of the results of sections~\ref{sec14} and \ref{sec15}.  Here the related $C^*$-algebra is the full $C^*$-algebra of a free group $C^*(F_k)$ and the related von~Neumann algebra is the free group factor $L(F_k)$.  In section~\ref{sec17}, we give a construction of pathological fully matricial analytic functions in the absence of boundedness conditions.

For a better perspective we give in Appendix~I the classical case of the difference quotient bialgebra duality on the Riemann sphere.  Checking the duality relations amounts to familiar applications of the Cauchy integral formula.  In a second appendix (Appendix~II) we briefly explain how the $B$-valued $R$-transform (\cite{11}) gives rise naturally to a fully matricial $R$-transform and this in turn to a noncommutative series $R$-transform which is precisely the unsymmetrized $R$-transform of \cite{4a}.  In this way all $B$-moments, not only the symmetric ones, are encoded in a generalization of the analytic $R$-transform.

\section{Preliminaries}
\label{sec2}

\subsection{}
\label{sec2.1} Let $E$ be a complex Banach algebra with $1$ and $1 \in B \subset E$ a Banach subalgebra and let $X \in E$.  The Banach algebra $B$ will play the role of ``scalars'' in the spectral analysis of $X$.  Note that $X$ and $B$ do not commute in general.  Also for the considerations in this section we could have replaced $B$ by a Banach subspace $\cU$, $1 \in \cU \subset E$.  Later on we will be especially interested in the case where $E = \cB(\cH)$ is the algebra of bounded operators on a Hilbert space $\cH$, $B$ a $C^*$-subalgebra of $E$ and $X$ is a self-adjoint operator.  The material in this section (with slight adaptations) is from \cite{14}.

\subsection{}
\label{sec2.2} The $n \x n$ matrices with entries in $B$, denoted $\fM_n(B)$ can be identified with $\fM_n \otimes B$, where $\fM_n$ is short for $\fM_n(\bC)$ and the tensor product is over $\bC$.  Then $I_n \otimes X \in \fM_n \otimes E \simeq \fM_n(E)$ is the diagonal matrix with diagonal entries equal $X$, which we will also denote sometimes by $\underset{n}{\underbrace{X \oplus \dots \oplus X}}$.

For each $n \ge 1$, the $n$-th $B$-resolvent set of $X$, is the set
\[
\rho_n(X;B) = \{\b \in\fM_n(B) \mid I_n\otimes X - \b \mbox{ invertible}\}.
\]
The sequence of sets $\rho(X;B) = (\rho_n(X;B))_{n \in \bN}$ will be called the full $B$-resolvent set of $X$.

\subsection{}
\label{sec2.3} The full $B$-resolvent set of $X$ is an example of a fully matricial $B$-set.  More generally, if $\cU$ is a Banach space over $\bC$, a fully matricial $\cU$-set $\O$ is a sequence $\O = (\O_n)_{n \in \bN}$, $\O_n \subset \fM_n(\cU)$ so that $\O_{m+n} \cap (\fM_m(\cU) \oplus \fM_n(\cU)) = \O_m \oplus \O_n$ and
\[
(\Ad S \otimes I_{\cU})(\O_n) = \O_n
\]
if $S \in GL(n;\bC) \subset \fM_n(\bC)$. Here $\Ad S$ denotes the inner automorphism $T \rightsquigarrow STS^{-1}$ of $\fM_n$ and $I_{\cU}$ is the identity operator on $\cU$.

\subsection{}
\label{sec2.4} On each set $\rho_n(X;B)$ there is an analytic function
\[
R_n(X;B): \rho_n(X;B) \to \fM_n(E)
\]
defined by
\[
R_n(X;B)(\b) = (I_n \otimes X - \b)^{-1}.
\]
The sequence $R(X;B) = (R_n(X;B))_{n \in \bN}$ is called the full $B$-resolvent of $X$.

\subsection{}
\label{sec2.5} The full $B$-resolvent of $X$ is an example of a fully matricial $E$-valued function.  More generally if $\O = (\O_n)_{n \in \bN}$ is a fully matricial $\cU$-set and $\cV$ is another Banach space, then a fully matricial $\cV$-valued function on $\O$ is a sequence $f = (f_m)_{n \in \bN}$ where 
\[
f_n: \O_n \to \fM_n(\cV)
\]
and
\[
f_{m+n} \mid \O_m \oplus \O_n = f_m \oplus f_n
\]
and
\[
(\Ad S \otimes I_v) \circ f_n = f_n\circ (\Ad S \otimes I_n)|_{\O_n},
\]
when $S \in GL(n;\bC)$.

\section{The Grassmannian completion}
\label{sec3}

\subsection{}
\label{sec3.1} For the affine theory it is sufficient to consider fully matricial functions and sets with respect to a Banach space with a specified non-zero element.

To construct the Grassmannian completion the Banach space must be upgraded to a Banach algebra over $\bC$ with unit.

Throughout this section $B$ is a Banach algebra over $\bC$ with unit $1 \in B$.  The Grassmannian $Gr_n(B)$ is defined to be the set
\[
GL_2(\fM_n(B))/\wt{\l n}
\]
where
\[
\bpm
a & b \\
c & d
\epm \wt{\l n} \bpm
a' & b' \\
c' & d'
\epm
\]
if there is $t \in GL_1(\fM_n(B))$ so that $bt = b'$, $dt = d'$.  We may describe alternatively $Gr_n(B)$ as the set of right $\fM_n(B)$ submodules of $\fM_n(B) \oplus \fM_n(B)$ which are isomorphic to $\fM_n(B)$ and have a complement of the same kind.  The submodule corresponding to $\left.\bpm a & b \\ c & d \epm\right/\wt{\l n}$ is then
\[
\left\{ \bpm
b \xi \\
d \xi
\epm \in \fM_n(B) \oplus \fM_n(B) \mid \xi \in \fM_n(B)\right\}.
\]
Note also that the equivalence relation $\wt{\rho n}$ on $GL_2(\fM_n(B))$
\[
\bpm
a & b \\
c & d
\epm \wt{\rho n} \bpm
a' & b' \\
c' & d'
\epm
\]
if $at = a'$, $ct = c'$ for some $t \in GL_1(\fM_n(B))$, provides an alternative realization of $Gr_n(B)$.

\subsection{}
\label{sec3.2} If
\[
\pi_j = \left.\bpm a_j & b_j \\ c_j & d_j \epm\right/\wt{\l n_j} \in Gr_{n_j}(B)
\]
$j = 1,2$, we define
\[
\pi_1 \oplus \pi_2 = \left.\bpm
a_1 \oplus a_2 & b_1 \oplus b_2 \\
c_1 \oplus c_2 & d_1 \oplus d_2
\epm\right/ \wt{\l(n_1+n_2)} \in Gr_{n_1+n_2}(B).
\]
Also if
\[
\pi = \left.\bpm 
a & b \\
c & d 
\epm\right/\wt{\l n} \in Gr_n(B)
\]
and $g \in GL_2(\fM_n(B))$ we see that $g\pi = \left.\left( g \bpm a & b \\ c & d \epm\right)\right/\l n$ is well-defined (i.e., depends only on $g$ and $\pi$).  We will use extensively the action of $GL_1(\fM_n(\bC))$ $(\simeq GL_n(\bC))$, denoted $s \cdot \pi$ where
\[
s \cdot \pi = \bpm s & 0 \\ 0 & s \epm \pi
\]
and
\[
\bpm
s & 0 \\
0 & s
\epm \in GL_2(\fM_n(\bC)) \subset GL_2(\fM_n(B)).
\]
Clearly
\[
s \cdot \pi = \left.\bpm
sa & sb \\
sc & sd
\epm\right/\wt{\l n} = \left. \bpm
sas^{-1} & sbs^{-1} \\
scs^{-1} & sds^{-1}
\epm\right/ \wt{\l n}.
\]

\subsection{}
\label{sec3.3} By definition a fully matricial $B$-set of the Grassmannian is a sequence of sets $\O = (\O_n)_{n \in \bN}$ where $\O_n \subset Gr_n(B)$ and the following two conditions are satisfied:
\begin{itemize}
\item[(i)] $\O_{m+n} \cap (Gr_m(B) \oplus Gr_n(B)) = \O_m \oplus \O_n$
\item[(ii)] $s \cdot \O_n = \O_n$ if $s \in GL(n;\bC)$ $(= GL_1(\fM_n(\bC)))$.
\end{itemize}

\subsection{}
\label{sec3.4} If $\O = (\O_n)_{n \in \bN}$ is a fully matricial $B$-set of the Grassmannian and $\cU$ is a Banach space over $\bC$, then a fully matricial $\cU$-valued function on $\O$ is a hierarchy of functions $f = (f_n)_{n \in \bN}$ where $f_n: \O_n \to \fM_n(\cU)$ are such that:
\begin{itemize}
\item[(i)] $f_{n_1+n_2}(\pi_1 \oplus \pi_2) = f_{n_1}(\pi_1) \oplus f_{n_2}(\pi_2)$.
\item[(ii)] $f_n(s \cdot \pi) = sf_n(\pi)s^{-1}$ if $s \in GL(n;\bC)$.
\end{itemize}

\subsection{}
\label{sec3.5} We can also define if $A,B$ are unital Banach algebras and if $\O,\Th$ are fully matricial $A$ and respectively $B$ sets of the Grassmannian, fully matricial maps $F: \O \to \Th$ as sequences $(F_n)_{n\ge 1}$ $F_n: \O_n \to \Th_n$ satisfying $F_{n_1+n_2}(\pi_1 \oplus \pi_2) = F_{n_1}(\pi_1) \oplus F_{n_2}(\pi_2)$ and $F_n(s \cdot \pi) = s \cdot F_n(\pi)$ $(s \in GL(n,\bC))$.

\subsection{}
\label{sec3.6} To emphasize the distinction between fully matricial $B$-sets of the Grassmannian and fully matricial $B$-sets we will sometimes refer to the latter as affine fully matricial $B$-sets.

\subsection{}
\label{sec3.7} The Grassmannian $Gr_n(B)$ is naturally a Banachic complex analytic manifold, so that a map defined on an open set $\o$ of $Gr_n(B)$ is analytic iff it is analytic when composed with the surjection
\[
p_n^{-1}(\o) \to \o
\]
where $p_n$ is the surjection
\[
GL_2(\fM_n(B) \to Gr_n(B).
\]
This can be checked along standard lines using charts
\[
\left\{\left.\left( \bpm
a & b \\
c & d 
\epm \bpm
1 & f \\
0 & 1
\epm \right)\right/\wt{\l n} \mid f \in \fM_n(B),\ \|f\| < \e\right\}
\]
for suitably small $\e > 0$ (a reference for analytic functions on Banach space is \cite{6}).

Remark that the largest affine fully matricial $B$-set $(\fM_n(B))_{n \in \bN}$ can be identified with a fully matricial $B$-set of the Grassmannian $\O = (\O_n)_{n \in \bN}$ where
\[
\O_n = \left\{ \left. \bpm
1 & f \\
0 & 1
\epm \right/\wt{\l n} \in Gr_n(B) \mid f \in \fM_n(B)\right\}.
\]
We will call $\O_n$ the affine part of $Gr_n(B)$.  This identification also turns every affine fully matricial $B$-set into a fully matricial $B$-set of the Grassmannian.  Fully matricial analytic maps of open affine matricial $B$-sets then become analytic fully matricial maps of the fully matricial $B$-set of the Grassmannian.

\subsection{}
\label{sec3.8} Note also, in a converse direction, that $\O = (\O_n)_{n \in \bN}$ is a fully matricial $B$-set of the Grassmannian iff $p^{-1}(\O) = (p_n^{-1}(\O_n))_{n \in \bN}$ is an affine fully matricial $\fM_2(B)$-set.  Here $p^{-1}(\O)$ is identified with a subset of $(GL_n(\fM_2(B)))_{n \in \bN}$ via the appropriate identifications
\[
GL_n(\fM_2(B)) \simeq GL_2(\fM_n(B)).
\]

\subsection{}
\label{sec3.9} Note that for $K$-theory reasons for instance $Gr_n(B)$ may not be connected and thus $\O_n$ the affine part of $Gr_n(B)$ not only will not be dense in $Gr_n(B)$ in this case, but even more it will not even be a set of uniqueness for analytic functions (consider for instance functions which are constant on the connected components).  This occurs if for instance there is $u \in GL(2;B)$ such that $[u]$ its $K_1$-class is not in the subgroup $\{[v] \in K_1(B) \mid v \in GL(1;B)\}$. Indeed, since $Gr_2(B) = GL(2;B)/T$ where $T = \left\{ \bpm a & 0 \\ b & d \epm \in GL(2;B) \mid a,d \in GL(1,B)\right\}$ we will have that $u/\wt{\l 2}$ is not in the connected component of $\left. \bpm 0 & 1 \\ 1 & 0 \epm\right/\wt{\l 2}$ in $Gr_2(B)$.

Note also that if $\bpm a & b \\ c & d \epm \in GL(2;B)$ is such that $b$ is a nonunitary isometry, then $\left. \bpm a & b \\ c & d \epm\right/\wt{\l 2}$ will not be in the closure of the affine part of $Gr_2(B)$.  If in $B$ the invertible elements are dense (i.e., $B$ has topological stable rank~$1$), then the affine part of $Gr_n(B)$ is dense in $Gr_n(B)$.

\subsection{}
\label{sec3.10} Let $B_1,\dots,B_p$ be unital Banach algebras over $\bC$.  A fully multimatricial $(B_1,\dots,B_p))$-set of the Grassmannians is a family of sets $\O = (\O_{n_1,\dots,n_p})_{(n_1,\dots,n_p) \in \bN^p}$ where
\[
\O_{n_1,\dots,n_p} \subset Gr_{n_1}(B_1) \x \dots \x Gr_{n_p}(B_p)
\]
and for each $1 \le j \le p$ and $\pi_k \in Gr_{n_k}(B_k)$ $k \in \{1,\dots,{\hat j},\dots,p\}$,
\[
(\O_{n_1,\dots,n_{j-1},n,n_{j+1},\dots,n_p} \cap (\{(\pi_1,\dots,\pi_{j-1})\} \x Gr_n(B_j) \x \{(\pi_{j+1},\dots,\pi_p)\}))_{n \in \bN}
\]
is a fully matricial $B_j$-set of the Grassmannian.  A fully multimatricial $\cU$-valued function on $\O$ is a family of functions $f = (f_{n_1,\dots,n_p})_{(n_1,\dots,n_p) \in \bN^p}$ where $f_{n_1,\dots,n_p}: \O_{n_1,\dots,n_p} \to \fM_{n_1} \otimes \dots \otimes \fM_{n_p} \otimes \cU$ and for each fixed $n_1,\dots,n_{j-1},n_{j+1},\dots,n_p$ and $\pi_k$ as above $f_{n_1,\dots,n_{j-1},n,n_{j+1},\dots,n_p}$ gives rise to a fully matricial
\[
\fM_{n_1} \otimes \dots \otimes \fM_{n_{j-1}} \otimes \fM_{n_{j+1}} \otimes \dots \otimes \fM_{n_p} \otimes \cU
\]
valued functions on
\[
(\O_{n_1,\dots,n_{j-1},n,n_{j+1},\dots,n_p} \cap (\{(\pi_1,\dots,\pi_{j-1})\} \x Gr_n(B_j) \x \{(\pi_{j+1},\dots,\pi_p)\}))_{n \in \bN}.
\]

\section{Grassmannian resolvents and resolvent sets}
\label{sec4}

\subsection{}
\label{sec4.1} Let $E$ be a complex Banach algebra with $1$ an $1 \in B \subset E$ a Banach subalgebra.  Let further $\pi= \left. \bpm a & b \\ c & d \epm\right/\wt{\l 1} \in Gr_1(E)$.  We define the $n$-th Grassmannian $B$-resolvent set of $\pi$ to be the set $\wt{\rho}_n(\pi;B) = \{\s \in Gr_n(B) \mid \s$ is a complement of $\underset{\mbox{$n$-times}}{\underbrace{\pi \oplus \dots \oplus \pi}}\}$.  If $\s = \left. \bpm \a & \b \\ \g & \d \epm\right/\wt{\l n}$ where $\a,\b,\g,\d \in \fM_n(B)$, then $\s \in \wt{\rho}_n(\pi;B)$ is equivalent to requiring that
\[
\bpm
\b & \begin{matrix} b & & O \\ & \ddots & \\ O & & b \end{matrix} \vspace{1\jot} \\
\d & \begin{matrix} d & & O \\ & \ddots & \\ O & & d \end{matrix}
\epm \in GL(2n;E).
\]

It is easily seen that $\wt{\rho}(\pi;B) = (\wt{\rho}_n(\pi;B))_{n \in \bN}$ is a fully matricial $B$-set of the Grassmannian.  The direct sum property is obvious and the similarity property follows from the fact that $s \cdot (\underset{n}{\underbrace{\pi \oplus \dots \oplus \pi}}) = \underset{n}{\underbrace{\pi \oplus \dots \oplus \pi}}$ if $s \in GL(n;\bC)$ and $\s$ is a complement of $\underset{n}{\underbrace{\pi \oplus \dots \oplus \pi}}$ iff $s \cdot \s$ is a complement of $s \cdot (\pi \oplus \dots \oplus \pi)$.  We shall call $\wt{\rho}(\pi;B)$ the {\em full Grassmannian $B$-resolvent of $\pi$}.

\subsection{}
\label{sec4.2} On $\wt{\rho}_n(\pi;B)$ we define the $\fM_n(E)$-valued analytic function $\wt{\cR}_n(\pi;B)$.  {\em If $\s \in \wt{\rho}_n(\pi;B)$ and $\pi$ are like in the preceding subsection, let $\z \in \fM_n(E)$ be so that
\[
\bpm
\begin{matrix} b & & O \\ & \ddots & \\ O & & b \end{matrix} & \b \vspace{1\jot} \\
\begin{matrix} d & & O \\ & \ddots & \\ O & & d \end{matrix} & \d
\epm^{-1} = \bpm * & * \\ * & \z \epm ,
\]
Then we define $\wt{\cR}_n(\pi;B)(\s) = \b\z$.}  If $\tau \in GL_1(\fM_n(B))$, $t \in GL_1(E)$, then replacing $\b,\d,b,d$ by $\b\tau,\d\tau,bt,dt$ will lead to replacing $\z$ by $\tau^{-1}\z$.  Since $\b\z = (\b\tau)(\tau^{-1}\z)$ we see that $\wt{\cR}_n(\pi;B)(\s)$ is well-defined.  We will call $\wt{\cR}_n(\pi;B)$ {\em the $n$-th Grassmannian $B$-resolvent of} $\pi$ and $\wt{\cR}(\pi;B) = (\wt{\cR}_n(\pi;B))_{n \in \bN}$ {\em the full Grassmannnian $B$-resolvent of} $\pi$.  It is easy to check that $\wt{\cR}(\pi;B)$ {\em is a fully matricial $E$-valued analytic function on} $\wt{\rho}(\pi;B)$.

\subsection{}
\label{sec4.3} As a first step toward fitting the ``affine'' resolvents into this framework, we shall see what happens if $\pi$ is the graph of an element $Y \in E$, that is, if
\[
\pi = \left. \bpm O & 1 \\ 1 & Y \epm\right/ \wt{\l 1} \in Gr_1(E).
\]
Remark that in this case
\[
\pi \oplus \dots \oplus \pi = \left. \bpm O & I_n \otimes 1 \\ I_n \otimes 1 & I_n \otimes Y \epm\right/ \wt{\l n}.
\]
{\em We shall denote the corresponding resolvents and resolvent sets with $\pi$ replaced by $Y$ and call them Grassmannian resolvents of $Y$}, so that $\wt{\rho}_n(Y;B)$ is the $n$-th Grassmannian resolvent set of $Y$ for instance.

\bigskip
\noindent
{\bf Lemma.} {\em 
We have
\[
\left. \bpm \a & \b \\ \g & \d \epm\right/ \wt{\l n} \in \wt{\rho}_n(Y;B)
\]
iff $\d - (I_n \otimes Y)\b$ is invertible.

Moreover, then
\[
\wt{\cR}_n(Y;B)\left( \left. \bpm \a & \b \\ \g & \d \epm\right/ \wt{\l n}\right) = \b(\d - (I_n \otimes Y)\b)^{-1}.
\]
}

\bigskip
\noindent
{\bf {\em Proof.}} The iff part follows from
\[
\bpm
I_n \otimes 1 & \b \\
I_n \otimes Y & \d
\epm = \bpm
I_n \otimes 1 & 0 \\
I_n \otimes Y & \d - (I_n \otimes Y)\b
\epm \bpm
I_n \otimes 1 & \b \\
0 & I_n \otimes 1
\epm .
\]
The factorization also implies
\[
\bpm
I_n \otimes 1 & \b \\
I_n \otimes Y & \d
\epm^{-1} = \bpm
* & * \\
* & (\d - (I_n \otimes Y)\b)^{-1}
\epm .
\]
\qed

\subsection{}
\label{sec4.4} If $\s = \left. \bpm \a & I_n \otimes 1 \\ \g & \d \epm\right/\wt{\l n} \in \wt{\rho}_n(Y;B)$ then $\wt{\cR}_n(Y;B)(\s) = (\d - I_n \otimes Y)^{-1}$ and $\d \in \rho_n(Y;B)$ (the affine $n$-th resolvent).  Thus the Grassmannian resolvent set and resolvent extend the affine ones.

\subsection{Transversality}
\label{sec4.5}

We will say that $\pi = \left. \bpm a & b \\ c & d \epm \right/\wt{\l n}$ and $\pi' = \left.\bpm a' & b' \\ c' & d' \epm\right/\wt{\l n}$ in $Gr_n(B)$ are {\em transversal} if $\bpm b' & b \\ d' & d \epm$ is invertible.  It is easily seen that the relation is symmetric and well-defined (i.e., depends only on the equivalence classes/$\wt{\l n}$).  Clearly, in terms of transversality, if $\th \in Gr_1(E)$ and $\s \in Gr_n(B)$ then $\s \in \wt{\rho}_n(\th;B)$ means precisely that $\s$ and $\th \oplus \dots \oplus \th$ are transversal.  We will frequently use in this paper the following rather obvious fact.

\bigskip
\noindent
{\bf Lemma.} {\em 
Let $\pi,\pi' \in Gr_n(B)$ and $g \in GL_2(\fM_n(B))$.  Then $\pi,\pi'$ are transversal iff $g\pi$ and $g\pi'$ are transversal.
}

\bigskip
\noindent
{\bf Corollary.} {\em 
Let $\g \in GL_2(B)$, $g = I_n \otimes \g \in GL_2(\fM_n(B))$ and $\pi \in Gr_1(E)$.  Then we have $\wt{\rho}_n(\g\pi;B) = g\wt{\rho}_n(\pi;B)$.
}

\section{The derivation $\wt{\p}$ on fully matricial functions of the Grassmannian}
\label{sec5}

\subsection{}
\label{sec5.1} Let $\O$ be a fully matricial open $B$-set of the Grassmannian.  We shall denote by $\cA(\O)$ the algebra of $\bC$-valued (that is scalar) fully matricial analytic functions on $\O$, under pointwise multiplication of the matricial values.  More generally we get an algebra $\cA(\O)$ for a fully multimatricial $(B_1,\dots,B_p)$-set of the Grassmannians.  In particular, if $\O$ is a fully matricial open $B$-set of the Grassmannian, then $\O \x \O$ is a fully multimatricial $(B,B)$-set and we shall denote the corresponding algebra by $\cA(\O;\O)$.  More generally we have algebras $\cA(\O;\dots;\O)$.  This extends the construction in the affine case \cite{14}.  The aim of this section will be to extend the construction of the derivation $\p$ from the affine case to a derivation $\wt{\p}$ in the Grassmannian framework.  Like in the affine case the construction rests on two technical lemmas.

\subsection{Lemma}
\label{sec5.2}  {\em 
Let $\O = (\O_n)_{n \in \bN}$ be an open fully matricial $B$-set of the Grassmannian and let
\[
\left. \bpm
a_j & b_j \\
c_j & d_j
\epm\right/\wt{\l n_j} \in \O_{n_j}\quad (j = 1,2).
\]
Then for all $x,y,z,t \in \fM_{n_1,n_2}(B)$
\[
\left.\bpm
a_1 & x & b_1 & y \\
0 & a_2 & 0 & b_2 \\
c_1 & z & d_1 & t \\
0 & c_2 & 0 & d_2
\epm\right/\wt{\l n_1+n_2} \in \O_{n_1+n_2}.
\]
}

\bigskip
\noindent
{\bf {\em Proof.}} Since $\O$ is open, for any given $x,y,z,t$ there is $\e \ne 0$ so that the conclusion of the lemma holds with $x,y,z,t$ replaced by $\e x,\e y,\e z,\e t$.  To obtain the result without $\e$, it suffices to use the $GL(n_1+n_2;\bC)$ invariance with $s = \bpm I_{n_1} & 0 \\ 0 & \e I_{n_2} \epm$.\qed

\subsection{Lemma}
\label{sec5.3}  {\em 
Let $\O$ be an open fully matricial $B$-set of the Grassmannian and $f \in \cA(\O)$ and let $a_j,b_j,c_j,d_j,t$ be like in the preceding lemma.  Then, there is $k \in \fM_{n_1,n_2}(\bC)$ so that
\begin{align*}
&f_{n_1+n_2} \left(\left. \bpm
a_1 & 0 & b_1 & 0 \\
0 & a_2 & 0 & b_2 \\
c_1 & 0 & d_1 & t \\
0 & c_2 & 0 & d_2
\epm\right/\wt{\l n_1+n_2}\right) \\
&= \bpm
f_{n_1}\left(\left.\bpm a_1 & b_1 \\ c_1 & d_1 \epm\right/\wt{\l n_1}\right) & k \\
0 & f_{n_2} \left(\left. \bpm a_2 & b_2 \\ c_2 & d_2 \epm\right/ \wt{\l n_2}\right)
\epm
\end{align*}
and $k$ depends linearly on $t$.  In fact we have
\[
\bpm
0 & k \\
0 & 0 
\epm = \frac {d}{d\e} f_{n_1+n_2} \left.\left(\left. \bpm
a_1 & 0 & b_1 & 0 \\
0 & a_2 & 0 & b_2 \\
c_1 & 0 & d_1 & \e t \\
0 & c_2 & 0 & d_2
\epm\right/\wt{\l n_1+n_2}\right)\right|_{\e=0} .
\]
}

\bigskip
\noindent
{\bf {\em Proof.}} Assume the right-hand side of the first equality is $\bpm u & k \\ h & v \epm$.  Then by the $GL(n_1+n_2;\bC)$ equivariance of $f_{n_1+n_2}$ applied to the similarity $\bpm \e I_{n_1} & 0 \\ 0 & I_{n_2} \epm$ we find that $\bpm u & \e k \\ \e^{-1}h & v \epm$ converges as $\e \to 0$ to $\bpm
f_{n_1}\left(\left. \bpm a_1 & b_1 \\ c_1 & d_1 \epm\right/\l n_1\right) & 0 \\
0 & f_{n_2}\left(\left. \bpm a_2 & b_2 \\ c_2 & d_2 \epm\right/\l n_2 \right)
\epm$.  This, then, implies $h = 0$ and that $\bpm 0 & k \\ 0 & 0 \epm$ is given by the second formula in the statement of the lemma, since $f$ as an analytic function is differentiable.  In turn, this formula which identifies the map taking $t$ to $k$ with a partial differential of $f_{n_1+n_2}$ shows that this map is a $\bC$-linear map.\qed

\subsection{}
\label{sec5.4} To define $\wt{\p}_{n_1,n_2}f_{n_1+n_2}$, we shall use the isomorphism
\[
\a_{n_1,n_2}: \fM_{n_1} \otimes \fM_{n_2} \to \cL(\fM_{n_1,n_2})
\]
which takes $A \otimes B$ to the linear map $X \to AXB$ in $\cL(\fM_{n_1,n_2})$.

\bigskip
\noindent
{\bf Definition.} Let $\O_1,f,a_j,b_j,c_j,d_j$ be like in \ref{sec5.3} and let $T \in \cL(\fM_{n_1,n_2})$ be the linear map, so that $T(t) = k$ when $t \in \fM_{n_1,n_2}(\bC) \subset \fM_{n_1,n_2}(B)$ and
\[
\bpm 
0 & k \\
0 & 0
\epm = \frac {d}{d\e} f_{n_1+n_2} \left.\left(\left. \bpm
a_1 & 0 & b_1 & 0 \\
0 & a_2 & 0 & b_2 \\
c_1 & 0 & d_1 & \e tb_2 \\
0 & c_2 & 0 & d_2
\epm\right/\wt{\l n_1+n_2}\right)\right|_{\e=0} .
\]
Then we define
\begin{align*}
&(\wt{\p}_{n_1,n_2}f_{n_1+n_2}) \left(\left. \bpm
a_1 & b_1 \\
c_1 & d_1
\epm\right/\wt{\l n_1};\left. \bpm
a_2 & b_2 \\
c_2 & d_2
\epm\right/\wt{\l n_2}\right) \\
&= \a_{n_1,n_2}^{-1}(T) \in \fM_{n_1} \otimes \fM_{n_2}.
\end{align*}

Note that if $z_j \in GL_1(\fM_{n_j}(B))$ then $\bpm 
z_1 & 0 \\
0 & z_2
\epm \in GL_1(\fM_{n_1+n_2}(B))$ and
\begin{align*}
&f_{n_1+n_2} \left(\left. \bpm
a_1 & 0 & b_1 z_1 & 0 \\
0 & a_2 & 0 & b_2 z_2 \\
c_1 & 0 & d_1 z_1 & \e tb_2 z_2 \\
0 & c_2 & 0 & d_2 z_2
\epm\right/\l n_1+n_2\right) \\
&= f_{n_1+n_2} \left( \left. \bpm
a_1 & 0 & b_1 & 0 \\
0 & a_2 & 0 & b_2 \\
c_1 & 0 & d_1 & \e tb_2 \\
0 & c_2 & 0 & d_2
\epm\right/ \l n_1+n_2\right)
\end{align*}
so that $\wt{\p}_{n_1,n_2}f_{n_1+n_2}$ is well-defined.

It is also easy to see that $\wt{\p}$ extends the definition of $\p$ in the affine case (\cite{14}).  Indeed if we take $a_j = I_{n_j} \otimes 1$, $c_j = 0$, $b_j = I_{n_j} \otimes 1$ in the preceding formulae we get exactly the formulae in the affine case, corresponding to the embedding
\[
\fM_n(B) \ni \b \to \left.\bpm 
I_n \otimes 1 & I_n \otimes 1 \\
0 & \b
\epm\right/ \l n \in Gr_n(B).
\]

\subsection{}
\label{sec5.5} Starting with this subsection and continuing in \ref{sec5.6} and \ref{sec5.7} we will check that $\wt{\p}$ turns $A(\O)$ into a ``topological'' infinitesimal bialgebra.  Since sections~\ref{sec5.5}--\ref{sec5.7} are just a technical extension of the affine case (sections~7.7--7.10 in \cite{14}), our exposition will be more compressed.

The first step is to check that
\[
\wt{\p}f = (\wt{\p}_{m,n}f_{m+n})_{(m,n) \in \bN^2} \in \cA(\O;\O).
\]
Since analyticity of the $\wt{\p}_{m,n}f_{m+n}$ is obvious, we are left with checking $GL(m) \x GL(n)$ equivariance and the direct sum properties.

In view of the equivariance property of $\a_{m,n}$ (see 7.7 in \cite{14}) it suffices to remark that if $S' \in GL(m)$ and $S'' \in GL(n)$ then assuming $t,k \in \fM_{m,n}$ and
\[
\bpm
0 & k \\
0 & 0 
\epm = \frac {d}{d\e} f_{m+n}\left. \left( \left. \bpm
a_1 & 0 & b_1 & 0 \\
0 & a_2 & 0 & b_2 \\
c_1 & 0 & d_1 & \e tb_2 \\
0 & c_2 & 0 & d_2
\epm\right/\wt{\l m+n}\right)\right|_{\e=0}
\]
we also have
\begin{align*}
&\bpm 
0 & S'kS''^{-1} \\
0 & 0
\epm \\
&= \frac {d}{d\e} f_{m+n} \left. \left( \left. \bpm 
S'a_1S'^{-1} & 0 & S'b_1S'^{-1} & 0 \\
0 & S''a_2S''^{-1} & 0 & S''b_2S''^{-1} \\
S'c_1S'^{-1} & 0 & S'd_1S'^{-1} & \e(S'tS''^{-1})(S''b_2S''^{-1}) \\
0 & S''c_2S''^{-1} & 0 & S''d_2S''^{-1}
\epm\right/\wt{\l m+n} \right)\right|_{\e = 0}.
\end{align*}
The last equality is a consequence of the $GL(m+n)$ equivalence of $f_{m+n}$ applied to $\bpm
S' & 0 \\
0 & S''
\epm$.  We thus have proved that $\wt{\p}_{m,n}f_{m+n}$ satisfies $GL(m) \x GL(n)$ equivariance.

The direct sum properties to be checked are:  if $\pi \in \O_m$, $\s \in \O_n$ and 
\begin{align*}
m &= m' + m'',\ \pi' \in \O_{m'},\ \pi'' \in \O_{m''} \\
n &= n' + n'',\ \s' \in \O_{n'},\ \s'' \in \O_{n''}
\end{align*}
then
\begin{align*}
(\wt{\p}_{m,n}f_{m+n})(\pi' \oplus \pi'',\s) &= \wt{\p}_{m',n}f_{m'+n}(\pi',\s) \oplus \wt{\p}_{m'',n}f_{m''+n}(\pi'',\s) \\
(\wt{\p}_{m,n}f_{m+n})(\pi,\s' \oplus \s'') &= \wt{\p}_{m,n'}f_{m+n'}(\pi,\s') \oplus \wt{\p}_{m,n''}f_{m+n''}(\pi,\s'').
\end{align*}
We will only discuss the first equality to be checked, the second being obtainable along similar lines.

Since the isomorphism $\a$ has the property: 
\[
\a_{m'+m'',n}^{-1}(T_1 \oplus T_2) = \a_{m',n}^{-1}(T_1) \oplus \a_{m'',n}(T_2)
\]
if $T_1 \in \cL(\fM_{m',n})$, $T_2 \in \cL(\fM_{m'',n})$ it is easily seen that what we must prove boils down to the following.

We have
\begin{align*}
&f_{m'+m''+n} \left( \left. \bpm
a'_1 & 0 & 0 & \b'_1 & 0 & 0 \\
0 & a''_1 & 0 & 0 & b''_1 & 0 \\
0 & 0 & a_2 & 0 & 0 & b_2 \\
{} & {} \\
c'_1 & 0 & 0 & d'_1 & 0 & t' \\
0 & c''_1 & 0 & 0 & d''_1 & t'' \\
0 & 0 & c_2 & 0 & 0 & d_2
\epm\right/\wt{\l m'+m''+n} \right) \\
&= \bpm
f_{m'}\left(\left. \bpm
a'_1 & b'_1 \\
c'_1 & d'_1
\epm \right/ \wt{\l m'} \right) & 0 & k' \\
0 & f_{m''}\left(\left. \bpm
a''_1 & b''_1 \\
c''_1 & d''_1
\epm \right/\wt{\l m''} \right) & k'' \\
0 & 0 & f_n\left(\left. \bpm
a_2 & b_2 \\
c_2 & d_2
\epm \right/ \wt{\l n} \right)
\epm
\end{align*}
where
\begin{align*}
&f_{m'+n} \left(\left. \bpm
a'_1 & 0 & b'_1 & 0 \\
0 & a_2 & 0 & b_2 \\
c'_1 & 0 & d'_1 & t' \\
0 & c_2 & 0 & d_2
\epm\right/\wt{\l m'+n}\right) \\
&= \bpm
f_{m'}\left(\left.\bpm
a'_1 & b'_1 \\
c'_1 & d'_1
\epm\right/\wt{\l m'}\right) & k' \\
0 & f_n\left(\left. \bpm
a_2 & b_2 \\
c_2 & d_2
\epm\right/\wt{\l n}\right)
\epm
\end{align*}
and
\begin{align*}
&f_{m''+n}\left(\left.\bpm
a''_1 & 0 & b''_1 & 0 \\
0 & a_2 & 0 & b_2 \\
c''_1 & 0 & d''_1 & t'' \\
0 & c_2 & 0 & d_2
\epm\right/\wt{\l m''+n}\right) \\
&= \bpm
f_{m''}\left(\left.\bpm
a''_1 & b''_1 \\
c''_1 & d''_1
\epm\right/\wt{\l m''}\right) & k'' \\
0 & f_n\left(\left.\bpm
a_2 & b_2 \\
c_2 & d_2
\epm\right/\wt{\l n}\right)
\epm .
\end{align*}
If we define $k'$ and $k''$ by the last two equalities (with Lemma~\ref{sec5.3} in mind) we get
\begin{align*}
&f_{m'+m''+n}\left(\left.\bpm
a'_1 & 0 & 0 & b'_1 & 0 & 0 \\
0 & a''_1 & 0 & 0 & b''_1 & 0 \\
0 & 0 & a_2 & 0 & 0 & b_2 \\
c'_1 & 0 & 0 & d'_1 & 0 & t' \\
0 & c''_1 & 0 & 0 & d''_1 & t'' \\
0 & 0 & c_2 & 0 & 0 & d_2
\epm\right/\wt{\l m'+m''+n}\right) \\
&= \bpm
f_{m'}\left(\left.\bpm
a'_1 & b'_1 \\
c'_1 & d'_1
\epm\right/\wt{\l m'}\right) & * & * \\
0 & f_{m''}\left(\left.\bpm
a''_1 & b''_1 \\
c''_1 & d''_1
\epm\right/\wt{\l m''}\right) & k'' \\
0 & 0 & f_n\left(\left.\bpm
a_2 & b_2 \\
c_2 & d_2
\epm\right/\wt{\l n}\right)
\epm
\end{align*}
and
\begin{align*}
&f_{m''+m'+n}\left(\left.\bpm
a''_1 & 0 & 0 & b''_1 & 0 & 0 \\
0 & a'_1 & 0 & 0 & b'_1 & 0 \\
0 & 0 & a_2 & 0 & 0 & b_2 \\
c''_1 & 0 & 0 & d''_1 & 0 & t'' \\
0 & c'_1 & 0 & 0 & d'_1 & t' \\
0 & 0 & c_2 & 0 & 0 & d_2
\epm\right/\wt{\l m''+m'+n}\right) \\
&= \bpm
f_{m''}\left(\left.\bpm
a''_1 & b''_1 \\
c''_1 & d''_1
\epm\right/\wt{\l m''}\right) & * & * \\
0 & f_{m'}\left(\left.\bpm
a'_1 & b'_1 \\
c'_1 & d'_1
\epm\right/\wt{\l m'}\right) & k' \\
0 & 0 & f_n\left(\left.\bpm
a_2 & b_2 \\
c_2 & d_2
\epm\right/\wt{\l n}\right)
\epm.
\end{align*}
Using a similarity which permutes the first two summands in $\bC^{m'} \oplus \bC^{m''} \oplus \bC^n$, we get that the 13-block in the formula for $f_{m'+m''+n}(\dots)$ is $k'$.  Thus all we must still do is to show that the 12-block in that formula is zero.  This in turn is immediate from Lemma~\ref{sec5.3} applied to $f(m'+m'')+n$ and $f_{m'+m''}$.  Thus we concluded checking that
\[
(\wt{\p}_{m,n}f_{m+n})_{(m,n) \in \bN^2} \in \cA(\O;\O).
\]

\subsection{}
\label{sec5.6} Our next task is to show that $\wt{\p}: \cA(\O) \to \cA(\O;\O)$ is a derivation.

\bigskip
\noindent
{\bf Lemma.} {\em 
Let $f,g \in \cA(\O)$ and let $\pi' = \left. \bpm
a' & b' \\
c' & d'
\epm\right/\wt{\l m} \in \O_m$, $\pi'' = \left.\bpm
a'' & b'' \\
c'' & d''
\epm\right/\wt{\l n} \in \O_n$ and $t \in \fM_{m,n}$.  Then we have
\begin{align*}
\a_{m,n}((\wt{\p}_{m,n}(fg)_{m+n})(\pi';\pi''))(t) &= f_m(\pi')\a_{m,n}((\wt{\p}_{m,n}g_{m+n})(\pi';\pi''))(t) \\
&+ \a_{m,n}((\wt{\p}_{m,n}f_{m+n})(\pi';\pi''))(t)g_n(\pi'').
\end{align*}
}

\bigskip
\noindent
{\bf {\em Proof.}} To simplify notations put
\begin{align*}
\xi &= \a_{m,n}((\wt{\p}_{m,n}f_{m+n})(\pi';\pi''))(t) \in \fM_{m,n}, \\
\eta &= \a_{m,n}((\wt{\p}_{m,n}g_{m+n})(\pi';\pi''))(t) \in \fM_{m,n}, \\
\z &= \a_{m,n}((\wt{\p}_{m,n}(fg)_{m+n})(\pi';\pi''))(t) \in \fM_{m,n}
\end{align*}
and
\[
\pi = \left. \bpm
a' & 0 & b' & 0 \\
0 & a'' & 0 & b'' \\
c' & 0 & d' & tb'' \\
0 & c' & 0 & d''
\epm\right/\wt{\l m+n} \in \O_{m+n}.
\]
Then, by Lemma~\ref{sec5.3} and Definition~\ref{sec5.4} we have
\begin{align*}
(fg)_{m+n}(\pi) &= \bpm
f_m(\pi')g_m(\pi') & \z \\
0 & f_m(\pi'')g_n(\pi'')
\epm \\
f_{m+n}(\pi) & \bpm
f_m(\pi') & \xi \\
0 & f_n(\pi'')
\epm
\end{align*}
and
\[
g_{m+n}(\pi) = \bpm
g_m(\pi') & \eta \\
0 & g_n(\pi'')
\epm .
\]
The lemma then follows from the equality of matrices derived from
\[
(fg)_{m+n}(\pi) = f_{m+n}(\pi)g_{m+n}(\pi).
\]
\qed

\bigskip
\noindent
{\bf Corollary.} {\em 
$\wt{\p}: \cA(\O) \to \cA(\O;\O)$ is a derivation.
}

\bigskip
\noindent
{\bf {\em Proof.}} Take into account that if $f,g \in \cA(\O)$ and $h \in \cA(\O;\O)$ then the $\cA(\O)$-bimodule structure $\cA(\O;\O)$ is given by the homomorphisms $f \to f \otimes 1$ and $g \to 1 \otimes g$ where $(f \otimes 1)_{m,n}(\pi',\pi'') = f_m(\pi') \otimes I_n$, $(1 \otimes g)_{m,n}(\pi',\pi'') = I_m \otimes g_n(\pi'')$, and that if $A \in \fM_m$, $B \in \fM_n$, $T \in \cL(\fM_{m,n})$ then
\[
\a_{m,n}^{-1}(AT(\cdot)B) = (A \otimes I_n)\a_{m,n}^{-1}(T(\cdot))(I_m \otimes B),
\]
the Corollary is immediately inferred from the Lemma.\qed

\subsection{}
\label{sec5.7} We pass now to the proof of the co-associativity property of $\wt{\p}$.  Like in the affine case (7.10 in \cite{14}) {\em since $\cA(\O;\O)$ and $\cA(\O;\O;\O)$ have not been identified with some topological tensor products of two and respectively three copies of $\cA(\O)$, we will have to define the maps $\id \otimes \wt{\p}: \cA(\O;\O) \to \cA(\O;\O;\O)$ and $\wt{\p} \otimes \id: \cA(\O;\O) \to \cA(\O;\O;\O)$.}

Let $\wt{k} \in \cA(\O;\O)$ and put $k = \wt{k}_{m,n+p}$, which is an analytic function on $\O_m \x \O_{n+p}$ with values in $\fM_m \otimes \fM_{n+p}$.  Let further $\pi = \left. \bpm
a & b \\
c & d
\epm\right/\wt{\l m} \in \O_m$, $\pi' = \left. \bpm
a' & b' \\
c' & d'
\epm\right/\wt{\l n} \in \O_n$, $\pi'' = \left. \bpm
a'' & b'' \\
c'' & d''
\epm\right/\wt{\l p} \in \O_p$.  We define:
\begin{align*}
&((\id \otimes \wt{\p})_{m,n,p}k)(\pi;\pi';\pi'') \\
&= \sum_{\substack{1 \le a,b \le m \\ 1 \le c,d \le n \\ 1 \le e,f \le p}} \left( \left. \frac {d}{d\e} k\left(\pi;\left. \bpm 
a' & 0 & b' & 0 \\
0 & a'' & 0 & b'' \\
c' & 0 & d' & \e(e_{d,e} \otimes 1)b'' \\
0 & c'' & 0 & d''
\epm\right/\wt{\l n+p}\right)\right|_{\e=0}\right)_{(a,b)(c,n+f)} \\
&\quad e_{ab}^{(m)} \otimes e_{cd}^{(n)} \otimes e_{ef}^{(p)}
\end{align*}
where $e_{ij}^{(r)}$ are the matrix-units in $\fM_r$ and the index $(a,b)(c,n+f)$ indicates the coefficient of $e_{ab}^{(m)} \otimes e_{c,n+f}^{(n+p)}$ of an element of $\fM_m \otimes \fM_{n+p}$.  It is easy to see that if $\wt{k} = f \otimes g$, where $f,g \in \cA(\O)$ then $(\id \otimes \wt{\p})(f \otimes g) = f \otimes \wt{\p}g$.  We also leave it to the reader to check that $\id \otimes \wt{\p}$ takes values in $\cA(\O;\O;\O)$.  This involves arguments of the type used in showing that $\wt{\p}$ takes values in $\cA(\O;\O)$.

Similarly, we define
\begin{align*}
&((\wt{\p} \otimes \id)_{m,n,p}k)(\pi;\pi';\pi'') \\
&= \sum_{\substack{1 \le a,b \le m \\ 1 \le c,d \le n \\ 1 \le e,f \le p}} \left( \left. \left( \frac {d}{d\e} k \left. \bpm
a & 0 & b & 0 \\
0 & a' & 0 & b' \\
c & 0 & d & \e(e_{b,c} \otimes 1)b' \\
0 & c' & 0 & d'
\epm\right/\wt{\l m+n};\pi''\right)\right|_{\e = 0}\right)_{(a,m+d),(e,f)} \\
&\quad e_{ab}^{(m)} \otimes e_{cd}^{(n)} \otimes e_{ef}^{(p)}.
\end{align*}

Checking that $(\id \otimes \wt{\p}) \circ \wt{\p} = (\wt{\p} \otimes \id) \circ \wt{\p}$, after all these questions are put aside, boils down, like in the affine case to permuting the order in which we take two derivatives.

\bigskip
\noindent
{\bf Lemma.} {\em 
If $\wt{h} \in A(\O)$ an $h = \wt{h}_{m+n+p}$, then
\[
(\id \otimes \wt{\p})_{m,n,p} \wt{\p}_{m,n+p}h = (\wt{\p} \otimes \id)_{m,n,p} \wt{\p}_{m+n,p}h.
\]
}

\bigskip
\noindent
{\bf {\em Proof.}} Using the notations already introduced in this subsection, we have:
\begin{align*}
&((\id \otimes \wt{\p})_{m,n,p} \circ \wt{\p}_{m,n+p}h)(\pi;\pi';\pi'')_{(a,b)(c,d)(e,f)} \\
&= \frac {d}{d\e_1} \left. \left(\left. \frac {d}{d\e_1} \left( h\left( \left. \bpm
a & 0 & 0 & b & 0 & 0 \\
0 & a' & 0 & 0 & b' & 0 \\
0 & 0 & a'' & 0 & 0 & b'' \\
c & 0 & 0 & d & \e_1(e_{b,c} \otimes 1)b' & 0 \\
0 & c' & 0 & 0 & d' & \e_2(e_{d,e} \otimes 1)b'' \\
0 & 0 & c'' & 0 & 0 & d''
\epm\right/\wt{\l m+n+p}\right)\right)_{(a,m+n+f)}\right|_{\e_1=0}\right)\right|_{\e_2=0}.
\end{align*}

Similarly we have:
\begin{align*}
&((\wt{\p} \otimes \id)_{m,n,p} \circ \wt{\p}_{m+n,p}h)(\pi;\pi';\pi'')_{(a,b)(c,d)(e,f)} \\
&= \frac {d}{d\e_2} \left. \left(\left. \frac {d}{d\e_1} \left(h \left( \left. \bpm
a & 0 & 0 & b & 0 & 0 \\
0 & a' & 0 & 0 & b' & 0 \\
0 & 0 & a'' & 0 & 0 & b'' \\
c & 0 & 0 & d & \e_2(e_{b,c} \otimes 1)b' & 0 \\
0 & c' & 0 & 0 & d' & \e_1(e_{d,e} \otimes 1)b'' \\
0 & 0 & c'' & 0 & 0 & d''
\epm\right/\wt{\l m+n+p}\right)\right)_{(a,m+n+f)}\right|_{\e_1=0}\right)\right|_{\e_2=0}.
\end{align*}
Clearly the two quantities are equal (the only difference is that inside the $6 \x 6$ matrix we have replaced $\e_1$ by $\e_2$ and $\e_2$ by $\e_1$, so that the equality is just a permutability of partial derivatives).\qed

\section{The resolvent equation and the duality transform}
\label{sec6}

\subsection{}
\label{sec6.1} We shall use the same framework as in sections~\ref{sec4} and \ref{sec5}, to carry out the computations which yield the functional equation for the Grassmannian resolvent $(\wt{\cR}_n(\pi,B)(\cdot))_{n \in \bN}$ where
\[
\pi = \left. \bpm
a & b \\
c & d
\epm\right/\wt{\l 1} = Gr_1(E).
\]

Let
\[
\s' = \left. \bpm
\a' & \b' \\
\g' & \d'
\epm\right/\wt{\l m} \in \wt{\rho}_m(\pi;B)
\]
and let
\[
\s'' = \left. \bpm
\a'' & \b'' \\
\g'' & \d''
\epm\right/\wt{\l n} \in \wt{\rho}_n(\pi;B).
\]
We then consider
\[
\s = \left. \bpm
\a' & 0 & \b' & 0 \\
0 & \a'' & 0 & \b'' \\
\g' & 0 & \d' & t\b'' \\
0 & \g'' & 0 & \d''
\epm\right/\wt{\l (m+n)} \in \wt{\rho}_{m+n}(\pi;B)
\]
where $t \in \fM_{m,n}(\bC) \subset \fM_{m,n}(B)$.  To compute $\wt{\cR}_{m+n}(\pi;B)(\s)$ we must examine the matrix
\[
\Xi = \bpm
I_{m+n} \otimes b & \begin{matrix} \b' & 0 \\ 0 & \b'' \end{matrix} \\
I_{m+n} \otimes d & \begin{matrix} \d' & t\b'' \\ 0 & \d'' \end{matrix}
\epm .
\]
Permuting indices $2$ and $3$ in the above matrix, viewed as a $4 \x 4$ block-matrix, we get
\begin{align*}
\Theta &= \bpm
I_m \otimes b & \b' & 0 & 0 \\
I_m \otimes d & \d' & 0 & t\b'' \\
0 & 0 & I_n \otimes b & \b'' \\
0 & 0 & I_n \otimes d & \d''
\epm \\
&= \bpm
\begin{matrix} * & * \\ * & \z' \end{matrix} & y \\
\begin{matrix} 0 & 0 \\ 0 & 0 \end{matrix} & \begin{matrix} * & * \\ * & \z'' \end{matrix}
\epm
\end{align*}
where
\begin{align*}
y &= -\bpm * & * \\ * & \z' \epm \bpm 0 & 0 \\ 0 & t\b'' \epm \bpm * & * \\ * & \z'' \epm \\
&= \bpm
* & * \\
* & -\z' + \b''\z''
\epm
\end{align*}
and $\b'\z' = \wt{\cR}_m(\pi;B)(\s')$, $\b''\z'' = \wt{\cR}_n(\pi;B)(\s'')$.  This gives that
\[
\Theta = \bpm
* & * & * & * \\
* & \z' & * & -\z't\b''\z'' \\
0 & 0 & * & * \\
0 & 0 & * & \z''
\epm
\]
so that switching indices $2$ and $3$ we get
\[
\Xi = \bpm
* & * & * & * \\
0 & * & 0 & * \\
0 & * & \z' & -\z' t\b''\z'' \\
0 & * & 0 & \z''
\epm .
\]
The last formula implies
\begin{align*}
\wt{\cR}_{m+n}(\pi;B)(\s) &= \bpm
\b' & 0 \\
0 & \b''
\epm \bpm
\z' & -\z't\b''\z'' \\
0 & \z''
\epm \\
&= \bpm
\wt{\cR}_m(\pi;B)(\s') & -\wt{\cR}_m(\pi;B)(\s')tR_n(\pi;B)(\s'') \\
0 & \wt{\cR}_n(\pi;B)(\s'')
\epm .
\end{align*}
Comparing this with the definition of $\wt{\p}_{m,n}\wt{\cR}(\pi;B)$ we find that we have proved the following result.

\bigskip
\noindent
{\bf Lemma.} {\em 
\[
(\id_E \otimes \wt{\p}_{m,n})\wt{\cR}_{m+n}(\pi;B)(\s';\s'') = -\wt{\cR}_m(\pi;B)(\s') \otimes_E \wt{\cR}_n(\pi;B)(\s'').
\]
}

\bigskip
In the statement of the Lemma $\id_E \otimes \wt{\p}_{m,n}$ refers to applying $\wt{\p}$ to a $E$-valued fully matricial analytic function.  The $\otimes_E$ among two matrices with entries in $E$ amounts to
\[
\left( \sum_{1 \le i,j \le m} c'_{ij} \otimes e_{ij}^{(m)}\right) \otimes_E \left( \sum_{1 \le k,l \le n} c''_{kl} \otimes e_{kl}^{(n)}\right) = \sum_{i,j,k,l} c'_{ij}c''_{kl} \otimes e_{ij}^{(m)} \otimes e_{kl}^{(n)}.
\]

We can write the resolvent equation also in a more compact form.

\bigskip
\noindent
{\bf Proposition.} {\em 
\[
(\id_E \otimes \wt{\p})\wt{\cR}(\pi;B) = -\wt{\cR}(\pi;B) \otimes_E \wt{\cR}(\pi;B).
\]
}

\subsection{Matrix entries of resolvents}
\label{sec6.2} 

An extension of the duality transform of \cite{14}, from the case of $Y \in E$ to the case of $\pi \in Gr_1(E)$, includes in particular also the possibility of working with ``unbounded operators $Y$'' represented by their graph and therefore the definition of the algebra $\cR\cA(Y;B)$ in \ref{sec9.1} of \cite{14}, which includes $Y$, must be replaced in our considerations here by the definition of an algebra where $Y$ does not appear.  {\em By $\cC\cR(\pi;B)$ we shall denote the set of matrix coefficients of $\{-\wt{\cR}_n(\pi;B)(\s) \mid n \in \bN, \s \in \wt{\rho}_n(\pi;B)\}$.  By $\cL\cR(\pi;B)$ we shall denote the linear span of $\cC\cR(\pi;B)$.}

\bigskip
\noindent
{\bf Lemma.} {\em 
$\cC\cR(\pi;B)$ is closed under multiplication.  In particular $\cL\cR(\pi;B)$ is a subalgebra of $E$.}

\bigskip
\noindent
{\bf {\em Proof.}} The lemma is a consequence of the computations in \ref{sec6.1}.  Indeed let $a,b$ be the $(i,j)$ and respectively the $(k,l)$ matrix-coefficient of $-\wt{\cR}_m(\pi;B)(\s')$ and $-\wt{\cR}_n(\pi;B)(\s'')$ and let $\s$ be defined like in \ref{sec6.1} with $t = e_{jk}$.  Then the computation of $-\wt{\cR}_{m+n}(\pi;B)(\s)$, we did, shows that its $(i,m+l)$-entry is exactly the $(i,l)$-entry of $(-\wt{\cR}_m(\pi;B)(\s'))e_{jk}(-\wt{\cR}_n(\pi;B)(\s''))$ which is $ab$.\qed

\subsection{The duality transform}
\label{sec6.3}

Let $E_1$ be the closure in $E$ of $\cL\cR(\pi;B)$.  We will define the duality transform associated with $\pi$ and $B$ on the topological dual $E_1^d$ of $E_1$.  In general, the bialgebra structure is only ``partially'' defined on $E_1^d$ for analysis reasons, which cannot be dealt in this generality, we will therefore often look for formulations which avoid such problems or we will introduce extra assumptions (as we did in \cite{14}).  Some important instances when these assumptions are satisfied will be shown in \S\ref{sec12}.

If $\var \in E_1^d$, we define $\cU(\var) \in \cA(\wt{\rho}(\pi;B))$ by $\cU(\var) = (\cU(\var)_n)_{n \in \bN}$ where
\[
\cU(\var)_n(\s) = (\id_{\fM_n} \otimes \var)(\wt{\cR}_n(\pi;B)(\s))
\]
for $\s \in \wt{\rho}_n(\pi;B)$.  Since $\id_{\fM_n} \otimes \var$ is $\fM_n$-linear on $\fM_n(B)$ we infer that $\cU(\var)$ is fully matricial since $\wt{\cR}(\pi;B)$ is fully matricial.  The continuity assumption on $\var$ is necessary to obtain the analyticity of $\cU(\var)$.

We also remark that $\cU(\var) = 0$ implies $\var = 0$, that is $\cU$ is injective.  Indeed, $\cU(\var) = 0$ implies $\var \mid \cL\cR(\pi;B) = 0$ and $E_1$ is the closure of $\cL\cR(\pi;B)$.

Up to now $\cL\cR(\pi;B)$ is only an algebra so we have only a coalgebra structure on the dual (modulo technical problems).  The behavior of $\cU$ with respect to this comultiplication is recorded in the next proposition.

\bigskip
\noindent
{\bf Proposition.} {\em 
If $\var \in E_1^d$, $\s' \in \wt{\rho}_m(\pi;B)$ and $\s'' \in \wt{\rho}_n(\pi;B)$, then we have
\begin{align*}
&(\id_{\fM_m} \otimes \id_{\fM_n} \otimes \var)(\wt{\cR}_m(\pi;B)(\s') \otimes_E \wt{\cR}_n(\pi;B)(\s'')) \\
&= -\wt{\p}_{m,n}(\cU(\var)_{m+n})(\s';\s'').
\end{align*}
}

\bigskip
\noindent
{\bf {\em Proof.}} The proposition is exactly what we obtain from Lemma~\ref{sec6.1} when we apply $\id_{\fM_m} \otimes \id_{\fM_n} \otimes \var$ to the equality there.\qed

\bigskip
To justify our assertion that the above proposition shows that the behavior of $\cU$ with respect to the comultiplication, note that the right-hand side is the $(m,n)$-component of $-\wt{\p}\cU(\var)$, while the left-hand side corresponds to $(\cU \otimes \cU)(\var \circ \mu)$ with $\mu$ denoting the multiplication on $\cL\cR(\pi;B)$ (see also the proof of Lemma~\ref{sec6.2}).

\subsection{}
\label{sec6.4} Further properties of the duality transform arise when there is an appropriate derivation-comultiplication on $\cL\cR(\pi;B)$.  To avoid questions such as the action of the derivation on elements of the Grassmannian, we will resort to a somewhat tautological (from the point of view of the duality transform) characterization of the derivation.

{\em We will assume there is a derivation}
\[
\p_{\pi:B}: \cL\cR(\pi;B) \to \cL\cR(\pi;B) \otimes \cL\cR(\pi;B)
\]
{\em such that}
\[
(\id_{\fM_n} \otimes \p_{\pi:B})\wt{\cR}_n(\pi;B)(\s) = \wt{\cR}_n(\pi;B)(\s) \otimes_{\fM_n} \wt{\cR}_n(\pi;B)(\s)
\]
{\em for all $n \in \bN$ and $\s \in \wt{\rho}_n(\pi;B)$.}

For the universal unitary and hermitian Grassmannian elements this will be proved in \S\ref{sec12}.

Remark that in view of Lemma~\ref{sec6.2} the linear map $\p_{\pi:B}$ is completely determined by the relation we assume.  Thus the assumption means that this unique linear map exists and that it is a derivation.  Note also that Lemma~\ref{sec6.2} similarly implies that $\p_{\pi:B}$, if it exists, is coassociative.

\bigskip
\noindent
{\bf Proposition.} {\em 
If $\var_1,\var_2,\var_3 \in E_1^d$ are such that $\var_1(a) = (\var_2 \otimes \var_3) \circ \p_{\pi:B}(a)$ if $a \in \cL\cR(\pi;B)$, then we have
\[
\cU(\var_1) = \cU(\var_2)\cU(\var_3).
\]
}

\bigskip
\noindent
{\bf {\em Proof.}} The proposition is almost obvious in view of the way we defined $\p_{\pi:B}$.\qed

\bigskip
Of course, as the reader probably already observed, the condition characterizing $\p_{\pi:B}$ replace in the Grassmannian context the conditions $\p B = 0$, $\p Y = 1 \otimes 1$ we required in the affine case (see 9.2 in \cite{14}), which corresponds to $\pi = \left.\bpm 0 & 1 \\ 1 & Y \epm\right/\wt{\l 1}$.

\subsection{The duality transform of traces}
\label{sec6.5}

In this section we return to the context of \ref{sec6.3}, that is we will not use the derivation-comultiplication of $\cL\cR(\pi;B)$.  We will record here that Proposition~9.5 of \cite{14} on transforms of traces in the affine case extends immediately to the Grassmannian setting.

\bigskip
\noindent
{\bf Proposition.} {\em 
An element $\var \in E_1^d$ satisfies the trace-condition $\var([E_1,E_1]) = 0$ if and only if
\[
\wt{\p}_{m,n}(\cU(\var))_{m+n}(\s_1;\s_2) = \e \circ \wt{\p}_{n,m}(\cU(\var))_{m+n}(\s_2;\s_1)
\]
for all $\s_1 \in \wt{\rho}_m(\pi;B)$, $\s_2 \in \wt{\rho}_n(\pi;B)$, $m \in \bN$, $n \in \bN$.  (Here $\e: \fM_m \otimes \fM_n \to \fM_n \otimes \fM_m$ permutes the two factors.)
}

\bigskip
\noindent
{\bf {\em Proof.}} The trace condition $\var([E_1,E_1])$ is equivalent to
\[
\var([\wt{\cR}_m(\pi;B)(\s_1)_{i,j},\wt{\cR}_n(\pi;B)(\s_2)_{k,l}])
\]
for all $\s_1 \in \wt{\rho}_m(\pi;B)$, $\s_2 \in \wt{\rho}_n(\pi;B)$ and indices $i,j,k,l$.  The last equality is then equivalent, by Proposition~\ref{sec6.3} to
\[
\wt{\p}_{m,n}(\cU(\var))_{m+n}(\s_1;\s_2) = \e \circ \wt{\p}_{n,m}(\cU(\var))_{m+n}(\s_2;\s_1).
\]
\qed

\section{More on the fully matricial affine space}
\label{sec7}

Roughly, a large part of this section is about the analogue of polynomials in the context of fully matricial analytic functions on the fully matricial affine space.  Besides providing a way to construct fully matricial analytic functions, this material will also underlie the series expansions in \S\ref{sec13}.

\subsection{The polynomial sub-bialgebra $\cZ(B^d)$ of $\cA(\fM(B))$}
\label{sec7.1}

Throughout \ref{sec7.1} it will suffice to assume that $B$ is a complex Banach space and $1 \in B$ is a non-zero vector (used in the definition of $\p$), there is no need for a multiplication on $B$ here.

The {\em fully matricial affine space over} $B$, that is the largest fully matricial $B$-set will be denoted $\fM(B) = (\fM_n(B))_{n \in \bN}$.

By $1\!\!1 \in \cA(\fM(B))$ we denote the unit element $1\!\!1 = (I_n \otimes 1)_{n \in \bN}$ (constant functions).  If $\var \in B^d$ (the topological dual of $B$) we define $z(\var) = (z(\var)_n)_{n \in \bN} \in \cA(\fM(B))$ by
\[
z(\var)_n\left( \bpm
b_{11} & \dots & b_{1n} \\
\vdots & & \vdots \\
b_{n1} & \dots & b_{nn}
\epm \right) = \bpm
\var(b_{11}) & \dots & \var(b_{1n}) \\
\vdots & & \vdots \\
\var(b_{n1}) & \dots & \var(b_{nn})
\epm \in \fM_n(\bC).
\]
Since $z(\var)_n$ is linear the definition of $\p$ immediately gives
\[
\p z(\var) = \var(1) 1\!\!1 \otimes 1\!\!1.
\]
We shall denote by $\cZ(B^d)$ the subalgebra of $\cA(\fM(B))$ generated by $1\!\!1$ and $\{z(\var) \mid \var \in B^d\}$.  It is easy to see that $\cZ(B^d)$ is isomorphic to the tensor-algebra $\cT(B^d)$ over the vector-space $B^d$.  Indeed, if $\var_1,\dots,\var_n$ are linearly independent in $B^d$ we can find $b_1,\dots,b_n \in B$ so that $\var_i(b_j) = \d_{ij}$.  If $P \in \bC\<X_1,\dots,X_n\>$ is a polynomial in the noncommuting indeterminates $X_1,\dots,X_n$ so that $P \ne 0$, then there is $N \in \bN$ so that we can find $N \x N$ matrices $A_k \in \fM_N(\bC)$, $1 \le k \le n$ so that $P(A_1,\dots,A_n) \ne 0$.  Then $P(z(\var_1),\dots,z(\var_n))$ evaluated at $b_1A_1 +\dots + b_nA_n \in \fM_N(B)$ is precisely $P(A_1,\dots,A_n) \in \fM_N(\bC)$.  Thus $z(\var_1),\dots,z(\var_n)$ are algebraically free.  This suffices to guarantee that the natural unital homomorphism $\cT(B^d) \to \cA(\fM(B))$ defined by the linear map $B^d \ni \var \to z(\var) \in \cA(\fM(B))$ is injective.

The fact that $\p z(\var) = \var(1)1\!\!1 \otimes 1\!\!1$ implies that $\cZ(B^d)$ is a subcoalgebra of $\cA(\fM(B))$, that is
\[
\p \cZ(B^d) \subset \cZ(B^d) \otimes \cZ(B^d).
\]
Also the structure of $\p$ on $\cZ(B^d)$ is easy to identify.  Let $1^{\perp} = \{\var \in B^d \mid \var(1) = 0\}$ and choose some element $\th \in B^d$ so that $\th(1) = 1$.  Let then $\cZ(1^{\perp}) \subset \cZ(B^d)$ be the subalgebra of $\cZ(B^d)$ generated by $\{z(\var) \mid \var \in 1^{\perp}\}$ and which is isomorphic to $\cT(1^{\perp})$.  Then clearly $\cZ(B^d)$ identifies with $(\cZ(1^{\perp}))\<z(\th)\>$ and $\cZ(1^{\perp})$ is in $\ker \p$ while $\p z(\th) = 1 \otimes 1$.  This means that the bialgebra $\cZ(B^d)$ with the structure induced from $\cA(\fM(B))$ is isomorphic to $(\cT(1^{\perp}))\<X\>$, $\p_{X:\cT(1^{\perp})})$.  Note in particular that
\[
\ker \p \cap \cZ(B^d) = \cZ(1^{\perp}).
\]

Moreover, if $B$ is a Banach space with a continuous conjugate-linear involution $b \to b^*$, then $\cZ(B^d)$ has an involution $(z(\var))^* =  z(\var^*)$ where $\var^*(b) = \overline{\var(b^*)}$.

Also, at the end of \ref{sec8.2} we will point out in a Remark an additional feature of $\cZ(B^d)$.

\subsection{Decomposable and reducible points in $\fM(B)$}
\label{sec7.2}

Like in \ref{sec7.1}, also in \ref{sec7.2}, we will only require that $B$ be a Banach space.

In view of the similarity and direct sum requirements for ``fully matricial'' objects, we are led to look at properties of points connected with these requirements.

\bigskip
\noindent
{\bf Definition.} An element $\b  \in \fM_n(B)$ is {\em decomposable} if there are $\b' \in \fM_p(B)$, $\b'' \in \fM_q(B)$ and $S \in GL(n;\bC)$, so that $n = p + q$, $p > 0$, $q > 0$ and $S\b S^{-1} = \b' \oplus \b''$.  An element $\b \in \fM_n(B)$ is {\em reducible} if there are $\b' \in \fM_p(B)$, $\b'' \in \fM_q(B)$, $\g \in \fM_{p,q}(B)$ and $S \in GL(n;\bC)$ so that
\[
S\b S^{-1} = \bpm
\b' & \g \\
0 & \b''
\epm
\]
and $p  > 0$, $q > 0$.  An element $\b \in \fM_n(B)$ is {\em approximately decomposable} (resp.\ {\em reducible}) if it is in the closure of the decomposable (resp.\ reducible) elements.  Elements which are not decomposable (reducible, approximately decomposable, approximately reducible) will be called indecomposable (resp.\ irreducible, strongly indecomposable, strongly irreducible).

\subsection{}
\label{sec7.3} To conclude this section of remarks about the fully matricial affine space, we should point out that there is a fully matricial action of the additive group $B$ on $\fM(B)$.  For each $b \in B$ there are fully matricial maps $T(b) = (T(b)_n)_{n \in \bN}: \fM(B) \to \fM(B)$ where $T(b)_n(\b) = \b + b \otimes I_n$ which give an action of $B$ on $\fM(B)$.

In case $B$ is a Banach algebra, there is also a multiplication action given by fully matricial maps $L(b) = (L(b)_n)_{n \in \bN}$, $R(b) = (R(b)_n)_{n \in \bN}$ so that $L(b)\b = (b \otimes I_n)\b$ and $R(b)\b = \b(b\otimes I_n)$.

Note also that even if $B$ is only a Banach space there is a multiplicative action of $\bC$ on $\fM(B)$.

\section{More on the fully matricial $B$-Grassmannian and on $\wt{\p}$}
\label{sec8}

In this section we present further properties of the fully matricial $B$-Grassmannian $Gr(B) = (Gr_n(B))_{n \in \bN}$.  This includes the action by fully matricial automorphisms of $GL(2;B)$ on $Gr(B)$ and the existence of a coderivation $\L$ such that $\L - \id$ plays the role of a grading of the bialgebras $\cA(\O)$.  We also discuss the properties of $\L$ in connection with the duality transform.

\subsection{The $GL(2;B)$ action on $Gr(B)$}
\label{sec8.1}

We recall that in \ref{sec3.2} we defined $g\pi$ if $g \in GL_2(\fM_n(B))$ and $\pi \in Gr_n(B)$ and $s \cdot \pi$ if $s \in GL(n;\bC)$.  An element $h = \bpm
b_{11} & b_{12} \\
b_{21} & b_{22}
\epm \in GL(2;B)$ gives rise to elements $h_n \in GL_2(\fM_n(B))$ where $h_n = \bpm
I_n \otimes b_{11} & I_n \otimes b_{12} \\
I_n \otimes b_{21} & I_n \otimes b_{22}
\epm$.  We define $C(h): Gr(B) \to Gr(B)$ by mapping $\pi_n \in Gr_n(B)$ to $h_n\pi_n$.  It is easy to check that $h_{m+n}(\pi_m \oplus \pi_n) = (h_m\pi_m) \oplus (h_n\pi_n)$ and that $h_n(s \cdot \pi_n) = s \cdot (h_n\pi_n)$ if $s \in GL(n;\bC)$.  This {\em establishes that $C(h)$ is a fully matricial map of $Gr(B)$ into $Gr(B)$.  It is immediate from the definition that $C(\cdot)$ is an action of $GL(2;B)$ by fully matricial automorphisms of $Gr(B)$}.

It is easily seen that $C(h)$ preserves transversality in each $Gr_n(B)$.

Clearly, when $B = \bC$ the $GL(2;\bC)$-action on $Gr_1(\bC)$ is the usual action on the Riemann sphere by fractional linear transformations.

\subsection{The coderivation $\L$}
\label{sec8.2}

Let $f = (f_n)_{n \in \bN} \in \cA(\O)$, where $\O = (\O_n)_{n \in \bN}$ is a fully matricial open $B$-set of the Grassmannian.  We define $\L f = (\L_nf_n)_{n \in \bN} \in \cA(\O)$, by
\[
\L f = \frac {d}{dt} \left.\left(e^tf \circ C\left( \bpm
1 & 0 \\
0 & e^t
\epm\right)\right)\right|_{t=0}
\]
which, componentwise, amounts to
\[
(\L_nf_n)(\pi_n) = \frac {d}{dt} \left.\left(e^tf\left( \bpm
1 & 0 \\
0 & e^t
\epm_n \pi_n\right)\right)\right|_{t=0}.
\]
Since
\[
\L f - f = \frac {d}{dt} \left.\left( f \circ C \left( \bpm
1 & 0\\
0 & e^t
\epm\right)\right)\right|_{t = 0}
\]
{\em it follows that $\L-\id$ is a} derivation of $\cA(\O)$.

To prove that $\L$ is a  coderivation amounts to proving that
\[
\wt{\p} \circ \L = (\L \otimes \id + \id \otimes \L) \circ \wt{\p}.
\]
This will be a consequence of the following lemma.

\bigskip
\noindent
{\bf Lemma.} {\em 
We have
\[
\wt{\p}\left( f \circ C\left( \bpm
1 & 0 \\
0 & e^t
\epm\right)\right) = e^t(\wt{\p} f) \circ \left( C \left( \bpm
1 & 0 \\
0 & e^t
\epm\right) \x C\left( \bpm
1 & 0 \\
0 & e^t
\epm\right)\right).
\]
}

\bigskip
\noindent
{\bf {\em Proof.}} Let $\pi_m = \left. \bpm
a_1 & b_1 \\
c_1 & d_1
\epm\right/\wt{\l m}$, $\pi_n = \left.\bpm
a_2 & b_2 \\
c_2 & d_2
\epm\right/\wt{\l n}$ and let $T$ and $T'$ be defined by
\begin{align*}
T &= \a_{m,n}\left( \wt{\p}_{m,n}\left( f \circ C \left( \bpm
1 & 0 \\
0 & e^t
\epm\right)\right)(\pi_m,\pi_n)\right), \\
T' &= \a_{m,n}\left( \wt{\p}_{m,n}(f)\left( \bpm
1 & 0 \\
0 & e^t
\epm_m \pi_m,\bpm
1 & 0 \\
0 & e^t
\epm \pi_n\right)\right).
\end{align*}
Since $\a_{m,n}$ is an isomorphism, it will suffice to prove that $T(s) = T'(e^ts)$ for all $s \in \fM_{m,n}(\bC)$. Indeed, we have:
\begin{align*}
&f_{m+n}\left( \bpm
1 & 0 \\
0 & e^t
\epm_{m+n} \left.\bpm
a_1 & 0 & b_1 & 0 \\
0 & a_2 & 0 & b_2 \\
c_1 & 0 & d_1 & sb_2 \\
0 & c_2 & 0 & d_2
\epm\right/\wt{\l m+n}\right) \\
&= f_{m+n} \left( \left.\bpm
a_1 & 0 & b_1 & 0 \\
0 & a_2 & 0 & b_2 \\
e^tc_1 & 0 & e^td_1 & e^tsb_2 \\
0 & e^tc_2 & 0 & e^td_2
\epm\right/\wt{\l m+n}\right) \\
&= \bpm
f_m\left( \bpm
1 & 0 \\
0 & e^t 
\epm_m \left. \bpm
a_1 & b_1 \\
c_1 & d_1
\epm\right/\wt{\l m}\right) & T(s) \\
0 & f_n\left( \bpm
1 & 0 \\
0 & e^t
\epm_n \left. \bpm
a_2 & b_2 \\
c_2 & d_2
\epm\right/\wt{\l n}\right)
\epm \\
&= \bpm
f_m\left( \left. \bpm
a_1 & b_1 \\
e^tc_1 & e^td_2
\epm\right/\wt{\l m}\right) & T'(e^ts) \\
0 & f_n\left( \left. \bpm
a_2 & b_2 \\
e^tc_2 & e^td_2
\epm\right/\wt{\l n}\right)
\epm,
\end{align*}
which implies $T(s) = T'(e^ts)$.\qed

\bigskip
To conclude the proof of the fact that $\L$ is a coderivation it will suffice to remark that taking the derivative $\frac {d}{dt}$ at $t = 0$ of the equality in the preceding lemma gives
\[
\wt{\p}(\L f-f) = \wt{\p} f + ((\L - \id) \otimes \id + \id \otimes (\L - \id))\wt{\p} f,
\]
which immediately implies
\[
\wt{\p} \L f = (\L \otimes \id + \id \otimes \L)\wt{\p} f.
\]

\bigskip
\noindent
{\bf Proposition.} {\em 
$\L-\id$ is a derivation of $\cA(\O)$ and $\L$ is also a coderivation, that is $\wt{\p} \circ \L = (\L \otimes \id + \id \otimes \L) \circ \wt{\p}$.
}

\subsection{The derivation $D$ of $\cL\cR(\pi;B)$}
\label{sec8.3}

In the next section we will show that the coderivation $\L$ discussed in the previous section is natural from the point of view of the duality transform.  This will involve describing what the natural coderivation $L$ on $\cL\cR(\pi;B)$ should be so that for the duality described in Theorem~$5.3$ of \cite{13}, the dual coderivation corresponds under the duality transform to $\L$.  Since in \ref{sec6.4} we assumed the existence of a derivation-comultiplication $\p_{\pi:B}$ on $\cL\cR(\pi;B)$, we will handle $L$ similarly based on an additional assumption.

\bigskip
\noindent
{\bf Remark.} In the affine case of $\fM(B)$, we have $\L\cZ(B^d) \subset \cZ(B^d)$ and
\[
\L(z(\var_1) \dots z(\var_n)) = (n+1)z(\var_1)\dots z(\var_n).
\]

\bigskip
Like in \ref{sec6.2} we let $\pi \in Gr_1(E)$ and we consider $\cL\cR(\pi;B)$.  The assumption about $L$ is roughly that on $\cL\cR(\pi;B)$ there is a linear map $D$ corresponding to the infinitesimal deformation of $\pi$ into $\bpm
1 & 0 \\
0 & e^t
\epm_1 \pi$ with $t \to 0$.  We will show that $D$ must then be a derivation of $\cL\cR(\pi;B)$ with values in itself.  (In case $\pi$ is the graph of an element $Y \in E$, the deformation is $Y \to e^t Y$ with $t \to 0$.)

More precisely our assumption can be formulated as follows:  {\em We assume there is a linear map $D: \cL\cR(\pi;B) \to E$ so that
\[
(\id_{\fM_n} \otimes D)\wt{\cR}_n(\pi;B)(\s) = \frac {d}{dt} \wt{\cR}_n\left.\left( \bpm
1 & 0 \\
0 & e^t
\epm_1 \pi;B\right)(\s)\right|_{t=0}
\]
for all $\s \in \wt{\rho}_n(\pi;B)$, $n \in \bN$.}  There is a simple identity which we will use to show that $D$ takes values in $\cL\cR(\pi;B)$.

\bigskip
\noindent
{\bf Lemma.} {\em 
We have $\s \in \wt{\rho}_n \left( \bpm
1 & 0 \\
0 & e^t
\epm_1 \pi;B\right)$ iff $\bpm
1 & 0 \\
0 & e^{-t}
\epm_n \s \in \wt{\rho}_n(\pi;B)$.  Moreover, then
\[
\wt{\cR}_n\left( \bpm
1 & 0 \\
0 & e^t
\epm_1 \pi;B\right)(\s) = e^{-t}\wt{\cR}_n(\pi;B)\left( \bpm
1 & 0 \\
0 & e^{-t}
\epm_n \s\right).
\]
}

\bigskip
\noindent
{\bf {\em Proof.}} Let $\pi = \left. \bpm
a & b \\
c & d
\epm\right/ \wt{\l 1}$, $\s = \left. \bpm
\a & \b \\
\g & \d
\epm\right/\wt{\l n}$ and let also $\pi' = \left. \bpm
I_n \otimes a & I_n \otimes b \\
I_n \otimes c & I_n \otimes d
\epm\right/\wt{\l n}$ and $a' = I_n \otimes a$, $b' = I_n \otimes b$, $c' = I_n \otimes c$, $d' = I_n \otimes d$.  Then $\s \in \wt{\rho}_n\left( \bpm
1 & 0 \\
0 & e^t
\epm_1 \pi;B\right)$ means $\bpm
b' & \b \\
e^td' & \d
\epm$ is invertible and this is obviously equivalent to $\bpm
b' & \b \\
d' & e^{-t}\d
\epm$ being invertible, which is that $\bpm
1 & 0 \\
0 & e^{-t}
\epm_n \s \in \wt{\rho}_n(\pi;B)$.

Moreover if $\bpm
b' & \b \\
e^td' & \d
\epm^{-1} = \bpm
* & * \\
* & \z
\epm$ and $\bpm
b' & \b \\
d' & e^t\d
\epm^{-1} = \bpm
* & * \\
* & \xi
\epm$, then $\z = e^{-t}\xi$.  The last part of the lemma follows from the two Grassmannian resolvents being equal to $\b\z$ and $\b\xi$ respectively.\qed

\bigskip
With the notations used in the proof of the lemma, to show that $D(\cL\cR(\pi;B) \subset \cL\cR(\pi;B))$ we must prove in view of the definition of $D$ that the entries of $\frac {d}{dt}(\b\z(t))|_{t=0}$ are in $\cL\cR(\pi;B)$ or equivalently the entries of $\frac {d}{dt}(e^{-t}\b\xi(t))|_{t=0}$.  Since $\b\xi(0) = \b\z(0)$ is a resolvent, its entries are in $\cL\cR(\pi;B)$, so we are left with showing $\b\left.\left( \frac {d}{dt} \xi(t)\right)\right|_{t=0}$ has entries in $\cL\cR(\pi;B)$.  We have
\[
\bpm
* & * \\
* & \frac {d}{dt} \xi(t)
\epm = \frac {d}{dt} \bpm
b' & \b \\
d' & e^{-t}\d
\epm^{-1} = \bpm
* & * \\
* & \xi
\epm \bpm 
0 & 0 \\
0 & e^{-t}\d
\epm \bpm
* & * \\
* & \xi
\epm .
\]
Hence we infer that
\[
\left.\left( \frac {d}{dt} \xi(t)\right)\right|_{t=0} = \xi(0)\d\xi(0)
\]
and we must show that $\b\xi(0)\d\xi(0)$ has entries in $\cL\cR(\pi;B)$.  It is easily seen that the $(2,4)$-block entry of the $4 \x 4$ block matrix
\[
\G^{-1} = \bpm
b' & \b & 0 & 0 \\
d' & \d & 0 & \d \\
0 & 0 & b' & \b \\
0 & 0 & d' & \d
\epm^{-1}
\]
is precisely $\xi(0)\d\xi(0)$.  On the other hand if $S$ is the permutation matrix
\[
S = \bpm
I_n & 0 & 0 & 0 \\
0 & 0 & I_n & 0 \\
0 & I_n & 0 & 0 \\
0 & 0 & 0 & I_n
\epm
\]
we see that
\[
S\G^{-1}S^{-1} = \bpm
b' & 0 & \b & 0 \\
0 & b' & 0 & \b \\
d' & 0 & \d & \d \\
0 & d' & 0 & \d
\epm^{-1} = \bpm
* & * \\
* & Z
\epm
\]
where $\bpm
\b & 0 \\
0 & \b
\epm Z$ is an $\wt{\cR}_{2n}(\pi;B)(\mu)$ for some $\mu = \wt{\rho}_{2n}(\pi;B)$.  Hence the entries of $\bpm
\b & 0 \\
0 & \b
\epm Z$ are in $\cL\cR(\pi;B)$.  Returning to $\G^{-1}$ we see that the $(2,4)$-block entry of $\G^{-1}$ coincides with the $(3,4)$-block entry of $S\G^{-1}S^{-1}$ which is the $(1,2)$-block entry of $Z$ (the blocks are $n \x n$).  This concludes the proof that $D$ maps $\cL\cR(\pi;B)$ into itself.

To prove that $D$ is a derivation we return to the proof of Lemma~\ref{sec6.2} where we showed $\cC\cR(\pi;B)$ is closed under multiplication.  With the notation of Lemma~\ref{sec6.2} we have
\[
-(\wt{\cR}_{m+n}(\pi;B)(\s))_{i,m+l} = (\wt{\cR}_m(\pi;B)(\s'))_{ij}(\wt{\cR}_n(\pi;B)(\s''))_{kl}
\]
where $\s',\s'',(i,j),(k,l)$ were given.  Then the definition of $D$ applied to the above equality shows that $D$ is a derivation.  Concluding we have proved

\bigskip
\noindent
{\bf Proposition.} {\em 
Under our assumptions $D$ is a derivation of $\cL\cR(\pi;B)$ into itself.
}

\subsection{The coderivation $L$ of $\cL\cR(\pi;B)$}
\label{sec8.4}

In this section {\em we assume the existence of $\p_{\pi:B}$ with the properties outlined in {\em \ref{sec6.4}} and we also assume the existence of the linear map $D$ like in {\em \ref{sec8.3}} and which implies that $D$ is a derivation of $\cL\cR(\pi;B)$.  In addition, we will assume that $\p_{\pi:B}$ is closed as an operator on $\cL\cR(\pi;B)$ endowed with the norm from $E$.}

{\em We define
\[
L: D + \id: \cL\cR(\pi;B) \to \cL\cR(\pi;B).
\]
Clearly $L - \id$ is a derivation of $\cL\cR(\pi;B)$.}

\bigskip
\noindent
{\bf Lemma.} {\em 
The map $L$ is a coderivation of $(\cL\cR(\pi;B),\p_{\pi:B})$, that is
\[
\p_{\pi:B} \circ L = (\id \otimes L + L \otimes \id) \circ \p_{\pi:B}.
\]
}

\bigskip
\noindent
{\bf {\em Proof.}} Since $\cL\cR(\pi;B)$ is the linear span of $\cC\cR(\pi;B)$ it suffices to check that the equality to be proved holds for the entries of $\wt{\cR}_n(\pi;B)(\s)$.  In view of the definitions of $\p_{\pi:B}$ and $L$ this boils down to showing that
\begin{align*}
&(\id_{\fM_n} \otimes \p_{\pi:B})\left.\left( \frac {d}{dt} \wt{\cR}_n(\pi;B)\left( \bpm 
1 & 0 \\
0 & e^{-t}
\epm_n \s\right)\right)\right|_{t=0} \\
&= \frac {d}{dt} \wt{\cR}_n(\pi;B)\left.\left( \bpm 
1 & 0 \\
0 & e^{-t}
\epm_n \s\right)\right|_{t=0} \otimes_{\fM_n} \wt{\cR}_n(\pi;B)(\s) \\
&+ \wt{\cR}_n(\pi;B)(\s) \otimes_{\fM_n} \frac {d}{dt} \wt{\cR}_n(\pi;B)\left.\left( \bpm
1 & 0 \\
0 & e^{-t}
\epm_n \s\right)\right|_{t=0}.
\end{align*}
It is immediate that the right-hand side equals
\begin{align*}
&\frac {d}{dt} \left. \left( \wt{\cR}_n(\pi;B)\left( \bpm
1 & 0 \\
0 & e^{-t}
\epm_n \s\right) \otimes_{\fM_n} \wt{\cR}_n(\pi;B)\left( \bpm
1 & 0 \\
0 & e^{-t}
\epm_n \s\right)\right)\right|_{t=0} \\
&= \frac {d}{dt} \left. (\id_{\fM_n} \otimes \p_{\pi:B})\wt{\cR}_n(\pi;B)\left( \bpm
1 & 0 \\
0 & e^{-t}
\epm_n \s\right)\right|_{t=0}.
\end{align*}
Thus the equality to be proved reduces to showing that
\begin{align*}
&\frac {d}{dt} (\id_{\fM_n} \otimes \p_{\pi:B})\wt{\cR}_n(\pi;B)\left.\left( \bpm 
1 & 0 \\
0 & e^{-t}
\epm_n \s\right)\right|_{t=0} \\
&= (\id_{\fM_n} \otimes \p_{\pi:B})\left(\frac {d}{dt} \wt{\cR}_n(\pi;B)\left.\left( \bpm
1 & 0 \\
0 & e^{-t}
\epm_n \s\right)\right|_{t=0}\right).
\end{align*}
Clearly, the last equality is a consequence of the assumption that $\p_{\pi:B}$ is closed.\qed

\subsection{The coderivations $L$ and $\L$ and the duality transform}
\label{sec8.5}

In this section {\em the same assumptions as in {\em \ref{sec8.4}} will hold throughout.}

Let $E_1$ be the closure of $\cL\cR(\pi;B)$ in $E$ and let $\var \in E_1^d$ {\em so that $\var$ is in the domain of $L^d$, that is $\var \circ L$ defined on $\cL\cR(\pi;B)$ is bounded} (extends to an element of $E_1^d$).

Recall that the $n$-th component of the duality transform is defined by
\[
\cU(\var)_n(\s) = (\id_{\fM_n} \otimes \var)(\wt{\cR}_n(\pi;B)(\s)).
\]

We have
\begin{align*}
\cU(L^d\var)_n(\s) &= (\id_{\fM_n} \otimes \var)(\id_{\fM_n} \otimes L)(\wt{\cR}_n(\pi;B)(\s)) \\
&= (\id_{\fM_n} \otimes \var) \frac {d}{dt} \wt{\cR}_n(\pi;B)\left.\left( \bpm
1 & 0 \\
0 & e^{-t} 
\epm_n(\s)\right)\right|_{t=0} \\
&= \frac {d}{dt} (\id_{\fM_n} \otimes \var)\left.\left(\wt{\cR}_n(\pi;B)\left(\bpm
1 & 0 \\
0 & e^{-t}
\epm_n(\s)\right)\right)\right|_{t=0} \\
&= -\frac {d}{dt} \cU(\var)_n\left.\left( \bpm
1 & 0 \\
0 & e^t
\epm_n \s\right)\right|_{t=0} \\
&= -((\L - \id)\cU(\var)_n)(\s) = (\id - \L)\cU(\var)_n(\s).
\end{align*}
Thus we have proved the following proposition.

\bigskip
\noindent
{\bf Proposition.} {\em 
We have $\cU(L^d\var) = (\id - \L)\cU(\var)$.
}

\bigskip
Note that the way the coderivation should be transformed under duality given in \cite[Thm.~5.3]{13} is in agreement with the above Proposition.

\section{The Grassmannian involution}
\label{sec9}

Throughout this section $B$ will be a unital Banach algebra with involution.  We will discuss the corresponding involutions on  $Gr(B)$ and bialgebras $\cA(\O)$, and the properties of the duality transform related to the involutions.

\subsection{The involution on $Gr(B)$}
\label{sec9.1}

On the affine fully matricial space the involution amounts simply to the conjugate-linear antiautomorphism $T \to T^*$ on $\fM_n(B)$, $n \in \bN$.  The extension to an antiholomorphic automorphism of the fully matricial $B$-Grassmannian has some additional technical points.

{\em We will first define the orthogonal $\pi^1$ of $\pi \in Gr_n(B)$ and then we shall define $\pi^* = \bpm
0 & 1 \\
-1 & 0
\epm_n \pi^{\perp}$.}

{\em If $\pi = \left. \bpm
a & b \\
c & d
\epm\right/\wt{\l n}$ we define $\pi^{\perp} = \left. \bpm
z^* & x^* \\
t^* & y^*
\epm\right/\wt{\l n}$ where $\bpm
x & y \\
z & t
\epm = \bpm
a & b \\
c & d
\epm^{-1}$.}  To check that $\pi^{\perp}$ is well-defined we begin with a simple algebraic lemma.

\bigskip
\noindent
{\bf Lemma.} {\em 
Suppose that $\bpm
a & b \\
c & d
\epm^{-1} = \bpm
x & y \\
z & t
\epm$ and $\bpm
a' & b \\
c' & d
\epm^{-1} = \bpm
x' & y' \\
z' & t'
\epm$.  Then $\bpm
x & y \\
z' & t'
\epm$ is invertible and there is $w$ invertible so that $wx = x'$, $wy = y'$.
}

\bigskip
\noindent
{\bf {\em Proof.}} Since
\[
\bpm
x & y \\
z' & t'
\epm \bpm
a & b \\
c & d
\epm = \bpm
1 & 0 \\
* & 1
\epm
\]
is invertible, we infer $\bpm
x & y \\
z' & t'
\epm$ is invertible.

On the other hand
\[
\bpm
x & y \\
z' & t'
\epm \bpm
a' & b' \\
c' & d'
\epm = \bpm
w & 0 \\
0 & 1
\epm
\]
is invertible, so that $w$ is invertible and we have
\[
\bpm
x & y \\
z' & t'
\epm = \bpm
w & 0 \\
0 & 1
\epm \bpm
x' & y' \\
z' & t'
\epm
\]
which gives $x = wx'$, $y = wy'$.\qed

\bigskip
\noindent
{\bf Corollary.} {\em 
The map $\pi \to \pi^{\perp}$ is well-defined.
}

\bigskip
\noindent
{\bf {\em Proof.}} We have two things to check.

First, using the same notation as in the lemma, since
\[
\pi = \left.\bpm
a & b \\
c & d
\epm\right/\wt{\l n} = \left.\bpm
a' & b \\
c' & d
\epm\right/\wt{\l n}
\]
we must show that
\[
\left.\bpm
z^* & x^* \\
t^* & y^*
\epm\right/\wt{\l n} = \left.\bpm
z'^* & x'^* \\
t'^* & y'^*
\epm\right/\wt{\l n}.
\]
This is indeed so, since $w^*$ is invertible and $x^*w^* = x'^*$, $y^*w^* = y'^*$.

Secondly if $u$ is invertible and
\[
\bpm
x'' & y'' \\
z'' & t''
\epm = \bpm
a & bu \\
c & du
\epm^{-1}
\]
then it is easily seen that $x = x''$, $y = y''$ and hence clearly
\[
\left.\bpm
z''^* & x''^* \\
t''^* & y''^*
\epm\right/\wt{\l n} = \left.\bpm
z^* & x^* \\
t^* & y^*
\epm\right/\wt{\l n}.
\]
\qed

\bigskip
Remark also that the definition of $\pi^{\perp}$ can also be written
\[
\pi^{\perp} = \left.\left( \bpm
a & b \\
c & d
\epm^{*-1}\bpm
0 & 1 \\
1 & 0
\epm\right)\right/\wt{\l n} = \left.\left( \bpm
a & b \\
c & d
\epm \bpm
0 & 1 \\
1 & 0
\epm\right)^{*-1}\right/\wt{\l n}.
\]

\bigskip
\noindent
{\bf Proposition.} {\em 
We have $\pi^{**} = \pi$ and $\pi^{\perp\perp} = \pi$.  The maps $\pi \to \pi^*$ and $\pi \to \pi^{\perp}$ are antiholomorphic automorphisms of $Gr_n(B)$.
}

\bigskip
\noindent
{\bf {\em Proof.}} That $\pi^{**} = \pi$ and $\pi^{\perp\perp} = \pi$ follows immediately using the formula for $\pi^*$ and the fact that $\pi^* = \bpm
0 & 1 \\
-1 & 0
\epm_n \pi^{\perp}$.  Because of the definition of $\pi^*$ via $\pi^{\perp}$, it suffices to show that $\pi \to \pi^{\perp}$ is antiholomorphic.

The antiholomorphicity needs only to be checked in charts
\[
\left\{\left.\left.\left( \bpm
a & b \\
c & d
\epm\bpm
1 & f \\
0 & 1
\epm\right)\right/\wt{\l n}\right| f \in \fM_n(B)\right\}.
\]
If $\pi = \left. \left( \bpm
a & b \\
c & d
\epm \bpm
1 & f \\
0 & 1
\epm\right)\right/\wt{\l n}$ then $\pi^{\perp} = \left.\left( \bpm
a & b \\
c & d
\epm^{*-1} \bpm
1 & 0 \\
-f^* & 1
\epm\bpm
0 & 1 \\
1 & 0
\epm\right)\right/\wt{\l n}$ which clearly is antiholomorphic as a function of $f \in \fM_n(B)$.\qed

\bigskip
That the definition of $\pi^*$ extends the definition of the involution on the affine space is easily seen.  Indeed, then $\pi = \left.\bpm
0 & 1 \\
1 & d
\epm\right/\wt{\l n}$, $\pi^{\perp} = \left.\bpm
1 & -d^* \\
0 & 1
\epm\right/\wt{\l n}$ and $\pi^* = \left. \bpm
0 & 1 \\
-1 & d^*
\epm\right/\wt{\l n}$.

We conclude this subsection remarking that in the formula for $\pi^{\perp}$ the matrix $\bpm
0 & 1 \\
1 & 0
\epm$ can be replaced by $\bpm
0 & 1 \\
-1 & 0
\epm$, since this does not affect the second column in the result.  Hence the formula for $\pi^*$ can be written also in the form
\[
\pi^* = \left.\left( -\bpm
0 & 1 \\
-1 & 0
\epm \bpm
a & b \\
c & d
\epm^{*-1}\bpm
0 & 1 \\
-1 & 0
\epm^{-1}\right)\right/\wt{\l n}.
\]

Note also that this gives $(C(g)\pi)^* = C(W g^{*-1}W^{-1})\pi^*$ where $W = \bpm
0 & 1 \\
-1 & 0 
\epm$ and $g = GL(2,B)$.

\subsection{The involution and the bialgebras $\cA(\O)$}
\label{sec9.2}

It is easy to see that $(\pi \oplus \s)^* = \pi^* \oplus \s^*$ and that $(s \cdot \pi)^* = s^{*-1} \cdot \pi^*$ where $\s \in Gr_m(B)$, $\pi \in Gr_n(B)$, $s \in GL(n;\bC)$.

It follows that if $\O = (\O_n)_{n \in \bN}$ is a fully matricial set of the $B$-Grassmannian then the same holds for $\O^* = (\O_n^*)_{n \in \bN}$, where $\O_n^* = \{\pi^* \mid \pi \in \O_n\}$.  Clearly $\O$ is open iff $\O^*$ is open.

If $f_n: \O_n \to \fM_n$, we define $f_n^*: \O_n^* \to \fM_n$ by $(f_n(\pi))^* = f_n^*(\pi^*)$, where $\pi \in \O_n$.  If $f_n$ is analytic then so is $f_n^*$ and if $f = (f_n)_{n \in \bN} \in \cA(\O)$ then $f^* = (f_n^*)_{n \in \bN} \in \cA(\O)$ and the map $f \to f^*$ is a conjugate-linear antiisomorphism.  More generally there is a conjugate-linear antiisomorphism $f \to f^*$ of $\cA(\O_1;\O_2;\dots;\O_p)$ and $\cA(\O_1^*,\O_2^*;\dots;\O_p^*)$ where
\[
f_{n_1,\dots,n_p}^*(\o_1^*;\dots;\o_p^*) = (f_{n_1,\dots,n_p}(\o_1;\dots;\o_p))^*.
\]

If $\O = \O^*$ then $\cA(\O)$ is an algebra with involution.  More generally $\cA(\O_1;\dots;\O_p)$ is an algebra with involution when $\O_j = \O_j^*$, $1 \le j \le p$.

To state the property of $\wt{\p}$ with respect to the involution we will denote by
\[
\s_{1,2}: \cA(\O;\O) \to \cA(\O;\O)
\]
the automorphism permuting the variables, that is
\[
(\s_{1,2}f)_{m,n}(\s,\pi) = \e_{m,n} \circ f_{n,m}(\pi,\s)
\]
with $\e_{m,n}: \fM_m \otimes \fM_n \to \fM_n \otimes \fM_m$ the tensorial permutation isomorphism.

\bigskip
\noindent
{\bf Proposition.} {\em 
If $f \in \cA(\O)$, then
\[
\wt{\p}f^* = \s_{1,2}(\wt{\p}f)^*.
\]
In particular if $\O = \O^*$, this is the compatibility of the involution and comultiplication of $\cA(\O)$.
}

\bigskip
\noindent
{\bf {\em Proof.}} If $L \in \cL(\fM_{m,n})$ and $L^* \in \cL(\fM_{n,m})$ is defined by $L^*(y) = (L(y^*))^*$ then
\[
\e_{m,n}((\a_{m,n}(L))^*) = \a_{n,m}(L^*).
\]
In view of this it is easily seen that the proposition will follow if we prove that:
\begin{align*}
&f_{m,n}^*\left(\left.\left( \G \bpm
a & b & 0 & 0 \\
c & d & 0 & tb' \\
0 & 0 & a' & b' \\
0 & 0 & c' & d'
\epm \G^{-1}\right)\right/\wt{\l m+n}\right) \\
&= S\left( f_{n+m}\left(\left.\left( \G' \bpm
\a' & \b' & 0 & 0 \\
\g' & \d' & 0 & t^*\b \\
0 & 0 & \a & \b \\
0 & 0 & \g & \d
\epm \G'^{-1}\right)\right/\wt{\l n+m}\right)\right)^* S^{-1}
\end{align*}
where $\G,\G',S$ are permutation matrices, the first two having the effect by conjugation of permuting second and third rows and columns in $4 \x 4$ block matrices and $S$ permuting first and second rows and columns in a $2 \x 2$ block matrix (the sizes of blocks corresponding to $m + m + n + n$, $n + n + m  + m$ and $n+m$ respectively).  The other notations used are $t \in \fM_{m,n}(\bC)$ and
\begin{align*}
\left(\left. \bpm
a & b \\
c & d
\epm\right/\wt{\l m}\right)^* &= \left. \bpm
\a & \b \\
\g & \d
\epm\right/\wt{\l m} \\
\left( \left. \bpm
a' & b' \\
c' & d'
\epm\right/\wt{\l n}\right)^* &= \left. \bpm
\a' & \b' \\
\g' & \d'
\epm\right/\wt{\l n}.
\end{align*}

Remark that the right-hand side of the equality to be proved is equal to
\[
\left( f_{n+m}\left(S \cdot \left.\left( \G' \bpm
\a' & \b' & 0 & 0 \\
\g' & \d' & 0 & t^*\b \\
0 & 0 & \a & \b \\
0 & 0 & \g & \d
\epm \G'^{-1}\right)\right/\wt{\l n+m} \right)\right)^*.
\]
Hence by the definition of $f_{m+n}^*$ it will suffice to show that
\begin{align*}
&\left(\left.\left( \G \bpm
a & b & 0 & 0 \\
c & d & 0 & tb' \\
0 & 0 & a' & b' \\
0 & 0 & c' & d'
\epm \G^{-1}\right)\right/\wt{\l m+n}\right)^* \\
&= S \cdot \left. \left( \G' \bpm
\a' & \b' & 0 & 0 \\
\g' & \d' & 0 & t^*\b \\
0 & 0 & \a & \b \\
0 & 0 & \g & \d
\epm \G'^{-1} \right)\right/\wt{\l m+n}.
\end{align*}
Writing the equality in the form $A/\wt{\l m+n} = B/\wt{\l m+n}$ the problem amounts to showing that $A^{-1}B$ is a lower triangular $2 \x 2$ block matrix with invertible diagonal blocks.  Denoting by $\Theta$ and $\Xi$ the $4 \x 4$ explicitly written matrices in $A$ and $B$ and by $\Sigma$ and $W$ the matrices $\bpm S & 0 \\ 0 & S \epm$ and $\bpm 0 & I \\ -I & 0 \epm$, we  have:
\begin{align*}
A &= -W\G\Theta^{*-1}\G^{-1}W^{-1} \\
B &= \Sigma \G'\Xi\G'^{-1}\Sigma^{-1}.
\end{align*}
Hence
\begin{align*}
A^{-1}B &= -W\G\Theta^*\G^{-1}W^{-1}\Sigma\G'\Xi\G'^{-1}\Sigma^{-1} \\
&= -W\G\Theta^*U\Xi U^{-1}\G^{-1}W^{-1}
\end{align*}
where $U = \G^{-1}W^{-1}\Sigma\G'$.  It is easily seen that
\[
U = \bpm
0 & 0 & 0 & -I_m \\
0 & 0 & I_m & 0 \\
0 & -I_m & 0 & 0 \\
I_n & 0 & 0 & 0
\epm
\]
and hence that
\[
U\Xi U^{-1} = \bpm
\d & -\g & 0 & 0 \\
-\b & \a & 0 & 0 \\
t^*\b & 0 & \d' & -\g' \\
0 & 0 & -\b' & \a'
\epm .
\]
To compute $\Xi^*U\Xi U^{-1}$ remark first that in view of the formula for $\pi^*$, we may assume
\[
\bpm
a & b \\
c & d
\epm^{-1} = \bpm
-\d^* & \b^* \\
\g^* & -\a^*
\epm
\]
as well as the primed analogue of this.  Equivalently we have
\begin{align*}
\bpm
\d & -\g \\
-\b & \a
\epm^{-1} &= -\bpm
a^* & c^* \\
b^* & d^*
\epm \\
&= -\bpm
a & b \\
c & d
\epm^*
\end{align*}
and its primed analogue.  We get
\begin{align*}
\Theta^*U\Xi U^{-1} &= \bpm
\bpm
a & b \\
c & d
\epm^* & \begin{matrix}
0 & 0 \\
0 & 0
\end{matrix} \\
\begin{matrix}
0 & 0 \\
0 & b'^*t^*
\end{matrix} & \bpm
a' & b' \\
c' & d'
\epm^*
\epm \bpm
\d & -\g & 0 & 0 \\
-\b & \a & 0 & 0 \\
t^*\b & 0 & \d' & -\g' \\
0 & 0 & -\b' & \a'
\epm \\
&= \bpm
I & 0 & 0 & 0 \\
0 & I & 0 & 0 \\
a'^*t^*\b & 0 & I & 0 \\
0 & b'^*t^*\a & 0 & I
\epm .
\end{align*}
This in turn gives
\begin{align*}
-W\G(\Theta^*U\Xi U^{-1})\G^{-1}W^{-1} &= -W \bpm
I & 0 & 0 & 0 \\
a'^*t^*\b & I & 0 & 0 \\
0 & 0 & I & 0 \\
0 & 0 & b'^*t^*\a & I
\epm W^{-1} \\
&= -\bpm
I & 0 & 0 & 0 \\
b'^*t^*\a & I & 0 & 0 \\
0 & 0 & I & 0 \\
0 & 0 & a'^*t^*\b & I
\epm
\end{align*}
which is a matrix of the desired kind.\qed

\subsection{The involution and the coderivation $\L$}
\label{sec9.3}

In this subsection we check the compatibility of $\L$ with the involution.

\bigskip
\noindent
{\bf Proposition.} {\em 
If $f \in \cA(\O)$ then we have
\[
\L f^* = (\L f)^*
\]
(the same $\L$ denotes the coderivations in $\cA(\O)$ and in $\cA(\O^*)$).
}

\bigskip
\noindent
{\bf {\em Proof.}} If $\pi = \left. \bpm
a & b \\
c & d
\epm \right/\wt{\l n} \in \O_n$, then
\begin{align*}
\left( \bpm
1 & 0 \\
0 & e^t
\epm_n \pi\right)^* &= \bpm
0 & 1 \\
-1 & 0
\epm_n \bpm
1 & 0 \\
0 & e^t
\epm_n^{*-1} \bpm
0 & 1 \\
-1 & 0
\epm^{-1} \pi^* \\
&= \bpm
e^{-t} & 0 \\
0 & 1
\epm_n \pi^* = \bpm
1 & 0 \\
0 & e^t
\epm_n \pi^*
\end{align*}
from which the proposition follows immediately using the formula for $\L$.\qed

\subsection{The involution and Grassmannian resolvents}
\label{sec9.4}

In this subsection we check the behavior of resolvents with respect to the involution.

We will need an algebraic lemma which provides explicit formulae for resolvents.

\bigskip
\noindent
{\bf Lemma.} {\em 
Assume that, with entries in some unital ring, we have $\bpm a & b \\ c & d \epm^{-1} = \bpm x & y \\ z & t \epm$ and $\bpm \a & \b \\ \g & \d \epm^{-1} = \bpm r & s \\ u & v \epm$.  Then the matrix $\bpm b & \b \\ d & \d \epm$ is invertible iff $x\b + y\d$ is invertible, which is also iff $rb + sd$ is invertible.  Moreover we then have:
\[
\bpm 
b & \b \\
d & \d
\epm^{-1} = \bpm
(rb+sd)^{-1}r & (rb+sd)^{-1}s \\
(x\b+y\d)^{-1}x & (x\b+y\d)^{-1}y
\epm .
\]
}

\bigskip
\noindent
{\bf {\em Proof.}}  Since
\[
\bpm
x & y \\
z & t
\epm \bpm
b & \b \\
d & \d
\epm = \bpm
0 & x\b + y\d \\
1 & z\b + t\d
\epm
\]
we get the ``iff $x\b + y\d$ is invertible'' part of the statement and
\begin{align*}
\bpm
b & \b \\
d & \d
\epm^{-1} &= \bpm
0 & x\b + y\d  \\
1 & z\b + t\d
\epm^{-1} \bpm
x & y \\
z & t
\epm \\
&= \bpm
* & * \\
(x\b+y\d)^{-1} & 0
\epm \bpm
x & y \\
z & t
\epm \\
&= \bpm
* & * \\
(x\b + y\d)^{-1}x & (x\b + y\d)^{-1}y
\epm .
\end{align*}
Similarly, since
\[
\bpm
r & s \\
u & v
\epm \bpm 
b & \b \\
d & \d
\epm = \bpm
rb + sd & 0 \\
ub + vd & 1
\epm
\]
we get the ``iff $rb + sd$ is invertible'' part of the statement and
\begin{align*}
\bpm
b & \b \\
d & \d
\epm^{-1} &= \bpm
rb+sd & 0 \\
ub+vd & 1
\epm^{-1} \bpm
r & s \\
u & v
\epm \\
&= \bpm
(rb+sd)^{-1} & 0 \\
* & * 
\epm \bpm
r & s \\
u & v
\epm \\
&= \bpm
(rb+sd)^{-1}r & (rb+sd)^{-1}s \\
* & *
\epm .
\end{align*}
\qed

\bigskip
The framework for resolvents will be a unital Banach algebra with involution $E$ and a Banach subalgebra with the same involution $1 \in B \subset E$.

\bigskip
\noindent
{\bf Proposition.} {\em 
Let $\pi \in Gr_1(E)$ and $\s \in Gr_n(B)$ be such that $\s \in \wt{\rho}_n(\pi;B)$.  Then $\s^* \in \wt{\rho}_n(\pi^*;B)$ and $(\wt{\cR}_n(\pi;B)(\s))^* = \wt{\cR}_n(\pi^*;B)(\s^*)$.
}

\bigskip
\noindent
{\bf {\em Proof.}} Remark that it suffices to prove the proposition when $B = E$ and $n = 1$.  Indeed, replacing $E$ by $\fM_n(E)$ we get the reduction to the case $n=1$.

Let $\pi = \left. \bpm
a & b \\
c & d
\epm\right/\wt{\l 1}$, $\s = \left. \bpm
\a & \b \\
\g & \d
\epm\right/\wt{\l 1}$ and use the notation for the inverses of the two matrices which we used in the lemma.  Then we have $\wt{\cR}_1(\pi;E)(\s) = \b(x\b + y\d)^{-1}y$ and $\s \in \wt{\rho}_1(\pi;E)$ is equivalent to the invertibility of $x\b + y\d$ by the lemma.

On the other hand we have
\begin{align*}
\pi^* &= \left. \bpm
-t^* & y^* \\
z^* & -x^*
\epm \right/\wt{\l 1} \\
\s^* &= \left. \bpm
-v^* & s^* \\
u^* & -r^*
\epm \right/\wt{\l 1}
\end{align*}
and
\[
\bpm
-t^* & y^* \\
z^* & -x^*
\epm^{-1} = \bpm
-d^* & b^* \\
c^* & -a^*
\epm .
\]
Applying again the lemma, to these new matrices, we get that $\s^* \in \wt{\rho}_1(\pi^*;E)$ is equivalent to the invertibility of $-d^*s^* - b^*r^*$ and $\wt{\cR}_1(\pi^*;E)(\s^*) = -s^*(d^*s^* + b^*r^*)^{-1}b^*$.  Since $-d^*s^* - b^*r^*$ is invertible iff $rb + sd$ is invertible, the equivalence of $\s^* \in \wt{\rho}_1(\pi^*;E)$ with $\s \in \wt{\rho}_1(\pi;E)$ is precisely the equivalence of the invertibility of $rb + sd$ and of $x\b + yd$.

To conclude the proof of the proposition we must show that
\[
(\b(x\b + y\d)^{-1}y)^* + s^*(d^*s^* + b^*r^*)^{-1}b^* = 0,
\]
or equivalently, that
\[
\b(x\b + y\d)^{-1}y + b(rb + sd)^{-1}s = 0.
\]
This is a consequence of the last assertion of the lemma, which gives that
\[
\bpm
b & \b \\
d & \d
\epm \bpm
(rb+sd)^{-1}r & (rb+sd)^{-1}s \\
(x\b+y\d)^{-1}x & (x\b+y\d)^{-1}y
\epm = \bpm
* & 0 \\
* & *
\epm .
\]
\qed

\subsection{The involution and the duality transform}
\label{sec9.5}

Like in the previous section $1 \in B \subset E$ will be Banach algebras with involution.  Since we will consider the duality transforms with respect to $\pi$, $B$ and with respect to $\pi^*$, $B$ we will use the notations $\cU_{\pi}(\cdot)$ and respectively $\cU_{\pi^*}(\cdot)$ to distinguish the two.

\bigskip
\noindent
{\bf Proposition.} {\em 
We have $(\cL\cR(\pi;B))^* = \cL\cR(\pi^*;B)$ and $(\cU_{\pi}(\var))^* = \cU_{\pi^*}(\var^*)$.
}

\bigskip
The proof is a straightforward consequence of Proposition~\ref{sec9.4} and of the definitions of $\cL\cR(\pi;B)$ and of the duality transform and will therefore be omitted.

\section{Dual Positivity}
\label{sec10}

\subsection{The Definition}
\label{sec10.1}

The Grassmannian extension of the notion of dual positivity is quite straightforward.  Here $B$ will be a unital Banach algebra with involution.

\bigskip
\noindent
{\bf Definition.} If $\O = \O^*$ an element $f \in \cA(\O)$ is {\em dual-positive} if $f = f^*$ and $\nabla_{n,n}f(\s,\s^*)$ is a positive map of $\fM_n$ into $\fM_n$ for all $\s \in \O_n$ and $n \in \bN$ ($\nabla_{m,n}f(\s',\s'')$ denotes the map $\a_{m,n}\wt{\p}_{m,n}f(\s',\s'')$).

\bigskip
Like in the affine case we have a few equivalent conditions.

\bigskip
\noindent
{\bf Proposition.} {\em 
If $\O = \O^*$ and $f \in \cA(\O)$, the following are equivalent:
\begin{itemize}
\item[{\em (i)}] $f$ is dual positive
\item[{\em (ii)}] $f = f^*$ and for any $\s^{(j)} \in \O_{n(j)}$, $1 \le j \le p$, $\bigoplus_{1 \le i,j \le p}(\nabla_{n(i),n(j)}f)(\s^{(i)},\s^{(j)*})$ is a positive linear map of $\bigoplus_{i,j} \fM_{n(i),n(j)}$, identiied with $\fM_{n(1)+\dots+n(p)}$ into itself.
\item[{\em (iii)}] $f = f^*$ and for any $\s \in \O_n$, the map $(\nabla_{n,n}f)(\s,\s^*): \fM_n \to \fM_n$ is completely positive.
\end{itemize}
}

\bigskip
The proof from the affine case \cite[Prop.~8.2]{14} immediately carries over to this more general case and will not be repeated.

\subsection{The duality transforms of the positive functionals}
\label{sec10.2}

Here $1 \in B \subset E$ will be an inclusion of unital Banach algebras with involution.  By $E_1$ we shall denote the closure of $\cL\cR(\pi;B)$ where $\pi = \pi^* \in Gr_1(E)$.  A functional $\var \in E_1^d$ is positive, denoted $\var \ge 0$, if $\var = \var^*$ and $\var(y^*y) \ge 0$ for all $y \in E_1$ (the hermiticity follows actually from the second requirement).

\bigskip
\noindent
{\bf Proposition.} {\em 
If $\var \in E_1^d$, then $\var \ge 0$ iff $-\cU(\var)$ is dual positive in $\cA(\wt{\rho}(\pi;B))$.
}

\bigskip
\noindent
{\bf {\em Proof.}} a) $\var \ge 0 \Rightarrow -\cU(\var)$ dual positive.  We shall use Proposition~\ref{sec6.3}, which implies that
\[
(\id_{\fM_n} \otimes \id_{\fM_n} \otimes \var)(\wt{\cR}(\pi;B)(\s) \otimes_E \wt{\cR}(\pi;B)(\s^*)) = -\wt{\p}_{n,n}\cU(\var)(\s,\s^*).
\]
Since $\pi = \pi^*$, we have
\[
\wt{\cR}(\pi;B)(\s) = (\wt{\cR}(\pi^*;B)(\s^*))^*.
\]
Hence, if $x_{ij} \in E_1$ are such that $\wt{\cR}(\pi;B)(\s) = \sum_{i,j} e_{ij} \otimes x_{ij}$ then
\[
-\wt{\p}_{n,n}\cU(\var)(\s;\s^*) = \sum_{1 \le i,j,k,l \le n} \var(x_{ij}x_{lk}^*)e_{ij} \otimes e_{kl}.
\]
We must check that
\[
\a_{n,n}(-\wt{\p}_{n,n}\cU(\var)(\s;\s^*))\left(\sum_{p,q} c_p{\bar c}_qe_{pq}\right) \ge 0.
\]
In view of the definition of $\a$, this is equivalent to
\[
\sum_{1 \le i,j,k,l \le n} c_j{\bar c}_k\var(x_{ij}x_{lk}^*)e_{il} \ge 0.
\]
Let $a_i = \sum_{1 \le j \le n} c_jx_{ij} \in E_1$.  Then the inequality we must prove, becomes
\[
\sum_{1 \le i,l \le n} \var(a_ia_l^*)e_{il} \ge 0
\]
or equivalently, for all $\l_1,\dots,\l_n \in \bC$
\[
\sum_{1 \le i,l \le n} \var(a_ia_l^*)\l_i{\bar \l}_l \ge 0.
\]
Putting $y = \sum_{1 \le i \le n} \l_ia_i$, we get $\var(yy^*) \ge 0$, which is indeed so.

b) $-\cU(\var)$ dual positive $\Rightarrow \var \ge 0$.  We have
\[
((-\nabla_{n,n}\cU(\var)(\s,\s^*))e_{jk})_{il} = \var((\wt{\cR}_n(\pi;B)(\s))_{ij}(\wt{\cR}_n(\pi;B)(\s))_{lk}^*).
\]
We must show that $\var(\xi\xi^*) \ge 0$ if $\xi \in \cL\cR(\pi;B)$.  Using instead of $\s',\s'',\s''',\dots$ the direct sum $\s = \s' \oplus \s'' \oplus \s''' \oplus \dots$ we can assume $\xi$ is a linear combination of matrix entries of $\wt{\cR}_n(\pi;B)(\s)$.  So, we must show that
\[
\sum_{1 \le i,j,k,l \le n} c_{ij}{\bar c}_{lk}\var((\wt{\cR}_n(\pi;B)(\s))_{ij}(\wt{\cR}_n(\pi;B)(\s))_{lk}^*) \ge 0
\]
which is equivalent to
\[
\sum_{1 \le i,j,k,l \le n} (-\nabla_{n,n} \cU(\var)(\s;\s^*))(e_{jk}))_{il}c_{ij}{\bar c}_{lk} \ge 0.
\]
Let $\Phi = -\nabla_{n,n} \cU(\var)(\s;\s^*)$.  The complete positivity of $\Phi$ gives that
\[
0 \le \Theta = \sum_{1 \le j,k \le n} \Phi(e_{jk}) \otimes e_{jk} = \sum_{1 \le i,j,k,l}(\Phi(e_{jk}))_{il}e_{il} \otimes e_{jk}.
\]
If $e_l$ is the canonical basis in $\bC^n$ and $\eta = \sum_{1 \le l,k \le n} {\bar c}_{lk}e_l \otimes e_k \in \bC^n \otimes \bC^n$, then
\[
0 \le \<\Theta \eta,\eta\> = \sum_{1 \le i,j,k,l \le n}(\Phi(e_{jk}))_{il}c_{ij}{\bar c}_{lk}
\]
which concludes the proof.\qed

\section{Stably Matricial Sets, Matricial Half-Planes and Unit Balls}
\label{sec11}

\subsection{Stably matricial sets}
\label{sec11.1}

Matricial half-planes, unit balls, which are fundamental in operator theory, don't satisfy the similarity requirement of fully matricial sets.  These matricial sets have somewhat weaker properties, and will be called {\em stably matricial}.  Stably matricial sets can easily be turned into fully matricial sets.  Since affine space is part of the Grassmannian we will only give the Grassmannian definition, though, the easy adaptation to the affine case, will often be easier to work with.

\bigskip
\noindent
{\bf Definition.} a) If $\Xi = (\Xi_n)_{n \in \bN}$, where $\Xi_n \subset Gr_n(B)$ we will call $\Xi$ a {\em stably matricial $B$-set of the Grassmannian} if the following two conditions are satisfied.
\begin{itemize}
\item[$1^{\circ}$] if $\pi \in Gr_m(B)$, $\s \in Gr_n(B)$ then $\pi \oplus \s \in \Xi_{m+n} \Leftrightarrow \pi \in \Xi_m$, $\s \in \Xi_n$.
\item[$2^{\circ}$] if $\pi \in Gr_m(B)$, $\s \in Gr_n(B)$ and $s \in GL(m+n;\bC)$ is such that $s \cdot (\pi \oplus \s) \in \Xi_{m+n}$, then there are $s' \in GL(m;\bC)$, $s'' \in GL(n;\bC)$, so that $s' \cdot \pi \in \Xi_m$, $s'' \cdot \s \in \Xi_n$.
\end{itemize}

b) If $\Xi = (\Xi_n)_{n \in \bN}$ is stably matricial and $f = (f_n)_{n \in \bN}$, where \linebreak $f_n: \Xi_n \to \fM_n$, then we will say $f$ is a stably matricial function on $\Xi$ if
\begin{itemize}
\item[$1^{\circ}$] $f_{m+n}(\pi \oplus \s) = f_m(\pi) \oplus f_n(\s)$ when $\pi \in \Xi_m$, $\s \in \Xi_n$.
\item[$2^{\circ}$] if $\pi,\pi' \in \Xi_n$ and $s \in GL_n(\bC)$ is such that $s \cdot \pi = \pi'$, then $\Ad s(f_n(\pi)) = f_n(\pi')$.
\end{itemize}

\bigskip
\noindent
{\bf Remark.} The reader can easily adapt part b) of the definition to functions with $f_n: \Xi_n \to \fM_n(\cX)$, where $\cX$ is some Banach space or even more generally $f_n: \Xi_n \to \Theta_n$, where $\Theta = (\Theta_n)_{n \in \bN}$ is a stably matricial $A$-set, where $A$ is some Banach algebra.  In another direction, there are several variables stably matricial functions.  For instance for two variables we consider $f = (f_{m,n})_{(m,n) \in \bN^2}$, $f_{m,n}: \Xi_m \x \Xi_n \to \fM_m \otimes \fM_n$ satisfying: $f_{m'+m'',n'+n''}(\pi' \oplus \pi'',\s' \oplus \s'') = f_{m',n'}(\pi',\s') \oplus f_{m'',n'}(\pi'',\s') \oplus f_{m',n''}(\pi',\s'') \oplus f_{m'',n''}(\pi'',\s'')$ and the similarity condition $(\Ad s_1 \otimes \Ad s_2)f_{m,n}(\pi,\s) = f_{m,n}(\pi',\s')$ when $\pi,\pi' \in \Xi_m$, $\s,\s' \in \Xi_n$ and $s_1 \in GL(m,\bC)$, $s_2 \in GL(n;\bC)$ are such that $s_1 \cdot \pi = \pi'$, $s_2 \cdot \s = \s'$.  The corresponding extensions of part b) of the next proposition is also an easy exercise left to the reader.

\bigskip
\noindent
{\bf Proposition.} {\em 
a)If $\Xi = (\Xi_n)_{n \in \bN}$ is a stably matricial $B$-set of the Grassmannian, then $\wt{\Xi} = (\wt{\Xi}_n)_{n \in \bN}$, where $\wt{\Xi}_n = GL(n;\bC) \cdot \Xi_n$, is a fully matricial $B$-set of the Grassmannian.  Moreover, if $\Xi$ is open, then $\wt{\Xi}$ is also open.

b)If $f = (f_n)_{n \in \bN}$ is a stably  matricial function ($f_n$ here is $\fM_n$-valued), then there is a unique extension $\wt{f} = (\wt{f}_n)_{n \in \bN}$ to $\wt{\Xi}$ so that $\wt{f}_n(s \cdot \pi) = \Ad s(\wt{f}_n(\pi))$.  The extension $\wt{f}$ is fully matricial.  Moreover if $\Xi$ is open and $f$ is analytic, then $\wt{f}$ is fully matricial analytic.  Similar statements hold for the more general $f$ considered in the remark preceding the proposition.

c) If $\Xi$ is a stably matricial $B$-set of the Grassmannian and $g \in GL(2;B)$, then $C(g)\Xi$ is also stably matricial and the map $\pi \to C(g)\pi$ defines a stably matricial isomorphism of $\Xi$ and $C(g)\Xi$, which is analytic when $\Xi$ is open.
}

\bigskip
The proof is quite straightforward checking and will be left to the reader.

\bigskip
\noindent
{\bf Corollary.} {\em 
If $\Xi,\wt{\Xi}$ are like in the proposition then the restriction maps $\cA(\wt{\Xi}) \to \cA(\Xi)$, $\cA(\wt{\Xi},\wt{\Xi}) \to \cA(\Xi,\Xi)$ (the $\cA(\Xi)$, $\cA(\Xi,\Xi)$ denoting the analytic stably matricial functions in the case of stably matricial sets) are bijective.  In particular $\cA(\Xi)$ becomes via the isomorphism a bialgebra with coderivation $\L$, isomorphic to $\cA(\wt{\Xi})$.
}

\subsection{The stably matricial unit disk $\cD_0(B)$}
\label{sec11.2}

Starting with this subsection and throughout the rest of section~\ref{sec11} we shall assume $B$ is a unital $C^*$-algebra.

Let $\cD_0(B)_n = \{T \in \fM_n(B) \mid \|T\| < 1\}$ be the open unit ball of $\fM_n(B)$ with respect to the $C^*$-norm and let $\cD_0^{cl}(B)_n$ be the closed unit ball.  We shall also denote by $\cU(n;B)$ the unitary group of $\fM_n(B)$.  Of course being subsets of $\fM_n(B)$ these sets are also subsets of $Gr_n(B)$.

\bigskip
\noindent
{\bf Proposition.} {\em 
The matricial sets $\cD_0(B) = (\cD_0(B)_n)_{n \in \bN}$, $\cD_0^{cl}(B) = (\cD_0(B)_n^{cl})_{n \in \bN}$ and $\cU(B) = (U(n;B))_{n \in \bN}$ are stably matricial $B$-sets.  Also, if $g_1,\dots,g_k \in GL(2;B)$, then more generally
\[
C(g_1)\cD_0(B) \cap \dots \cap C(g_k)\cD_0(B) = (C(g_1)\cD_0(B)_n \cap \dots \cap C(g_k)\cD_0(B)_n)_{n \in \bN}
\]
is an open stably matricial set.
}

\bigskip
\noindent
{\bf {\em Proof.}} The only part of the statement which is perhaps not immediately clear is property $2^{\circ}$ of part a) of Definition~\ref{sec11.1}.  In view of Proposition~\ref{sec11.1} part c), the statement about $C(g_1)\cD_0(B) \cap \dots \cap C(g_k)\cD_0(B)$ needs only to be checked in case $g_1 = I_2$, since we may move the set by the automorphism defined by $g_1^{-1}$, which replaces $g_1,\dots,g_k$ by $1,g_1^{-1}g_2,\dots,g_1^{-1}g_k$.  Thus, the considered sets will all be affine and it will simplify notations to work in the affine  context.  Assume $s(x \oplus y)s^{-1} \in \Xi_{m+n}$, where $s \in GL(m+n;\bC)$, $x \in \fM_m(B)$, $y \in \fM_n(B)$ and where $\Xi$ denotes one of the matricial sets we are considering.  We must show there is $s' \in GL(m;\bC)$ such that $s'xs'^{-1} \in \Xi_m$ (the statement is symmetric in $x$ and $y$, so only one half needs to be proved).  Let $V = s(\bC^m \oplus 0_n) \subset \bC^{m+n}$ and assume $B$ has been identified with a $C^*$-algebra of operators on some Hilbert space $H$.  Then $\fM_n(B)$ becomes a $C^*$-algebra of operators acting on $\bC^n \otimes H$ and $V \otimes H$ is an invariant subspace for $s(x \oplus y)s^{-1}$ and $T = s(x \oplus y)s^{-1} \mid V \otimes H$ is an operator of norm $< 1$, $\le 1$ or unitary, depending on which of our first three possible $\Xi$ we consider (for the unitary case the restriction is isometric, but being similar to $x$ which is invertible, it must be unitary).  In the fourth case we already know $\|T\| < 1$ and if $g_j^{-1} = \bpm
b_1 & b_2 \\
b_3 & b_4
\epm$ any $2 \le j \le k$, then $I_V \otimes b_1 + (I_V \otimes b_2)T$ is invertible because it is similar to $I_m \otimes b_1 + (I_m \otimes b_2)x$ which  is invertible and we have
\begin{align*}
&\|(I_V \otimes b_3 + (I_V \otimes b_4)T)(I_V \otimes b_1 + (I_V \otimes b_2)T)^{-1}\| \\
&\le \|(I_n \otimes b_3 + (I_n \otimes b_4)(s(x \oplus y)s^{-1}))(I_n \otimes b_1 + (I_n \otimes b_2)(s(x \oplus y)s^{-1}))^{-1}\| < 1.
\end{align*}
If $u: V \to \bC^m$ is a unitary operator and if $s' = u(_V|s|_{\bC^m \oplus 0_n})$, then $s'xs'^{-1}$ is unitarily equivalent to $s(x \oplus y)s^{-1}|_{V \otimes H}$ and hence in $\Xi_n$.\qed

\bigskip
An obvious characterization of $C(g)\cD_0(B)_n$ which we used in the proof and which will be useful also in what follows, is that:
\[
\left. \bpm
\a & \b \\
\g & \d
\epm\right/\wt{\l n} \in C(g)\cD_0(B)_n,
\]
where $g \in GL(2;B)$ and $g^{-1} = \bpm
b_1 & b_2 \\
b_3 & b_4
\epm$ if and only if $(b_1 \otimes I_n)\b + (b_2 \otimes I_n)\d$ is invertible and
\[
\|((b_3 \otimes I_n)\b + (b_4 \otimes I_n)\d)((b_1 \otimes I_n)\b + (b_2 \otimes I_n)\d)^{-1}\| < 1.
\]

Similarly $\left. \bpm
\a & \b \\
\g & \d
\epm\right/\wt{\l n} \in C(g)\cU(n;B)$ iff $(b_1 \otimes I_n)\b + (b_2 \otimes I_n)\d$ is invertible and
\[
((b_3 \otimes I_n)\b + (b_4 \otimes I_n)\d)((b_1 \otimes I_n)\b + (b_2 \otimes I_n)\d)^{-1} \in \cU(n;B).
\]

\subsection{The stably matricial upper and lower half-planes $(\cH_{\pm}(B))$}
\label{sec11.3}

We will define some of the matricial sets which underly the noncommutative spectral analysis of a hermitian element with respect to $B$.  This will be done via transformations $C(g)$ from $\cD_0(B)$ and $\cU(B)$.

It will be convenient to have also a notation for the open stably matricial disk at infinity $\cD_{\i}(B) = C\left( \bpm
0 & 1 \\
1 & 0
\epm\right) \cD_0(B)$.  The open upper and lower stably matricial half planes $\cH_+(B) = (\cH_+(B)_n)_{n \in \bN}$ and $\cH_-(B) = (\cH_-(B)_n)_{n \in \bN}$ are then defined by
\begin{align*}
\cH_+(B) &= C\left(\bpm
-i & i \\
1 & 1
\epm\right)\cD_0(B) \\
\cH_-(B) &= C\left(\bpm
i & -i \\
1 & 1
\epm\right)\cD_0(B)
\end{align*}
(which is just the familiar Cayley transforms in the notations we use).  Equivalently we have
\begin{align*}
\cH_+(B) &= C\left(\bpm
i & -i \\
1 & 1
\epm\right)\cD_{\i}(B) \\
\cH_-(B) &= C\left(\bpm
-i & i \\
1 & 1
\epm\right)\cD_{\i}(B)
\end{align*}
and
\[
\cD_0(B) = \bpm
i & 1 \\
-i & 1
\epm \cH_+(B) = \bpm
-i & 1 \\
i & 1
\epm \cH_-(B).
\]
Since $\cD_0(B)^* = \cD_0(B)$ we can use the formula 
\[
(C(g)\pi)^* = C(Wg^{*-1}W^{-1})\pi^*
\]
where $W = \bpm
0 & 1 \\
-1 & 0
\epm$.  This then easily gives 
\begin{align*}
\cD_{\i}(B)^* &= \cD_{\i}(B) \mbox{ and} \\
\cH_+(B)^* &= \cH_-(B).
\end{align*}
The hermitian stably matricial set $\cH(B) = (\cH(B)_n)_{n \in \bN}$ is obtained from $\cU(B)$ by ``Cayley transform''
\[
\cH(B) = C\left(\bpm
i & -i \\
1 & 1
\epm\right)\cU(B).
\]
Since $\cU(B) = \cU(B)^*$ and $\cU(B) = C\left( \bpm
0 & 1 \\
1 & 0
\epm\right)\cU(B)$, etc., it is easy to see that
\[
\cH(B) = \cH(B)^*
\]
and also
\[
\cH(B) = C\left( \bpm
-i & i \\
1 & 1
\epm\right)\cU(B).
\]

\bigskip
\noindent
{\bf Remark.} We have left out in the above discussion, because it seemed too well-known, the important fact that the affine part $\fM_n(B) \cap \cH_+(B)_n$ of $\cH_+(B)_n$ consists of those $\b \in \fM_n(B)$ such that $\Im \b \ge 0$ and $\Im \b$ is invertible, that is $\Im \b \ge \e I_n \otimes 1$ for some $\e > 0$.  To take this to the Grassmannian leads to Krein spaces (see \cite{5}).

\subsection{Mixed unit balls $\D_{p,q}(B)$ and half-planes $X_{p,q}(B)$ in the matricial resolvents of unitary and hermitian elements}
\label{sec11.4}

Let $E$ be a unital $C^*$-algebra so that $1 \in B \subset E$ and let $u \in \cU(1;E)$ and $\chi  \in \cH(E)_1$.  We will show that in $\wt{\rho}(u;B)$ and $\wt{\rho}(\chi;B)$, the subsets $\cD_0(B)_n$, $\cD_{\i}(B)_n$, $\cH_{\pm}(B)_n$ are part of some more general families of ``unit balls'' and ``half-planes''.  Of course the questions for $\chi$ and $u$ are equivalent via Cayley transform.

If $x \in \cD_0(B)_n$ then $x \in \wt{\rho}_n(u;B)$ and $I_n \otimes u - x = (I_n \otimes u)(I_n \otimes 1 - (I_n \otimes u^*)x)$ is invertible because $\|(I_n \otimes u^*)x\| < 1$.  Since $C\left( \bpm
0 & 1 \\
1 & 0  
\epm\right)u = u^{-1} \in  \cU(1;E)$ we also have $C\left( \bpm
0 & 1 \\
1 & 0
\epm\right)\cD_0(B) = \cD_{\i}(B) \subset \wt{\rho}(u;B)$.  Hence, $\wt{\rho}(u;B)$ being fully  matricial, if $x \in \cD_0(B)_p$, $y \in \cD_0(B)_q$ then $x \oplus C\left(\bpm 
0 & 1 \\
1 & 0
\epm\right)y \in \wt{\rho}_{p+q}(u;B)$.  Since $\wt{\rho}_{p+q}(u;B)$ is open, but $\cD_0(B)_p \oplus C\left(\bpm
0 & 1 \\
1 & 0
\epm\right)\cD_0(B)_q$ isn't, if $p > 0$, $q > 0$, there should be a larger open set in $\wt{\rho}_{p+q}(u;B)$ containing it.  Such a set is $\D_{p,q}(B) \subset Gr_{p+q}(B)$.

Using the action of $GL(2(p+q),\bC)$ on $Gr_{p+q}(B)$ we define
\[
\D_{p,q}(B) = g_{p,q}\cD_0(B)_{p+q}
\]
where
\[
g_{p,q} = \bpm
I_p & 0 & 0 & 0 \\
0 & 0 & 0 & I_q \\
0 & 0 & I_p & 0 \\
0 & I_q & 0 & 0
\epm \in GL(2(p+q),\bC).
\]
In particular $\D_{p,0}(B) = \cD_0(B)_p$ and $\D_{0,q}(B) = \cD_{\i}(B)_q$.

From here we also define
\begin{align*}
\wt{\D}_{p,q}(B) &= GL(p+q;\bC) \cdot \D_{p,q}(B), \\
X_{p,q}(B) &= C\left(\bpm
-i & i \\
1 & 1
\epm\right)\D_{p,q}(B), \\
\wt{X}_{p,q}(B) &= GL(p+q;\bC)\cdot X_{p,q}(B) \\
&= C\left(\bpm
-i & i \\
1 & 1
\epm\right) \wt{\D}_{p,q}(B).
\end{align*}

Since the definition of $\D_{p,q}(B)$ is tied to the special decomposition $\bC^{p+q} = \bC^p \oplus \bC^q$, it is natural to consider also sets like
\begin{align*}
\wh{\D}_{p,q}(B) &= U(p+q;\bC) \cdot \D_{p,q}(B) \\
\wh{X}_{p,q}(B) &= U(p+q;\bC) \cdot X_{p,q}(B).
\end{align*}
Clearly $X_{p,0}(B) = \cH_+(B)_p$, $\wt{X}_{p,0}(B) = \wt{\cH}_+(B)_p$, $X_{0,q}(B) = \cH_-(B)_q$, $\wt{X}_{0,q}(B) = \wt{\cH}_-(B)_q$.

\bigskip
\noindent
{\bf Lemma.} {\em 
We have $\s \in \D_{p,q}(B)$ if and only if for some
\[
\bpm
x & y \\
z & t
\epm \in \cD_0(B)_{p+q}
\]
we have
\[
\s = \left. \bpm
0 & 0 & I_p \otimes 1 & 0 \\
0 & I_q \otimes 1 & z & t \\
I_p \otimes 1 & 0 & x & y \\
0 & 0 & 0 & I_q \otimes 1
\epm\right/\wt{\l p+q} .
\]
}

\bigskip
The proof is a straightforward computation which will be omitted.

\bigskip
\noindent
{\bf Proposition.} {\em 
We have $\wt{\D}_{p,q}(B) \subset \wt{\rho}_{p+q}(u;B)$ and $\wt{X}_{p,q}(B) \subset \wt{\rho}_{p+q}(x;B)$.  Moreover $C\left(\bpm
0 & 1 \\
1 & 0
\epm\right) \wh{\D}_{p,q} = \wh{\D}_{q,p}$, $C\left( \bpm
0 & 1 \\
1 & 0
\epm\right) \wt{\D}_{p,q} = \wt{\D}_{q,p}$ and $\wh{X}_{p,q}^* = \wh{X}_{q,p}$, $\wt{X}_{p,q}^* = \wt{X}_{q,p}$.
}

\bigskip
\noindent
{\bf {\em Proof.}}  To establish the two inclusions it suffices to show $\D_{p,q}(B) \subset \wt{\rho}_{p+q}(u;B)$, since $\wt{\rho}(u;B)$ is fully matricial and we can use Cayley transform to pass from $u$ to $\chi$.  In view of the lemma, we must show that
\[
V = \bpm
I_p\otimes 1 & 0 & I_p \otimes 1 & 0 \\
z & t & 0 & I_q \otimes 1 \\
x & y & I_p \otimes u & 0 \\
0 & I_q \otimes 1 & 0 & I_q \otimes u
\epm
\]
is invertible when
\[
\bpm
x & y \\
z & t 
\epm \in \cD_0(B)_{p+q}.
\]
This is equivalent to the invertibility of
\[
g_{p,q}V \bpm
I_{p+q+p} \otimes 1 & 0 \\
0 & I_q \otimes u^{-1}
\epm = \bpm
I_p \otimes 1 & 0 & I_p \otimes 1 & 0 \\
0 & I_q \otimes 1 & 0 & I_q \otimes 1 \\
x & y & I_p \otimes u & 0 \\
z & t & 0 & I_q \otimes u^{-1}
\epm .
\]
Multiplying to the right by
\[
\bpm
I_{p+q} \otimes 1 & -I_{p+q} \otimes 1 \\
0 & I_{p+q} \otimes 1
\epm
\]
we see that the invertibility of $V$ is equivalent to the invertibility of
\begin{align*}
&-\bpm
x & y \\
z & t
\epm + \bpm
I_p \otimes u & 0 \\
0 & I_q \otimes u^{-1}
\epm \\
&= \bpm
I_p \otimes u & 0 \\ 
0 & I_q \otimes u^{-1}
\epm \left( I_{p+q} \otimes 1 - \bpm
I_p \otimes u^{-1} & 0 \\
0 & I_q \otimes u
\epm \bpm
x & y \\
z & t
\epm \right)
\end{align*}
which follows from the fact that
\[
\left\| \bpm
I_p \otimes u^{-1} & 0 \\
0 & I_q \otimes u
\epm \bpm
x & y \\
z & t
\epm \right\| < 1.
\]

Let $s_{p,q} = \bpm
0 & I_p \\
I_q & 0
\epm \in U(p+q;\bC)$.  The transformation of $\wh{\D}_{p,q}$ and $\wt{\D}_{p,q}$ into $\wh{\D}_{q,p}$ and $\wt{\D}_{q,p}$ via $C\left( \bpm
0 & 1 \\
1 & 0 
\epm\right)$ is a consequence of
\[
C \left(\bpm 
0 & 1 \\
1 & 0
\epm\right) \D_{p,q} = s_{p,q} \cdot \D_{q,p}
\]
which in turn follows from the easy to check matrix equality
\[
\bpm
0 & I_{p+q} \\
I_{p+q} & 0
\epm g_{p,q} = \bpm
s_{p,q} & 0 \\
0 & s_{p,q}
\epm g_{q,p} \bpm
s_{q,p} & 0 \\
0 & s_{q,p}
\epm .
\]

Passing to the action of the involution on $\wt{X}_{p,q}$, $\wh{X}_{p,q}$ we observe that if $\s \in Gr_n(B)$ and $\g \in GL(2n;\bC)$ then $(\g\s)^* = (W_n\g^{*-1}W_n^{-1})\s^*$ where $W_n = \bpm
0 & I_n \\
-I_n & 0
\epm$.  If $\pi \in X_{p,q}$ then there is $\s \in \cD_0(B)_{p+q}$ so that
\[
\pi = C_{p+q}g_{p,q}\s
\]
where $C_{p+q}$ stands for $\bpm
-iI_{p+q} & iI_{p+q} \\
I_{p+q} & I_{p+q}
\epm$.  It is easily seen that $g_{p,q}^{*-1} = g_{p,q}$, $C_{p+q}^{*-1} = \frac {1}{2} C_{p+q}$ and $W_{p+q}C_{p+q}W_{p+q}^{-1} = -iC_{p+q} \bpm
0 & I_{p+q} \\
I_{p+q} & 0
\epm$.  Hence
\[
\pi^* = C_{p+q}\G_{p+q}W_{p+q}g_{p,q}W_{p+q}^{-1} \s^*
\]
where $\G_n = \bpm
0 & I_n \\
I_n & 0
\epm$.  Let $d_{p,q} = \bpm
I_p & 0 \\
0 & -I_q
\epm \in U(p+q;\bC)$.  It is easy to see that $W_{p+q}g_{p,q}W_{p+q}^{-1} = g_{p,q} \bpm
d_{p,q} & 0 \\
0 & d_{p,q}
\epm$ so that
\[
\pi^* = C_{p+q}\G_{p+q}g_{p,q}(d_{p,q} \cdot \s^*).
\]
Since $d_{p,q} \s^* \in \cD_0(B)_{p+q}$ we have
\begin{align*}
\pi^* \in C\left( \bpm
-i & i \\
1 & 1
\epm\right) C\left( \bpm
0 & 1 \\
1 & 0
\epm\right) \D_{p,q}
&= C\left( \bpm
-i & i \\
1 & 1
\epm \right) s_{p,q} \cdot \D_{q,p} \\
&= s_{p,q} \cdot X_{q,p}.
\end{align*}
Thus we have proved
\[
X_{p,q}^* \subset s_{p,q}X_{q,p}.
\]
This easily implies
\begin{align*}
&\wh{X}_{p,q}^* \subset \wh{X}_{q,p} \mbox{ and} \\
&\wt{X}_{p,q}^* \subset \wt{X}_{q,p}
\end{align*}
which by symmetry must be equalities.\qed

\subsection{Extending the invertibility of stably matricial functions}
\label{sec11.5}

Inverse function theorems usually provide local inverses and this raises the question whether we can use them to get fully matricial inverses.  We prove here the results on stably matricial functions which make such uses of inverse function theorems possible.

\bigskip
\noindent
{\bf Lemma.} {\em 
Let $\Xi$ be an open stably matricial $B$-set and let $f = (f_n)_{n \in \bN}$ be an $\cX$-valued continuous stably matricial function on $\Xi$ ($\cX$ a Banach space).  If $f$ is injective then also its fully matricial extension $\wt{f}$ to $\wt{\Xi}$ is injective.
}

\bigskip
\noindent
{\bf {\em Proof.}}  Assume $\pi,\pi' \in \Xi_n$ and $s,s' \in GL(n;\bC)$ are such that $\wt{f}_n(s \cdot \pi) = \wt{f}_n(s'\cdot \pi')$.  Let $s'' = \bpm
0 & s^{-1}s' \\
s'^{-1}s & 0
\epm \in GL(2n;\bC)$ and let $\pi'' = \pi \oplus \pi' \in \Xi_{2n}$.  We have $\wt{f}_{2n}(s'' \cdot \pi'') = \wt{f}_{2n}(\pi'')$ or equivalently $(\Ad s'')f_{2n}(\pi'') = f_{2n}(\pi'')$.  This implies that the matrix $f_{2n}(\pi'')$ with entries in $\cX$ commutes with the scalar matrix $s''$.  It is easily seen that this implies $f_{2n}(\pi'')$ commutes with any matrix $\s$ in the algebra generated by $I_{2n}$ and $s''$.  We may then choose $\s$ to be a logarithm of $s''$.  We will then have
\[
(\Ad(\exp t\s))f_{2n}(\pi'') = f_{2n}(\pi'')
\]
for all $t \in \bC$.

Since $\Xi$ is open, there is $\e > 0$ such that $(\exp t\s) \cdot \pi'' \in \Xi$ if $|t| < \e$.  Since $f$ is injective we infer that $(\exp t\s) \cdot \pi'' = \pi''$ if $|t| < \e$.  If $m \in \bN$ is such that $m > \e^{-1}$ then we have $s'' \cdot \pi'' = (\exp m^{-1}\s) \cdot \pi'' = \pi''$.  This in turn implies that $s^{-1}s' \cdot \pi' = \pi$ and hence we have $s \cdot \pi = s' \cdot \pi'$.\qed

\bigskip
\noindent
{\bf Proposition.} {\em 
Let $\Xi$ be an affine open stably matricial $B$-set and let $f$ be a $B$-valued injective stably matricial function on $\Xi$.  Assume $f$ is continuous and open.  Then $\wt{f}(\wt{\Xi}) = (\wt{f}_n(\wt{\Xi}_n))_{n \in \bN}$ is an open fully matricial $B$-set and each map $\wt{f}_n$ is a homeomorphism of $\wt{\Xi}_n$ onto $\wt{f}_n(\wt{\Xi}_n)$.  Moreover the inverse map $\wt{f}^{-1} = (\wt{f}_n^{-1})_{n \in \bN}$ is fully matricial.  If additionally $f$ is analytic, then so is $\wt{f}^{-1}$.
}

\bigskip
\noindent
{\bf {\em Proof.}} The lemma guarantees that each $\wt{f}_n$ is a bijection of $\wt{\Xi}_n$ onto $\wt{f}_n(\wt{\Xi}_n)$.  Since $f$ is continuous and open and $\wt{f}_n$ is $GL(n;\bC)$-equivariant it follows $\wt{f}_n$ is a homeomorphism onto an open set.  Also since $\wt{f}$ preserves the structures involved (direct sum, equivariance, topology) its inverse is also fully matricial.  In the analytic case $\wt{f}^{-1}$ is clearly analytic.\qed

\bigskip
A typical application of the preceding proposition would run as follows.  If $\cX$ is a Banach space and $f: \{x \in \cX \mid \|x\| < R\} \to \cX$, $f(0) = 0$, is a holomorphic map so that $\|f(x)\| \le C$ when $\|x\| < R$ and $\|(Df)(0)^{-1}\| = M < \i$, then the usual inverse function results give that there exist $\e_1 > 0$ and $\e_2 > 0$ which depend only on $R^{-1}CM$ such that the restriction of $f$ to $\{x \in \cX \mid \|x\| < \e_1R\}$ is an analytic isomorphism onto an open subset of $\cX$ which contains $\{x \in \cX \mid \|x\| < \e_2 RM^{-1}\}$.  (Replacing $f$ by $R^{-1}((Df)(0))^{-1} \circ f(R\cdot)$ this reduces to the case $R = 1$ and $(Df)(0) = I$.)  Assuming $B$ is a $C^*$-algebra (for control of matricial norms) and remarking that the differential of a stably matricial $B$-valued function in $R\cD_0(B)$ has as components at the origin the multiples of the differential of the first component, we can apply the inverse function result to the components of the function and combine this with the Proposition to get the following corollary.

\bigskip
\noindent
{\bf Corollary.} {\em 
Assume $B$ is a $C^*$-algebra and assume $f$ is a $B$-valued analytic stably matricial map on $R\cD_0(B)$ such that $f_1(0) = 0$ and $f(R\cD_0(B)) \subset C\cD_0(B)$.  Suppose $(Df_1)(0)^{-1}$ exists and $\|(Df_1)(0)^{-1}\| < M$, then there are $\e_1 > 0$, $\e_2 > 0$ which depend only on $R^{-1}CM$ so that $\wt{f}$ is a fully matricial isomorphism of $\e_1 R\wt{\cD_0(B)}$ onto an open fully matricial affine $B$-set which contains $\e_2 RM^{-1}\wt{\cD_0(B)}$.  In particular the inverse of $\wt{f}$ is defined on a fully matricial set containing $\e_2 RM^{-1}\wt{\cD_0(B)}$ and is fully matricial holomorphic.
}

\section{More about duality transforms for unitary and hermitian elements of the Grassmannian}
\label{sec12}

\subsection{Universal duality transforms for unitary and hermitian elements}
\label{sec12.1}

The universal duality transforms which we construct in this subsection are the duality transforms for some universal unitary and hermitian Grassmannian elements relative to a given unital $C^*$-algebra $B$.  We will then prove in \ref{sec12.2}--\ref{sec12.5} that certain technical assumptions we made in connection with the duality transform in \ref{sec6.4} and in \ref{sec8.3} are fulfilled in the universal case (\ref{sec8.3} only for the unitary case).

Let $E(B) = B * C(\bT)$ be the unital full free product $C^*$-algebra (that is, the $C^*$-algebra free product with amalgamation over $\bC$ 1) of $B$ and $C(\bT)$.  We shall identify $B$ and $C(\bT)$ with the corresponding $C^*$-subalgebras of $E(B)$.  Identifying $\bT$ with the unit circle in $\bC$, the function $\bT \to \bC$ giving the embedding, is then a unitary element $u \in E(B)$, so that $C(\bT) \ni f \to f(u) \in E(B)$ is the inclusion of $C(\bT)$ into $E$.  We shall also denote by $u$ the element in $Gr_1(E(B))$ corresponding to $u$, that is $\left.\bpm
0 & 1 \\
1 & u
\epm\right/\wt{\l 1}$.  The universal hermitian element is then $\chi = C\left(\bpm
-i & i \\
1 & 1
\epm\right) u = \left. \bpm
i & i(u-1) \\
1 & u+1
\epm\right/\wt{\l 1} \in Gr_1(E(B))$.  If $E$ is a unital $C^*$-algebra with a given unitary element $v$ and containing $B$ as a unital $C^*$-subalgebra, then there is a unique unital $*$-homomorphism $j: E(B) \to E$ acting as the identity on $B$ and so that $j(u) = v$.  Then the duality transform $\cU_v: E^d \to \cA(\wt{\rho}(v;B))$ is so that $\wt{\rho}(u;B) \subset \wt{\rho}(v;B)$ and $\cU_v(\var)|_{\wt{\rho}(v;B)} = \cU_u(\var \circ j)$.  Along the same lines, in the hermitian case, if $E$ is now a unital $C^*$-algebra containing $B$ and $h \in Gr_1(E)$ is a hermitian element, then there is $v \in E$ a unitary element so that $h = C\left( \bpm
-i & i \\
1 & 1
\epm\right) v$ and there is $j: E(B) \to E$, $j|_B = \id_B$, $j(u) = v$.  Then $\cU_h(\var)|_{\wt{\rho}(\chi;B)} = \cU_{\chi}(\var \circ j)$.

Note that several universal objects arise from working in $E(B)$ with $u$ and $\chi$.  The resolvent sets $\wt{\rho}(u;B)$, $\wt{\rho}(\chi;B)$ are universal, in the sense that any resolvents $\wt{\rho}(v;B)$ and $\wt{\rho}(h;B)$ as above contain these.  Also the algebras $\cL\cR(u;B)$ and $\cL\cR(\chi;B)$ are universal as they map into $\cL\cR(v;B)$'s and $\cL\cR(h;B)$'s respectively.

\subsection{A strengthened assumption \ref{sec6.4} is satisfied by the universal unitary element}
\label{sec12.2}

We shall prove assumption \ref{sec6.4} for the universal unitary element $u \in E(B)$ in a stronger form where $\cL\cR(u;B)$ is replaced by the larger algebra $\cQ\cR(u;B)$ which is the inverse closed subalgebra of $E(B)$ generated by $B \cup \cL\cR(u;B)$.  Note that since $0 \in \rho_1(u;B)$ we have $u^{-1} \in \cL\cR(u;B)$ and hence $B\<u\> \subset \cQ\cR(u;B)$.

That assumption \ref{sec6.4} in its original form, holds for the universal unitary element, will then be obtained as a corollary.

\bigskip
\noindent
{\bf Proposition.} {\em
There is a derivation
\[
\p_{u:B}^{\cQ}: \cQ\cR(u;B) \to \cQ\cR(u;B) \otimes \cQ\cR(u;B)
\]
such that
\[
(\id_{\fM_n} \otimes \p_{u:b}^{\cQ})\wt{\cR}_n(u;B)(\s) = \wt{\cR}_n(u;B)(\s) \otimes_{\fM_n} \wt{\cR}_n(u;B)(\s)
\]
for all $n \in \bN$ and $\s \in \wt{\rho}_n(u;B)$.  Moreover we have
\[
\p_{u:b}^{\cQ}|_B = 0 \mbox{ and } \p_{u:B}^{\cQ}u = 1 \otimes 1.
\]
}

\bigskip
\noindent
{\bf {\em Proof.}} Let $\mu$ be a faithful representation of $E(B)$ on some Hilbert space $\cH$.  The universal property of the full free product implies that for any fixed $H = H^* \in \cB(\cH)$ and all $\e \in \bR$ there are representations $\mu_{\e}: E(B) \to \cB(\cH)$ so that $\mu_{\e} \mid B = \mu \mid B$ and $\mu_{\e}(u) = UU(\e)$ where $U = \mu(u)$ and $U(\e) = \exp(i\e H)$.  It will be convenient to identify $B$ with $\mu(B)$ (that is, to assume $B$ is a $C^*$-subalgebra of $\cB(\cH)$).  If $\s = \left. \bpm \a & \b \\ \g & \d \epm\right/\wt{\l n} \in \wt{\rho}_n(u;B)$, then $(\d - (I_n \otimes UU(\e))\b)^{-1}$ exists and
\[
(\id_{\fM_n} \otimes \mu_{\e})(\wt{\cR}_n(u;B)(\s)) = \b(\d - (I_n \otimes UU(\e))\b)^{-1}.
\]
It follows that $\mu_{\e}(x)$ is differentiable at $\e = 0$, as a function of $\e \in \bR$, when $x \in \cL\cR(u;B)$ and since $\mu_{\e}(b)$ is constant for $b \in B$, differentiability actually holds for $x \in B \cup \cL\cR(u;B)$ and hence for $x \in \cQ\cR(u;B)$.  This yields a derivation
\[
d_H: \cQ\cR(u;B) \to \cB(\cH)
\]
with respect to the $\cQ\cR(u;B)$-bimodule structure on $\cB(\cH)$ defined by $\mu$, where
\[
d_H(x) = \left. -i \frac {d}{d\e} \mu_{\e}(x)\right|_{\e = 0}.
\]
We then have
\begin{align*}
&(\id_{\fM_n} \otimes d_H)(\wt{\cR}_n(u;B)(\s)) = \wt{\cR}_n(U;B)(\s)(I_n \otimes UH)\wt{\cR}_n(U;B)(\s) \\
&= (\id_{\fM_n} \otimes \mu)(\wt{\cR}_n(u;B)(\s))(I_n \otimes UH)(\id_{\fM_n} \otimes \mu)(\wt{\cR}_n(u;B)(\s))
\end{align*}
and moreover
\[
d_H|_B = 0 \mbox{ and } d_Hu = UH.
\]
We may extend by complex linearity the map $H \rightsquigarrow d_H$ to $\cB(\cH)$ and get derivations $d_T: \cQ\cR(u;B) \to \cB(\cH)$ with $T \in \cB(\cH)$.  Putting these together yields a derivation
\[
d: \cQ\cR(u;B) \to \cB(\cH)^{\cB(\cH)}
\]
where $d(x) = (d_T(x))_{T \in \cB(\cH)}$ and the $\cQ\cR(u;B)$-bimodule structure on $\cB(\cH)^{\cB(\cH)}$ is deduced from $\mu^{\cB(\cH)}$.  We also have a $\cQ\cR(u;B)$-bimodules map
\[
\var: \cQ\cR(u;B) \otimes \cQ\cR(u;B) \to \cB(\cH)^{\cB(\cH)}
\]
where
\[
\var(x \otimes y) = (\mu(x)UT \mu(y))_{T \in \cB(\cH)}.
\]
It is easily seen that $\var$ is injective and that we have
\[
(\id_{\fM_n} \otimes d)(\wt{\cR}_n(u;B)(\s)) = (\id_{\fM_n} \otimes \var)(\wt{\cR}_n(u;B)(\s) \otimes_{\fM_n} \wt{\cR}_n(u;B)(\s)),
\]
\[
d|_B = 0 \mbox{ and } d(u) = \var(1 \otimes 1).
\]
The formula for $(\id_{\fM_n} \otimes d)(\wt{\cR}_n(u;B)(\s))$ implies that
\[
d(\cL\cR(u;B)) \subset \var(\cL\cR(u;B) \otimes \cL\cR(u;B)).
\]
The set
\[
\{x \in \cQ\cR(u;B) \mid d(x)\in \var(\cQ\cR(u;B) \otimes \cQ\cR(u;B))\}
\]
is an inverse-closed subalgebra of $\cQ\cR(u;B)$ and since it contains $\cL\cR(u;B)$ and $B$ it is equal to $\cQ\cR(u;B)$.  From
\[
d(\cQ\cR(u;B)) \subset \var(\cQ\cR(u;B) \otimes \cQ\cR(u;B))
\]
we infer the existence of a linear map $\p_{u:B}^{\cQ}: \cQ\cR(u;B) \to \cQ\cR(u;B) \otimes \cQ\cR(u;B)$, such that $\var \circ \p_{u:B}^{\cQ} = d$.  Since $d$ is a derivation we easily infer that $\p_{u:B}^{\cQ}$ is also a derivation and the properties of $\p_{u:B}^{\cQ}$ are easily obtained from those of $d$.\qed

\bigskip
\noindent
{\bf Corollary.} {\em 
There is a derivation
\[
\p_{u:B}: \cL\cR(u;B) \to \cL\cR(u;B) \otimes \cL\cR(u;B)
\]
such that
\[
(\id_{\fM_n} \otimes \p_{u:B})\wt{\cR}_n(u;B)(\s) = \wt{\cR}_n(u;B)(\s) \otimes_{\fM_n} \wt{\cR}_n(u;B)(\s)
\]
for all $n \in \bN$ and $\s \in \wt{\rho}_n(u;B)$.
}

\bigskip
\noindent
{\bf {\em Proof.}} We can take
\[
\p_{u:B} = \p_{u:B}^{\cQ} \mid \cL\cR(u;B)
\]
and remark, like we did in the proof of the Proposition that the formula for $(\id_{\fM_n} \otimes \p_{u:b}^{\cQ})(\wt{\cR}_n(u;B))(\s))$ implies that
\[
\p_{u:B}^{\cQ}(\cL\cR(u;B)) \subset \cL\cR(u;B) \otimes \cL\cR(u;B).
\]
\qed

\bigskip
\noindent
{\bf Remark.} It is easy to see that the larger algebra $\fM\cQ\cR(u;B) \supset \cQ\cR(u;B)$, which is the closure under taking entries of the inverses of invertible square matrices with entries in the algebra, is still contained in the domain of definition of the derivations $d_H$.  This yields a proof of the modified Proposition with $\cQ\cR(u;B)$ replaced by $\fM\cQ\cR(u;B)$ and $\p_{u:B}^{\cQ}$ by a derivation
\[
\p_{u:B}^{\fM\cQ}: \fM\cQ\cR(u;B) \to \fM\cQ\cR(u;B) \otimes \fM\cQ\cR(u;B).
\]

\subsection{Assumption \ref{sec8.3} is satisfied by the universal unitary element}
\label{sec12.3}

We shall prove the following result (assumption \ref{sec8.3}) about the universal unitary element $u \in E(B)$.

\bigskip
\noindent
{\bf Proposition.} {\em 
There is a linear map $D: \cL\cR(u;B) \to E(B)$ so that
\[
(\id_{\fM_n} \otimes D)(\wt{\cR}_n(u;B)(\s)) = \left. \frac {d}{dt} \wt{\cR}_n(e^tu;B)(\s)\right|_{t=0}
\]
for all $n \in \bN$ and $\s \in \wt{\rho}_n(u;B)$.
}

\bigskip
\noindent
{\bf {\em Proof.}} Since $E(B)$ is a free product of $B$ and $C(\bT)$, there is a one-parameter group of automorphisms $\a(t)$, $t \in \bR$ of $E(B)$ such that $\a(t)|_B = \id_B$ and $\a(t)(u) = e^{it}u$.  It follows that $\wt{\rho}_n(u;B) = \wt{\rho}_n(\a(t)u;B)$.

On the other hand it is easily seen that $\wt{\cR}_n(e^tu;B)$ is holomorphic in $t$ in a neighborhood of zero.  This gives
\begin{align*}
\left.\frac {d}{dt} \wt{\cR}_n(e^tu;B)(\s)\right|_{t=0}
&= \left. -i \frac {d}{dt} \wt{\cR}_n(\a(t)u;B)(\s)\right|_{t=0} \\
&= \left. -i \frac {d}{dt} (\id_{\fM_n} \otimes \a(t))\wt{\cR}_n(u;B)(\s)\right|_{t=0}
\end{align*}
taking into account that the automorphism of the Grassmannian induced by $\a(t)$ leaves $\s$ fixed since it lies in $\Gr_n(B)$.

By the definition of $\wt{\cR}_n(u;B)(\s)$, since $I_n \otimes u$ is in the domain of definition of the infinitesimal generator of $\id_{\fM_n} \otimes \a(t)$ we infer that $\wt{\cR}_n(u;B)(\s)$ is also in the domain of definition of the infinitesimal generator of $\id_{\fM_n} \otimes \a(t)$.  This in turn gives that the matrix-coefficients of $\wt{\cR}_n(u;B)(\s)$ are in the domain of definition of the infinitesimal generator $P = \left. \frac {d}{dt} \a(t)\right|_{t=0}$ of $\a(t)$.  Thus $\cL\cR(u;B)$ is in the domain of definition of $P$ and we can define $D$ to be $-iP$.\qed

\subsection{The relation of $\cL\cR(\chi;B)$ and $\wt{\cR}_n(\chi;B)$ to $\cQ\cR(u;B)$, $\wt{\cR}_n(u;B)$ and $\wt{\cR}_n(u^{-1};B)$}
\label{sec12.4}

Since $\chi$ was defined as a ``fractional-linear transform'' of $u$ certain results about $\chi$ can be derived from results about $u$.  We collect in this section some technical facts underlying the passage from $u$ to $\chi$.

Since $\chi = C\left(\bpm -i & i \\ 1 & 1 \epm\right)u$ it is immediate that $\wt{\rho}_n(\chi;B) = C_n\left(\bpm -i & i \\ 1 & 1 \epm\right)\wt{\rho}_n(u;B)$.  Also, since there is an automorphism $\g$ of $E(B)$ such that $\g \mid B = \id_B$ and $\g(u) = u^{-1}$, we have $\wt{\rho}_n(u;B) = \wt{\rho}_n(u^{-1};B)$.  On the other hand $C\left(\bpm 0 & 1 \\ 1 & 0 \epm\right)u = u^{-1}$ and hence $C_n\left( \bpm 0 & 1 \\ 1 & 0 \epm\right)\wt{\rho}_n(u;B) = \wt{\rho}_n(u^{-1};B)$ and $C_n\left( \bpm 0 & 1 \\ 1 & 0 \epm\right)\wt{\rho}_n(u;B) = \wt{\rho}_n(u;B)$.

\bigskip
\noindent
{\bf Lemma.} {\em 
Let $\s = \left. \bpm \a & \b \\ \g & \d \epm\right/\wt{\l n} \in \wt{\rho}_n(u;B)$, $\s^{-1} = C_n\left(\bpm 0 & 1 \\ 1 & 0 \epm\right)\s = \left.\bpm \g & \d \\ \a & \b \epm\right/\wt{\l n}$ and let $\nu = C_n\left( \bpm -i & i \\ 1 & 1 \epm\right)\s \in \wt{\rho}_n(\chi;B)$.  We have:

a) $(I_n \otimes u)(\wt{\cR}_n(u;B)(\s))(I_n \otimes u) + I_n \otimes u = -\wt{\cR}_n(u^{-1};B)(\s^{-1})$

b) $\cL\cR(u^{-1};B) \subset \cQ\cR(u;B)$

c) $\wt{\cR}_n(\chi;B)(\nu) = -\frac {i}{2} \wt{\cR}_n(u^{-1};B)(\s^{-1})(I_n \otimes 1 - I_n \otimes u^{-1}) \\ + \frac {i}{2} \wt{\cR}_n(u;B)(\s)(I_n \otimes 1 - I_n \otimes u)$

d) $\cL\cR(\chi;B) \subset \cQ\cR(u;B)$.
}

\bigskip
\noindent
{\bf {\em Proof.}} a) implies b), since $\{u\} \cup \cL\cR(u;B) \subset \cQ\cR(u;B)$.

On the other hand, a) is equivalent to
\[
(I_n \otimes u)\wt{\cR}_n(u;B)(\s) + \wt{\cR}_n(u^{-1};B)(\s^{-1})(I_n \otimes u^{-1}) = -I_n \otimes 1.
\]
This follows from the computation of the left hand side, which is equal to
\begin{align*}
&(I_n \otimes u)\b(\d - (I_n \otimes u)\b)^{-1} + \d(\b - (I_n \otimes u^{-1})\d)^{-1}(I_n \otimes u^{-1}) \\
&= (I_n \otimes u)\b(\d-(I_n \otimes u)\b)^{-1} + \d((I_n \otimes u)\b-\d)^{-1} \\
&= ((I_n \otimes u)\b-\d)(\d-(I_n \otimes u)\b)^{-1} = -I_n \otimes 1.
\end{align*}

Since
\[
\{u,u^{-1}\} \cup \cL\cR(u;B) \cup \cL\cR(u^{-1};B) \subset \cQ\cR(u;B)
\]
we see that d) is a consequence of c).  So, the only thing we still must prove is c).

In view of the definition of $\chi$ and $\nu$, to compute $\wt{\cR}_n(\chi;B)(\s)$ we must first compute the $22$-entry of the block-matrix
\begin{align*}
&\left(\bpm -iI_n \otimes 1 & iI_n\otimes 1 \\
I_n \otimes 1 & I_n \otimes 1
\epm\bpm I_n \otimes 1 & \b \\
I_n \otimes u & \d
\epm\right)^{-1} \\
&= \frac {1}{2} \left( \bpm I_n \otimes 1 & O \\
I_n \otimes u & \d - (I_n \otimes u)\b
\epm\bpm I_n \otimes 1 & \b \\
O & I_n \otimes 1
\epm\right)^{-1} \bpm iI_n \otimes 1 & I_n \otimes 1 \\
-iI_n \otimes 1 & I_n \otimes 1
\epm \\
&= \frac {1}{2} \bpm I_n \otimes 1 & -\b \\
O & I_n \otimes 1
\epm\bpm I_n \otimes 1 & O \\
-(\d-(I_n \otimes u)\b)^{-1}(I_n \otimes u) & (\d - (I_n \otimes u)\b)^{-1}
\epm\bpm iI_n \otimes 1 & I_n \otimes 1 \\
-iI_n \otimes 1 & I_n \otimes 1
\epm \\
&= \bpm * & O \\
-\frac {1}{2} (\d-(I_n \otimes u)\b)^{-1}(I_n \otimes u) & \frac {1}{2}(\d - (I_n \otimes u)\b)^{-1}
\epm \bpm iI_n \otimes 1 & I_n \otimes 1 \\
-iI_n \otimes 1 & I_n \otimes 1 \epm \\
&= \bpm * & * \\
* & \frac {1}{2}(\d-(I_n \otimes u)\b)^{-1}(I_n \otimes (1-u))
\epm .
\end{align*}
To get $\wt{\cR}_n(\chi;B)(\nu)$, we must multiply this entry to the left by the $12$-entry of $\nu$, which is $i(\b-\d)$.  Thus we have
\begin{align*}
&\wt{\cR}_n(\chi;B)(\nu) = \frac {i}{2}(\b-\d)(\d-(I_n \otimes u)\b)^{-1}(I_n \otimes (1-u)) \\
&= \frac {i}{2}(\wt{\cR}_n(u;B)(\s))(I_n \otimes (1-u)) - \frac {i}{2}\d((I_n \otimes u^{-1})\d-\b)^{-1}(I_n \otimes u^{-1})(I_n \otimes (1-u)) \\
&= \frac {i}{2}(\wt{\cR}_n(u;B)(\s))(I_n \otimes (1-u)) - \frac {i}{2}(\wt{\cR}_n(u^{-1};B)(\s^{-1}))(I_n \otimes (1-u^{-1})).
\end{align*}
\qed

\subsection{Assumption \ref{sec6.4} is satisfied by the universal hermitian element $\chi$}\ 
\label{sec12.5}

\bigskip
\noindent
{\bf Proposition.} {\em 
Let $\p_{\chi:B}\xi$ be defined for $\xi \in \cL\cR(\chi;B)$ by
\[
\p_{\chi:B}\xi = \frac {i}{2} (1 \otimes (1-u))(\p_{u:b}^{\cQ}\xi)((1-u) \otimes 1).
\]
Then $\p_{\chi:B}$ is a derivation of $\cL\cR(\chi;B)$ into $\cL\cR(\chi;B) \otimes \cL\cR(\chi;B)$ such that
\[
(\id_{\fM_n} \otimes \p_{\chi:B})\wt{\cR}_n(\chi;B)(\nu) = \wt{\cR}_n(\chi;B)(\nu) \otimes_{\fM_n} \wt{\cR}_n(\chi;B)(\nu)
\]
for all $n \in \bN$ and $\nu \in \wt{\rho}_n(\chi;B)$.
}

\bigskip
\noindent
{\bf {\em Proof.}} Remark that $\p_{\chi:B}$ is well-defined as a derivation of $\cL\cR(\chi;B)$ into $\cQ\cR(u;B)$ in view of the results in \ref{sec12.2} and \ref{sec12.4} and of the fact that the map $\eta \rightsquigarrow -i(1 \otimes (1-u))\eta((1-u) \otimes 1)$ is a $\cQ\cR(u;B)$-bimodule map of $\cQ\cR(u;B) \otimes \cQ\cR(u;B)$ into itself.  It will suffice to prove the formula for $(\id_{\fM_n} \otimes \p_{\chi:B})\wt{\cR}_n(\chi;B)(\nu)$, since this formula completely determines $\p_{\chi:B}$ and establishes that its range lies in $\cL\cR(\chi;B) \otimes \cL\cR(\chi;B)$.

Using Lemma~$12.4$ a) and c) we have:
\begin{align*}
&\wt{\cR}_n(\chi;B)(\nu) = \frac {i}{2}(\wt{\cR}_n(u;B)(\s))(I_n \otimes (1-u)) \\
&+ \frac {i}{2}(I_n \otimes u)(\wt{\cR}_n(u;B)(\s))(I_n \otimes (u-1)) \\
&+ \frac {i}{2}(I_n \otimes (u-1)) \\
&= \frac {i}{2}((I_n \otimes (1-u))(\wt{\cR}_n(u;B)(\s))(I_n \otimes (1-u)) - I_n \otimes (1-u)).
\end{align*}
Since
\begin{align*}
&(\id_{\fM_n} \otimes \p_{\chi:B})(\wt{\cR}_n(\chi;B)(\nu)) \\
&= \frac {i}{2}(I_n \otimes 1 \otimes (1-u))(\id_{\fM_n} \otimes \p_{u:B}^{\cQ})(\wt{\cR}_n(\chi;B)(\nu))(I_n \otimes (1-u) \otimes 1))
\end{align*}
and
\begin{align*}
&(\id_{\fM_n} \otimes \p_{u:B}^{\cQ})(\wt{\cR}_n(\chi;B)(\nu)) \\
&= \frac {i}{2}(-(I_n \otimes 1) \otimes_{\fM_n} ((\wt{\cR}_n(u;B)(\s))(I_n \otimes (1-u))) \\
&-((I_n \otimes (1-u))(\wt{\cR}_n(u;B)(\s)) \otimes_{\fM_n} (I_n \otimes 1) + I_n \otimes 1 \otimes 1 \\
&+((I_n \otimes (1-u))(\wt{\cR}_n(u;B)(\s)) \otimes_{\fM_n} ((\wt{\cR}_n(u;B)(\s))(I_n \otimes (1-u)))
\end{align*}
we get that
\begin{align*}
&(\id_{\fM_n} \otimes \p_{\chi:B})(\wt{\cR}_n(\chi;B)(\nu)) \\
&= \frac {1}{4}((I_n \otimes (1-u)) \otimes_{\fM_n} (I_n \otimes (1-u))(\wt{\cR}_n(u;B)(\s))(I_n \otimes (1-u)) \\
&+(I_n \otimes (1-u))(\wt{\cR}_n(u;B)(\s))(I_n \otimes (1-u)) \otimes_{\fM_n} (I_n \otimes (1-u)) \\
&-(I_n \otimes (1-u)) \otimes_{\fM_n} (I_n \otimes (1-u)) \\
&-(I_n \otimes (1-u))(\wt{\cR}_n(u;B)(\s))(I_n \otimes (1-u)) \\
&\otimes_{\fM_n} (I_n \otimes (1-u))(\wt{\cR}_n(u;B)(\s))(I_n \otimes (1-u))) \\
&= -\frac {1}{4}((I_n \otimes (1-u))\wt{\cR}_n(u;B)(\s)(I_n \otimes (1-u)) - I_n \otimes (1-u)) \\
&\otimes_{\fM_n} ((I_n \otimes (1-u))\wt{\cR}_n(u;B)(\s)(I_n \otimes (1-u)) - I_n \otimes (1-u)) \\
&= -\frac {1}{4} (-2i\wt{\cR}_n(\chi;B)(\nu)) \otimes_{\fM_n} (-2i\wt{\cR}_n(\chi;B)(\nu)) \\
&= (\wt{\cR}_n(\chi;B)(\nu)) \otimes_{\fM_n} (\wt{\cR}_n(\chi;B)(\nu)).
\end{align*}
\qed

\section{The series expansion at the origin of a fully matricial analytic function}
\label{sec13}

\subsection{}
\label{sec13.1} Since we will deal with a local question, we will work in the affine framework and we {\em will use only the Banach space structure of} $B$ (like in section~\ref{sec7}).  We will assume throughout sections \ref{sec13.1}--\ref{sec13.8} that $B$ is finite-dimensional and then in section \ref{sec13.9} give the general nonsense argument extending the result about the series expansion to the case of general $B$.

Let $\O = (\O_n)_{n \ge 1}$ be a fully matricial open $B$-set so that $\O_n \ni 0_n$ (the zero element of $\fM_n(B)$) and let $f = (f_n)_{n \ge 1} \in \cA(\O)$.  The series expansion of $f_n$ at $0_n$ is
\[
f_n(\b) = f_n(0_n) + \sum_{k \ge 1} \frac {1}{k!} d^k f_n (0_n) \underset{\mbox{$k$-times}}{\underbrace{[\b,\dots,\b]}}
\]
where $d^kf_n(0_n)$ is viewed as a symmetric $k$-linear map of $\fM_n(B) \x \dots \x \fM_n(B)$ ($k$-times) into $\fM_n$ and where it is assumed that some disk $\{z\b \mid z \in \bC, |z| < 1 + \e\}$ is in $\O_n$.  Our aim here is to identify the $(d^kf_n(0_n))_{n \ge 1,k \ge 1}$ which occur taking into account the ``fully matricial''---conditions (convergence question aside).

\subsection{}
\label{sec13.2} Let $F_{n,k}(\b) = d^kf_n(0_n)[\b,\dots,\b]$, where $\b \in \fM_n(B)$.

\bigskip
\noindent
{\bf Lemma.} {\em 
We have $(F_{n,k})_{n \ge 1} \in \cA(\fM_n(B))$.
}

\bigskip
\noindent
{\bf {\em Proof.}} The lemma is an immediate consequence of taking $k$-th order differentials at the origin of the two equalities (equivariance and direct sum conditions) which are satisfied by the fully matricial function $f$.\qed

\subsection{}
\label{sec13.3} Since $F_{n,k}$ was derived from $d^kf_n(0_n)$, there are linear maps $\Phi_{n,k}: (\fM_n)^{\otimes k} \otimes B^{\otimes k} \to \fM_n$, satisfying the symmetry condition
\begin{align*}
&\Phi_{n,k}(A_1 \otimes \dots \otimes A_k \otimes b_1 \otimes \dots \otimes b_k) \\
&= \Phi_{n,k}(A_{\s(1)} \otimes \dots \otimes A_{\s(k)} \otimes b_{\s(1)} \otimes \dots \otimes b_{\s(k)})
\end{align*}
for all $\s \in \fS_k$, so that
\[
F_{n,k}(A \otimes b) = \Phi_{n,k}(A \otimes \dots \otimes A \otimes b \otimes \dots \otimes b).
\]
To recover $\Phi_{n,k}$ from $F_{n,k}$ one uses the fact that $k! \Phi_{n,k}(A_1 \otimes \dots \otimes A_k \otimes b_1 \otimes \dots \otimes b_k)$ is the coefficient of $\l_1 \dots \l_k$ in $F_{n,k}(\l_1A_1 \otimes b_1 + \dots + \l_kA_k \otimes b_k)$ which is a polynomial of degree $k$ in $\l_1,\dots,\l_k \in \bC$.  The properties of $(F_{n,k})_{n \ge 1}$ in Lemma~\ref{sec13.2}, then translate immediately into properties of the $\Phi_{n,k}$ which we record in the next lemma.

\bigskip
\noindent
{\bf Lemma.} {\em 
a) We have
\begin{align*}
&\Phi_{n,k}(SA_1S^{-1} \otimes \dots \otimes SA_kS^{-1} \otimes b_1 \otimes \dots \otimes b_k) \\
&= S\Phi_{n,k}(A_1 \otimes \dots \otimes A_k \otimes b_1 \otimes \dots \otimes b_k)S^{-1},
\end{align*}
where $S \in GL(n,\bC)$.

b) If $n = n' + n''$ and if we put
\begin{align*}
&\wt{\Phi}_{n,k}(A_1 \otimes b_1 \otimes A_2 \otimes b_2 \otimes \dots \otimes A_k \otimes b_k) \\
&= \Phi_{n,k}(A_1 \otimes \dots \otimes A_k \otimes b_1 \otimes \dots \otimes b_k),
\end{align*}
we then have
\begin{align*}
&\wt{\Phi}_{n,k}((A'_1 \otimes  b'_1) \oplus (A''_1 \otimes b''_1)) \otimes \dots \otimes ((A'_k \otimes b'_k) \oplus (A''_k \otimes b''_k)) \\
&= \wt{\Phi}_{n',k}(A'_1 \otimes b'_1 \otimes \dots \otimes A'_k \otimes b'_k) \oplus \wt{\Phi}_{n'',k}(A''_1 \otimes b''_1 \otimes \dots \otimes A''_k \otimes b''_k)
\end{align*}
where $A'_j \in \fM_{n'}$, $A''_j \in \fM_{n''}$, $b'_j,b''_j \in B$.
}

\subsection{}
\label{sec13.4} Let $\rho_{n,k} \in \fS_k \to \cL((\bC^n)^{\otimes k}) \simeq \fM_n^{\otimes k}$ be the representation of $\fS_k$ which permutes the $\bC^n$ factors in the $k$-fold tensor product.

\bigskip
\noindent
{\bf Lemma.} {\em 
Let $\Phi_{n,k}: (\fM_n)^{\otimes k} \otimes B^{\otimes k} \to \fM_n$ be linear maps satisfying conditions a) and b) of Lemma~{\em \ref{sec13.3}} and let $\Psi_{n,k}: (\fM_n)^{\otimes (k+1)} \otimes B^{\otimes k} \to \bC$ denote the map $\Psi_{n,k}(A_1 \otimes \dots \otimes A_{k+1} \otimes b_1 \otimes \dots \otimes b_k) = \Tr(\Phi_{n,k}(A_1 \otimes \dots \otimes A_k \otimes b_1 \otimes \dots \otimes b_k)A_{k+1})$.  Then there are $\var_{\s} \in (B^{\otimes k})^d$ so that 
\begin{align*}
&\Psi_{n,k}(A_1 \otimes \dots \otimes A_{k+1} \otimes b_1 \otimes \dots \otimes b_k) \\
&= \sum_{\s \in \fS_{k+1}} \Tr((A_1 \otimes \dots \otimes A_{k+1}) \rho_{n,k+1}(\s)) \var_{\s}(b_1 \otimes \dots \otimes b_k)
\end{align*}
for all $n \in \bN$.  If $n \ge k+1$ the $\var_{\s}$ for which the previous formula holds for that value of $n$ are unique.
}

\bigskip
\noindent
{\bf {\em Proof.}} Condition a) implies that 
\begin{align*}
&\Psi_{n,k}(SA_1S^{-1} \otimes \dots \otimes SA_{k+1}S^{-1} \otimes b_1 \otimes \dots \otimes b_k) \\
&= \Psi(A_1 \otimes \dots \otimes A_{k+1} \otimes b_1 \otimes \dots \otimes b_k).
\end{align*}
Keeping $n$ and $b_1,\dots,b_k$ fixed, this is a linear map $\Tr(\cdot X)$ for some $X \in (\fM_n)^{k+1}$ and we must have $S^{\otimes (k+1)}X(S^{\otimes (k+1)})^{-1} = X$ for $S \in GL(n;\bC)$.  By the theorem of Weyl $X = \sum_{\s \in \fS_{k+1}} c_{\s}\rho_{n,k+1}(\s)$.  It is easy to see that for fixed $n \in \bN$ we can find $\var_{\s} \in (B^{\otimes k})^d$ so that we will have $X = \sum_{\s \in \fS_{k+1}} \var_{\s}(b_1 \otimes \dots \otimes b_k)\rho_{n,k+1}(\s)$ ($\var_{\s}$ is the linear map interpolating between the coefficient $c_{\s}$ for the $b_1 \otimes \dots \otimes b_k$ with the $b_j$'s running over a fixed basis of $B$).

We still have to prove that the functionals $\var_{\s}$ can be chosen independent of $n$.  This can be seen as follows.  On one hand if $n \ge k+1$ the $\rho_{n,k+1}(\s)$, $\s \in \fS_{k+1}$, are linearly independent and the coefficients $c_{\s}$ and functionals $\var_{\s}$ are unique.  On the other hand condition b) in Lemma~\ref{sec13.3} applied to $b'_j = b''_j$ and $A''_j = 0$ shows that
\begin{align*}
&\Phi_{n,k}((A'_1 \oplus 0) \otimes \dots \otimes (A'_k \oplus 0) \otimes b_1 \otimes \dots \otimes b_k) \\
&= \Phi_{n',k}(A'_1 \otimes \dots \otimes A'_k \otimes b_1 \otimes \dots \otimes b_k) \oplus 0_{n''},
\end{align*}
which implies 
\begin{align*}
&\Psi_{n',k}(A'_1 \otimes \dots \otimes A'_{k+1} \otimes b_1 \otimes \dots \otimes b_k) \\
&= \Psi_{n,k}((A'_1 \oplus 0) \otimes \dots \otimes (A'_{k+1} \oplus 0) \otimes b_1 \otimes \dots \otimes b_k).
\end{align*}
Since
\begin{align*}
&\Tr((A'_1 \otimes \dots \otimes A'_{k+1})\rho_{n',k+1}(\s)) \\
&= \Tr(((A'_1 \oplus 0) \otimes \dots \otimes (A'_{k+1} \oplus 0))\rho_{n,k+1}(\s)),
\end{align*}
we infer that the $\var_{\s}$ which work for $n$ can also be used for $n' < n$.  Thus the $\var_{\s}$ in the formula for some $n \ge k+1$ can be used for all $n \in \bN$.\qed

\subsection{}
\label{sec13.5} The result of the preceding section did not draw all the consequences from condition b) of Lemma~\ref{sec13.3} as we shall see in this section.

Remark first that
\[
\Tr((A_1 \otimes \dots \otimes A_{k+1})\rho_{n,k+1}(\s)) = \prod_{\{(i_1,\dots,i_p) \mid (i_1,\dots,i_p) \mbox{\scriptsize{ cycle of }} \s\}} \Tr(A_{i_p} \dots A_{i_1}),
\]
(the product is over disjoint cycles of $\s$, that is without repetitions).  Let $C_{k+1} \subset \fS_{k+1}$ be the cycles of length $k+1$.

\bigskip
\noindent
{\bf Lemma.} {\em 
The conclusion of Lemma~{\em \ref{sec13.4}} holds with the sum taken only over $C_{k+1}$ (that is $\var_{\s} = 0$ if $\s \notin C_{k+1}$).
}

\bigskip
\noindent
{\bf {\em Proof.}} We begin by remarking that Lemma~\ref{sec13.3}~b) implies
\begin{align*}
&\Psi_{nm,k}(\underset{m}{\underbrace{(A_1 \oplus \dots \oplus A_1)}} \otimes \dots \otimes \underset{m}{\underbrace{(A_{k+1} \oplus \dots \oplus A_{k+1})}} \otimes b_1 \otimes \dots \otimes b_k) \\
&= \Tr(\wt{\Phi}_{nm,k}((A_1 \oplus \dots \oplus A_1) \otimes b_1 \otimes \dots (A_k \oplus \dots \oplus A_k) \otimes b_k \otimes (A_{k+1} \oplus \dots \oplus A_{k+1})) \\
&= m \Tr(\Phi_{n,k}(A_1 \otimes \dots \otimes A_k \otimes b_1 \otimes \dots \otimes b_k) A_{k+1}) \\
&= m \Psi_{n,k}(A_1 \otimes \dots \otimes A_{k+1} \otimes b_1 \otimes \dots \otimes b_k).
\end{align*}
On the other hand, if $c(\s)$ is the number of cycles of the permutation $\s$, we have
\begin{align*}
&\Tr(((A_1 \oplus \dots \oplus A_1) \otimes \dots \otimes (A_{k+1} \oplus \dots \oplus A_{k+1}))\rho_{mn}(\s)) \\
&= \prod_{(i_1 \dots i_p) \mbox{\scriptsize{ cycle of }} \s} \Tr((A_{i_p} \oplus \dots \oplus A_{i_p}) \dots (A_{i_1} \oplus \dots \oplus A_{i_1})) \\
&= \prod_{(i_1,\dots,i_p) \mbox{\scriptsize{ cycle of }} \s} m \Tr(A_{i_p} \dots A_{i_1}) \\
&= m^{c(\s)}\Tr((A_1 \otimes \dots \otimes A_{k+1})\rho_n(\s)).
\end{align*}
Hence Lemma~\ref{sec13.4} applied to $nm$ and $\underset{m}{\underbrace{A_j \oplus \dots \oplus A_j}}$ instead of $n$ and $A_j$ gives that
\begin{align*}
&m \Psi_{n,k}(A_1 \otimes \dots \otimes A_{k+1} \otimes b_1 \otimes \dots \otimes b_k) \\
&= \sum_{\s \in \fS_{k+1}} m^{c(\s)} \Tr((A_1 \otimes \dots \otimes A_{k+1})\rho_{n,k+1}(\s)) \var_{\s}(b_1 \otimes \dots \otimes b_k).
\end{align*}
By the uniqueness of the $\var_{\s}$ when $n \ge k+1$ we infer that $\var_{\s} = 0$ if $c(\s) > 1$.\qed

\subsection{}
\label{sec13.6} We shall now translate the result of \ref{sec13.5} in terms of the $\Phi_{n,k}$.  Since $C_{k+1}$ is parametrized by $\fS_k$ by putting each cycle in the form $\pi = (k+1,\theta(k),\dots,\theta(1))$ for some $\theta \in \fS_k$, we have that
\[
\Tr((A_1 \otimes \dots \otimes A_{k+1})\rho_{n,k+1}(\pi)) = \Tr(A_{\theta(1)} \dots A_{\theta(k)} A_{k+1}).
\]

\bigskip
\noindent
{\bf Lemma.} {\em 
Let $\Phi_{n,k}: (\fM_n)^{\otimes k} \otimes B^{\otimes k} \to \fM_n$ be linear maps satisfying conditions a) and b) of Lemma~{\em \ref{sec13.3}}.  Then there are $\psi_{\theta} \in (B^{\otimes k})^d$ for $\theta \in \fS_k$ so that
\[
\Phi_{n,k}(A_1 \otimes \dots \otimes A_k \otimes b_1 \otimes \dots \otimes b_k) = \sum_{\theta \in \fS_k} \psi_{\theta}(b_1 \otimes \dots \otimes b_k)A_{\theta(1)} \dots A_{\theta(k)}
\]
for all $n \in \bN$.  If $n \ge k+1$ the $\psi_{\theta}$ are unique.
}

\bigskip
\noindent
{\bf {\em Proof.}} In view of the remarks preceding the Lemma, this is an immediate consequence of Lemma~\ref{sec13.5} with $\psi_{\theta} = \var_{\pi}$, where $\pi = (k+1,\theta(k),\dots,\theta(1))$.\qed

\subsection{}
\label{sec13.7} Since $B$ is finite-dimensional, we can express $\wt{\Phi}_{n,k}$ in terms of a basis $\var_1,\dots,\var_N$ of $B^d$.  We will use for this the fully matricial functions $z(\var)$ defined in \ref{sec7.1}.

\bigskip
\noindent
{\bf Lemma.} {\em 
Let $\Phi_{n,k}: (\fM_n)^{\otimes k} \otimes B^{\otimes k} \to \fM_n$ be linear maps satisfying conditions a) and b) of Lemma~{\em \ref{sec13.3}} and let $\wt{\Phi}_{n,k}$ be defined as in Lemma~{\em \ref{sec13.3}}.  Then there are numbers $a(\theta;j(1),\dots,j(k))$, where $\theta \in \fS_k$ and $1 \le j(1),\dots,j(k) \le N$, so that 
\begin{align*}
&\wt{\Phi}_{n,k}(A_1 \otimes b_1 \otimes \dots \otimes A_k \otimes b_k) \\
&= \sum_{\theta \in \fS_k} \sum_{\substack{1 \le j(p) \le N \\ 1 \le p \le k}} a(\theta;j(1),\dots,j(k))z(\var_{j(1)})_n(A_{\theta(1)} \otimes b_{\theta(1)}) \dots z(\var_{j(k)})(A_{\theta(k)} \otimes b_{\theta(k)})
\end{align*}
for all $n \in \bN$.  If $n \ge k+1$, the numbers $a(\theta;j(1),\dots,j(k))$ are unique.
}

\bigskip
\noindent
{\bf {\em Proof.}} This is essentially a reformulation of Lemma~\ref{sec13.6}.  We define the numbers $a(\theta;j(1),\dots,j(k))$ so that
\[
\psi_{\theta}(b_1 \otimes \dots \otimes b_k) = \sum_{\substack{1 \le j(p) \le N \\ 1 \le p \le k}} a(\theta;j(1),\dots,j(k))\var_{j(1)}(b_{\theta(1)}) \dots \var_{j(k)}(b_{\theta(k)}).
\]
Remark that given $\psi_{\theta}$, these numbers are  unique, which will give uniqueness of the numbers $a(\dots)$ once the $\psi_{\theta}$ are unique for $n \ge k+1$.  The fact that this lemma is a consequence of \ref{sec13.6} becomes obvious after noting that
\[
z(\var)_n(A \otimes b) = \var(b)A.
\]
\qed

\bigskip
\noindent
{\bf Corollary.} {\em 
If in the preceding lemma the $\wt{\Phi}_{n,k}$ satisfy the symmetry condition
\[
\wt{\Phi}_{n,k}(A_{\a(1)} \otimes b_{\a(1)} \otimes \dots \otimes A_{\a(k)} \otimes b_{\a(k)}) = \wt{\Phi}_{n,k}(A_1 \otimes b_1 \otimes \dots \otimes A_k \otimes b_k)
\]
for all $\a \in \fS_k$, then the numbers $a(\theta;j(1),\dots,j(k))$ will not depend on $\theta$.
}

\bigskip
\noindent
{\bf {\em Proof.}} This follows immediately from the uniqueness part of the lemma.\qed

\subsection{}
\label{sec13.8} We are now ready to return to the context of \ref{sec13.1} and to describe the series expansion of $f \in \cA(\O)$ at the origin.  We will assume $\var_1,\dots,\var_N$ is a basis of $B^d$.

\bigskip
\noindent
{\bf Theorem.} {\em 
If $\O$ is a fully matricial affine $B$-set containing the origin and $f \in \cA(\O)$, then for each $k \in \bN$ there are numbers $a_k(j(1),\dots,j(k))$ so that
\begin{align*}
&d^kf_n(0_n)[\b_1,\dots,\b_k] \\
&= \sum_{\theta \in \fS_k} \sum_{\substack{1 \le j(p) \le N \\ 1 \le p \le k}} a_k(j(1),\dots,j(k))z(\var_{j(1)})_n(\b_{\theta(1)})\dots z(\var_{j(k)})(\b_{\theta(k)}).
\end{align*}
Equivalently if $F_k = (F_{n,k})_{n \ge 1}$ where 
\[
F_{n,k}(\b) = d^kf_n(0_n)[\b,\dots,\b]
\]
then
\[
(k!)^{-1}F_k = \sum_{\substack{1 \le j(p) \le N \\ 1 \le p \le k}} a_k(j(1),\dots,j(k))z(\var_{j(1)}) \dots z(\var_{j(k)}).
\]
The formula for $d^kf_n(0_n)$ determines the numbers $a_k$ uniquely if $n \ge k+1$.  The series expansion of $f$ at the origin is
\[
f = (f_n(0_n))_{n \ge 1} + \sum_{k \ge 1} \sum_{\substack{1 \le j(p) \le N \\ 1 \le p \le k}} a_k(j(1),\dots,j(k)) z(\var_{j(1)}) \dots z(\var_{j(k)}).
\]
}

\bigskip
\noindent
{\bf {\em Proof.}} This follows easily from an application of Lemma~\ref{sec13.7} and Corollary~\ref{sec13.7} to 
\[
\Phi_{n,k}(A_1 \otimes \dots \otimes A_k \otimes b_1 \otimes \dots \otimes b_k) = d^kf_n(0_n)[A_1 \otimes b_1,\dots,A_n \otimes b_n]
\]
(that these $\Phi_{n,k}$ have the required properties is a consequence of \ref{sec13.2} and \ref{sec13.3}).\qed

\bigskip
\noindent
{\bf Remark.} The constant term in the series expansions $(f_n(0_n))_{n \ge 1}$ has a particularly simple form since $f$ is fully maticial $f_n(0_n) = f_1(0) I_n \otimes 1$ or with the notations in \ref{sec7.1} we can write $f_1(0)1\!\!1$ for $(f_n(0_n))_{n \ge 1}$.

\subsection{The series expansion for $\cX$-valued fully matricial analytic functions and general $B$}
\label{sec13.9}

After having dealt with the algebraic questions assuming that $\dim B < \i$ and that the fully matricial analytic function $f$ is scalar-valued, it is now easy to remove the restriction on the dimension and get the result for an arbitrary Banach space $B$ and additionally allow the function $f$ to be $\cX$-valued, where $\cX$ is some other complex Banach space.

If $\a \in \cL(B^{\wh{\otimes} k},\cX)$ where $\wh{\otimes} k$ is the $k$-th projective tensor power of $B$, then $\a$ is a bounded $k$-linear map of $B \x \dots \x B$ ($k$-times) into $\cX$.  We define a fully matricial analytic function $\cZ_k(\a)$ on $\fM(B)$ with values in $\cX$ by
\[
\cZ_k(\a)_n(\b) = (\id_{\fM_n} \otimes \a)(\b \otimes_{\fM_n} \dots \otimes_{\fM_n} \b)
\]
where $\b \in \fM_n(B)$ and the last $\otimes_{\fM_n}$ is $k$-fold.  If $\b = \sum_{1 \le i, j \le n} e_{ij} \otimes b_{ij}$ this means
\begin{align*}
\cZ_k(\a)_n(\b) &= \sum_{1 \le i_1,\dots,i_{k+1} \le n} e_{i_1i_2} \dots e_{i_ki_{k+1}} \otimes \a(b_{i_1i_2} \otimes \dots b_{i_ki_{k+1}}) \\
&= \sum_{1 \le i_1,\dots,i_{k+1} \le n} e_{i_1i_{k+1}} \otimes \a(b_{i_1i_2} \otimes b_{i_2i_3} \otimes \dots \otimes b_{i_ki_{k+1}}).
\end{align*}
It is easily seen that $\cZ_k(\a) = (\cZ_k(\a)_n)_{n \ge 1}$ is a fully matricial analytic function on $\fM(B)$ with values in $\cX$.

If $k = 1$ and $\cX = \bC$ then $z(\var)$ defined in \ref{sec7.1} coincides with $\cZ_1(\var)$.  Also if $\var_1,\dots,\var_k \in B^d$ then $\var_1 \otimes \dots \otimes \var_k \in \cL(B^{\otimes k},\bC)$ and $z(\var_1)\dots z(\var_k) = \cZ_k(\var_1 \otimes \dots \otimes \var_k)$.  If $\dim B < \i$ then Theorem~13.8 is precisely the fact that $(d^kf(0))_n[\b,\dots,\b]$ equals $(\cZ_k(\a))_n(\b)$ for some $\a \in \cL(B^{\otimes k},\bC)$.  Note also that taking $n = k + 1$ and $\b = e_{12} \otimes b_1 + \dots + e_{kk+1} \otimes b_k$ we have
\[
\cZ_k(\a)_{k+1}(\b) = e_{1k+1} \otimes \a(b_1 \otimes \dots \otimes b_k).
\]

\bigskip
\noindent
{\bf Theorem.} {\em 
If $\O$ is a fully matricial affine $B$-set containing the origin and $f$ is a fully matricial $\cX$-valued analytic function on $\O$, then for each $k \in \bN$ there is a unique $\a_k \in \cL(B^{\wh{\otimes} k},\cX)$ so that
\[
d^kf_n(0_n)[\b,\dots,\b] = \cZ_k(\a_k)_n(\b)
\]
for $\b \in \fM_n(B)$.  The series expansion of $f$ at the origin is
\[
f_1(0)1\!\!1 + \sum_{k \ge 1} (k!)^{-1}\cZ_k(\a_k).
\]
}

\bigskip
\noindent
{\bf {\em Proof.}} The case $\dim B < \i$ and $\cX = \bC$ is just Theorem~13.8 as already noted.  Also if $\dim B < \i$ and $N = \dim \cX < \i$ we have that $f = {}_1f \otimes x_1 + \dots + {}_Nf \otimes x_N$ where ${}_jf \in A(\O)$, $1 \le j \le N$ and $x_1,\dots,x_N \in \cX$.  The theorem in this case follows from the theorem applied to each of the ${}_jf$.

Next, assume $\dim B < \i$, but place no restriction on $\cX$.  For each $\var \in \cX^d$ there is $\a_{\var} \in \cL(B^{\otimes k},\bC)$ so that
\[
d^k(\fM(\var) \circ f)(0) = \cZ_k(\a_{\var}).
\]
Both sides being linear in $\var$, since $\a_{\var}$ is unique, we infer $\var \rightsquigarrow \a_{\var}$ is a linear map $\cX^d \to \cL(B^{\otimes k},\bC)$.  It follows that there can be at most $(\dim B)^k$ linearly independent $\a_{\var}$'s and this implies that the linear span in $\cX$ of the union of ranges of $(\th \otimes \id_{\cX})((d^kf)(0))_n$ with $n \in \bN$ is a finite-dimensional subspace $\cX_0 \subset \cX$.  Thus $(d^kf)(0)$ is essentially a fully-matricial $\cX_0$-valued analytic function on $\fM(B)$, which is equal to its $k$-th differential at the origin.  Hence $(d^kf)(0) = \cZ_k(\a)$ where $\a \in \cL(B^{\otimes k},\cX_0)$ is viewed as $\cX$-valued by composition with the inclusion $\cX_0 \subset \cX$.

Finally we pass to the general case, i.e., no restrictions on the dimensions of $B$ and $\cX$.  For each finite dimensional subspace $C$ of $B$ let $f_C = f \mid \fM(C)$ and let $\a_C \in  \cL(C^{\otimes k},\bC)$ be such that $d^kf_C(0) = \cZ_k(\a_C)$.  If $C_1 \subset C_2$ we have $\a_{C_2} \mid C_1^{\otimes k} = \a_{C_1}$.  Putting these together yields a linear map $\a': B^{\otimes k} \to \bC$.  Clearly to conclude the proof it will suffice to show that $\a'$ extends by continuity to $B^{\wh{\otimes} k}$ or equivalently that $\a'$ corresponds to a bounded $k$-linear map $B \x \dots \x B \to \cX$.  If $b_1,\dots,b_k \in B$ let $\b = e_{12} \otimes b_1 + e_{23} \otimes b_2 + \dots + e_{kk+1} \otimes b_k$ and recall the remark preceding the statement of the theorem, that $(\cZ_k(\a_C))_{k+1}(\b) = e_{1k+1} \otimes \a_C(b_1 \otimes \dots \otimes b_k)$ where $C = \bC b_1 + \dots + \bC b_k \subset B$.  It follows that
\[
(d^kf_{k+1}(0))(\b) = (d^kf(0))_{k+1}(\b) = e_{1k+1} \otimes \a'(b_1 \otimes \dots \otimes b_k).
\]
For any crossnorms on $\fM_n(B) = \fM_n \otimes B$ and $\fM_n(\cX) = \fM_n \otimes \cX$ there is a constant $K$ such that
\[
\|(d^kf_{k+1}(0))(\b)\| \le K\|\b\|^k
\]
which gives for some other constant $K'$ that
\[
\|\a'(b_1 \otimes \dots \otimes b_k)\| \le K'(\|b_1\| + \dots + \|b_k\|)^k.
\]
If $\|b_1\| = \dots = \|b_k\|$ we get
\[
\|\a'(b_1 \otimes \dots \otimes b_k)\| \le K'k^k \|b_1\| \dots \|b_k\|.
\]
Obviously this implies the inequality also without the assumption on the equality of norms.\qed

\subsection{The series expansion of a composition of fully matricial analytic functions}
\label{sec13.10}

This section records the fact that the noncommutative series expansion of a composition of fully matricial holomorphic maps coincides with the composition of the series.  More precisely we have the following proposition.

\bigskip
\noindent
{\bf Proposition.} {\em
Let $B(j)$ $(j = 1,2,3)$ be Banach spaces and let $\O(j)$ be affine fully matricial $B(j)$-sets containing the origin.  Let further $g: \O(1) \to \O(2)$ and $f: \O(2) \to \O(3)$ be fully matricial holomorphic maps which send the origin to the origin.  With the notations of {\em \ref{sec13.9}} let
\begin{align*}
\a_k &\in \cL(B(1)^{\wh{\otimes}k},\ B(2)), \\
\b_k &\in \cL(B(2)^{\wh{\otimes}k},\ B(3)), \\
\g_k &: \cL(B(1)^{\wh{\otimes}k},\ B(3))
\end{align*}
be multilinear maps so that
\begin{align*}
(d^kg_n)(O_n)[\b,\dots,\b] &= k!\cZ(\a_k)_n(\b) \\
(d^kf_n)(O_n)[\b,\dots,\b] &= k!\cZ(\b_k)_n(\b) \\
(d^k(f_n \circ g_n))(O_n)[\b,\dots,\b] &= k!\cZ(\g_k)_n(\b).
\end{align*}
Then we have
\[
\g_k = \sum_{\substack{i_1 + \dots + i_l = k \\ i_j \ge 1}} \b_l \circ (\a_{i_1} \otimes \dots \otimes \a_{i_l}).
\]
}

\bigskip
The proposition follows in a straightforward way from the formulae for $d^k(f \circ g)_n$ in terms of the differentials of $f_n$ and $g_n$.  Indeed we have
\[
(k!)^{-1} d^k(f \circ g)_n = \sum_{\substack{1 \le l \le k \\ i_1 + \dots + i_l = k \\ i_j \ge 1}} (l!i_1!\dots i_l!)^{-1} d^lf_n \circ (d^{i_1}g_n \otimes \dots \otimes d^{i_l}g_n)
\]
when applied to $n$ copies of
\[
\b = \sum_{1 \le i, j \le n} e_{ij} \otimes b_{ij}
\]
the left-hand side gives
\[
\sum_{1 \le p_1,\dots,p_{k+1} \le n} e_{p_1p_{k+1}} \otimes \g_k(b_{p_1p_2} \otimes \dots \otimes b_{p_kp_{k+1}})
\]
while the right-hand side is
\[
\sum_{1 \le p_1,\dots,p_k \le n} e_{p_1p_{k+1}} \otimes S_{p_1\dots p_{k+1}}
\]
where
\begin{align*}
S_{p_1\dots p_{k+1}} = &\sum_{\substack{1 \le l \le k \\ i_1 + \dots + i_l = k \\ i_j \ge 1}} \b_l(\a_{i_1}(b_{p_1p_2} \otimes \dots \otimes b_{p_{i_1}p_{i_1+1}}) \cdots \\
&\otimes \a_{i_l}(b_{p_{i_1+\dots+i_{l-1}+1},p_{i_1+\dots+i_{l-1}+2}} \otimes \dots \otimes b_{p_{i_1+\dots+i_l},p_{i_1+\dots+i_l+1}})).
\end{align*}
Comparing the two proves the proposition.

\section{The Asymptotic Integral Formula for the Coefficients in the Unit Disk when $B = \fM_k$}
\label{sec14}

\subsection{}
\label{sec14.1} We will build on the series expansion results to study stably matricial analytic functions in the unit disk.  We will prove stably matricial analogues of the classical facts about density of polynomials and about the identification of the coefficients with the Fourier coefficients of the restriction to the unit circle.  Here the role of the unit circle will be taken by the stably matricial unitary group $\cU(B)$.  Integration on $\cU(B)$ will be the large $N$ limit of integrations with respect to Haar measures.  {\em We will assume that $B = \fM_k(\bC)$ in {\em \ref{sec14.3}} and {\em \ref{sec14.4}}.  The case $B = \bC^k$ will be the subject of section}~\ref{sec16}.

\subsection{Polynomial approximation for totally bounded analytic functions}
\label{sec14.2}\ \ 

\bigskip
\noindent
{\bf Definition.} A $B_2$-valued stably matricial function $f$ on the stably matricial $B_1$-set $\Xi$, where $B_k$ are $C^*$-algebras is {\em totally bounded} if for some $t > 0$, $tf$ factors through the inclusion $\cD_0^{cl}(B_2) \subset \fM(B_2)$.  We define the uniform norm $\|f\|_{\i,\Xi}$ of $f$ to be the infimum of the $t^{-1}$ or equivalently to be
\[
\sup_{n \in \bN} \sup_{\b \in \Xi_n} \|f_n(\b)\|
\]
where $\fM_n(B_2)$ is endowed with the $C^*$-norm.  If $\Xi = R\cD_0(B_1)$ we shall also write $\|f\|_{\i,R}$ for $\|f\|_{\i,\Xi}$.

If $\Xi$ is an open stably matricial set, we shall put
\[
\bH^{\i}(\Xi) = \{f \in \cA(\Xi) \mid \|f\|_{\i,\Xi} < \i\}.
\]

\bigskip
\noindent
{\bf Proposition.} {\em 
If $B$ is a finite-dimensional $C^*$-algebra and $f \in \bH^{\i}(R'\cD_0(B))$ then given $0 < R < R'$ and $\e > 0$, there is $Z \in \cZ(B^d)$ so that
\[
\|f-Z\|_{\i,R} < \e.
\]
}

\bigskip
\noindent
{\bf {\em Proof.}}  If $g: \{z \in \bC \mid |z| < R'\} \to \cX$, where $\cX$ is a Banach space, is an analytic function and if $\sup_{|z| < R'} \|g(z)\| = M$ then we have the well-known bound
\[
\sup_{|z| \le R} \left\| g(z) - \sum_{0 \le n < N} c_nz^n\right\| \le M(R/R')^N(1-R/R')^{-1}
\]
where $\sum_{n \ge 0} c_nz^n$ is the series expansion of $g$ at zero.

Since the series expansion at the origin of $f_n$ restricted to any one-dimensional subspace of $\fM_n(B)$ is the series expansion of the restriction of $f_n$ to that subspace, we infer that if $\b \in \fM_n(B)$, $\|\b\| \le R$ then
\begin{align*}
&\left\|f_n(\b) - f_1(0)I_n - \sum_{1 \le k < N} (k!)^{-1}d^kf_n(0_n)[\b,\dots,\b]\right\| \\
&\le C(R/R')^N(1-R/R')^{-1}
\end{align*}
where $C$ is the uniform norm of $f$ on $R'\cD_0(B)$.

By Theorem~\ref{sec13.8} and Remark~\ref{sec13.8} this means that
\[
\|f-Z\|_{\i,R} \le C(R/R')^N(1-R/R')^{-1}
\]
for a polynomial function $Z \in \bZ(B^d)$.  Choosing $N$ large enough the right-hand side will be $< \e$.\qed

\subsection{The asymptotic integral formula for the coefficients}
\label{sec14.3}

We precede the proof of the integral formula by some free probability preparations.

\bigskip
\noindent
{\bf Lemma.} {\em 
Let $(A,\var)$ be a noncommutative probability space and let $U \in A$ be an invertible element and let $1 \in M \subset A$ be a subalgebra of $A$, so that $M$ and $\{U,U^{-1}\}$ are free.  Assume further that $\var(U^k) = \d_{k,0}$ for $k \in \bZ$ and that $a_1,\dots,a_m,b_1,\dots,b_n$, $c \in M$.  Then we have
\begin{align*}
&\var(a_1Ua_2U\dots a_mUcU^{-1}b_n\dots U^{-1}b_2U^{-1}b_1) \\
&= \begin{cases}
0 &\mbox{if $m \ne n$} \\
\var(a_1b_1)\dots \var(a_mb_m)\var(c) &\mbox{if $m = n$.}
\end{cases}
\end{align*}
}

\bigskip
\noindent
{\bf {\em Proof.}} We shall denote by $\Phi$ the left-hand side of the equality to be proved.  Replacing $U$ by $e^tU$ does not change the assumptions of the lemma, but $\Phi$ changes to $e^{t(m-n)}\Phi$.  It follows that $\Phi = 0$ if $m \ne n$.

Next, observe that it suffices to show that $\var(c) = 0 \Rightarrow \Phi = 0$.  Indeed, assuming we have proved this, the formula for $\Phi$ follows by induction over $m=n$.  Indeed, if $c - \var(c)1 = c'$ we have $\var(c') = 0$ and
\begin{align*}
&\var(a_1U\dots a_nUcU^{-1}b_n \dots U^{-1}b_1) \\
&= \var(c)\var(a_1U \dots a_{n-1}Ua_nb_nU^{-1}b_{n-1} \dots U^{-1}b_1) + 0
\end{align*}
and this equals
\[
\var(c)\var(a_nb_n)\var(a_1b_1) \dots \var(a_{n-1}b_{n-1})
\]
if we have proved the formula up to $n-1$.

To conclude the proof, assume that $\var(c) = 0$ and replace each $a_j$ by $\var(a_j)1 + a'_j$ with $\var(a'_j) = 0$ and similarly $b_j = \var(b_j)1+b'_j$ with $\var(b'_j) = 0$.  This reduces the proof to showing that
\[
\var(a_1U\dots a_nUcU^{-1}b_n \dots U^{-1}b_1) = 0
\]
when $\var(c) = 0$ and $a_j,b_j \in \{1\} \cup (M \cap \ker \var)$.  In this case
\[
\var(a_1U \dots a_nUcU^{-1}b_n \dots U^{-1}b_1)
\]
is the expectation of an alternating product of elements from the two sets $M \cap \ker \var$ and $\{U^k \mid k \in \bZ\backslash\{0\}\}$ and this by freeness equals zero.\qed

\bigskip
Another ingredient in the proof of the integral formula will be an asymptotic freeness result for random Haar unitary matrices which can be found in \cite{16}, \cite{12} or \cite{3}.  We record it here as the next proposition.

\bigskip
\noindent
{\bf Proposition.} {\em 
Let $U_N$ be the tautological function on $U(Nk) \simeq \cU(N;\fM_k)$ with values in $\fM_{Nk}$ viewed as an element of the noncommutative probability space $(A_N,\Phi_N)$, where $A_N = \fM_{Nk}(L^{\i}(U(Nk);d\mu_N))$ and where
\[
\Phi_N(T) = \int_{U(Nk)} (Nk)^{-1}\Tr T d\mu_N
\]
with $\mu_N$ denoting Haar measure.  Let further $\rho_N: \fM_k \to A_N$ be the unital inclusion which identifies $a \in \fM_k$ with the constant matrix $a \otimes I_N \in \fM_k \otimes \fM_N \simeq \fM_{Nk}$.  Then, as $N \to \i$ the sets $\{U_N,U_N^{-1}\}$ and $\rho_N(\fM_k)$ are asymptotically free in $(A_N,\Phi_N)$.
}

\bigskip
To deal with the polynomials $\cZ(B^d)$ when $B = \fM_k$ we will  need to identify certain elements in $\fM_k^d$.  If $A \in \fM_k$ we shall note by $\var_A \in \fM_k^d$ the functional $\var_A(X) = \Tr(XA^t)$ ($t$ denotes the transpose).  In particular if by $e_{pq}$ we denote the matrix-units in $\fM_k$ then we shall also write $\var_{pq}$ for $\var_{e_{pq}}$.

\bigskip
\noindent
{\bf Theorem.} {\em 
If $B = \fM_k$ and $\a_1,\dots,\a_m,\b_1,\dots,\b_n \in \fM_k$ then we have
\begin{align*}
&\lim_{N \to \i} \int_{\cU(N;B)} N^{-1} \Tr(z(\var_{\a_1})(\o)\dots z(\var_{\a_m})_N(\o)(z(\var_{\b_1})_N(\o) \dots z(\var_{\b_n})_N(\o))^*d\mu_N(\o) \\
&= \begin{cases}
0 &\mbox{if $m \ne n$} \\
\prod_{1 \le j \le m} k^{-1} \Tr(\a_j\b_j^*) &\mbox{if $m = n$.}
\end{cases}
\end{align*}
}

\bigskip
\noindent
{\bf {\em Proof.}} It is easily seen that it suffices to prove the theorem in case $\a_i = e_{p_iq_i}$ and $\b_j = e_{r_js_j}$.

Note also that with the notations of the preceding proposition the function
\[
\cU(Nk) \ni \o \to z(\var_{pq})_N(\o)
\]
can be identified with the first $N \x N$ block of $\rho_N(e_{1p})U\rho_N(e_{q1})$.  Using this observation, the left-hand side of the equality we want to prove becomes 
\[
\lim_{N \to \i} k\Phi_N(\rho_N(e_{1p_1})U\rho_N(e_{q_11})\dots \rho_N(e_{1p_m})U\rho_N(e_{q_m1}) \rho_N(e_{1s_n})U^*\rho_N(e_{r_n1})\dots \rho_N(e_{1s_1})U^*\rho_N(e_{r_11})).
\]
By the asymptotic freeness recorded in the proposition the limit is equal to
\[
k\Phi(Ue_{q_1p_2}Ue_{q_2p_3} \dots Ue_{q_ms_n}U^{-1}e_{r_ns_{n-1}}U^{-1}e_{r_{n-1}s_{n-2}} \dots U^{-1}e_{r_1p_1})
\]
where $\Phi$ is the free product trace-state on $C(\bT) * \fM_k$, where $C(\bT)$ is given the Haar state and $\fM_k$ its unique trace-state.  We may now invoke the lemma at the beginning of this section to get that this is further equal to $0$ if $m \ne n$ and if $m = n$ equal to
\begin{align*}
&k^{1-(m+1)}\Tr e_{r_1p_1} \Tr e_{q_ms_m} (\Tr e_{q_{m-1}p_m} e_{r_ms_{m-1}}) \dots (\Tr e_{q_1p_2}e_{r_2s_1}) \\
&= k^{-m}\d_{r_1p_1}\d_{q_ms_m}(\d_{q_{m-1}s_{m-1}}\d_{p_mr_m}) \dots (\d_{q_1s_1}\d_{p_2r_2}).
\end{align*}
On the other hand the right-hand side of the equality to be proved is
\[
k^{-m} \prod_{1 \le j \le m} \Tr e_{p_jq_j}e_{s_jr_j} = k^{-m} \prod_{1 \le j \le m} \d_{p_jr_j}\d_{q_js_j}
\]
which concludes the proof.\qed

\bigskip
\noindent
{\bf Corollary.} {\em 
Assume $f \in \cA((1+\e)\cD_0(\fM_k))$ for some $\e > 0$ is totally bounded.  Then the Taylor expansion of $f$ at the origin is
\[
f = a_01\!\!1 + \sum_{m \ge 1} \sum_{\substack{1 \le p_i,q_i \le k \\ 1 \le i \le m}} a_m(p_1,q_1;\dots;p_m,q_m)z(\var_{p_1q_1})\dots z(\var_{p_mq_m})
\]
where
\[
a_0 = \int_{\cU(N;\fM_k)} N^{-1}\Tr f_N(\o)d\mu_N(\o)
\]
for all $N \ge 1$, and where
\begin{align*}
&a_m(p_1,q_1;\dots;p_m,q_m) \\
&= k^m \lim_{N \to \i} \int_{\cU(N;\fM_k)} N^{-1} \Tr f_N(\o)(z(\var_{p_1q_1})_N(\o)\dots z(\var_{p_mq_m})_N(\o))^*d\mu_N(\o).
\end{align*}
In particular we have
\begin{align*}
&d^mf_N(0_N)[\b,\dots,\b] \\
&= m! \sum_{\substack{1 \le p_i,q_i \le 1 \\ 1 \le i \le m}} a_m(p_1,q_1;\dots;p_mq_m)z(\var_{p_1q_1})_N(\b) \dots z(\var_{p_mq_m})_N(\b).
\end{align*}
}

\bigskip
\noindent
{\bf {\em Proof.}} That the Taylor series of $f$ has the form written in the statement of the corollary for some constants $a_0,a_m(\dots)$ is the content of \ref{sec13.8}.  To check the integral formulae for the coefficients we can reduce the proof to the case when $f$ is given by a polynomial function in $\cZ(B^d)$.  Indeed by the proof of Proposition~\ref{sec14.2} the sequence of polynomial functions given by the $n$-th order Taylor expansion convergences in the norm $\|\cdot\|_{\i,1}$ to $f$ and hence the limits of integrals for these converge to those for $f$.  From polynomial functions in $\cZ(B^d)$ the proof then reduces to the case when $f$ is $1\!\!1$ or some product $z(\var_{p_1q_1})\dots z(\var_{p_mq_m})$.  The statement is then a consequence of the preceding theorem and of the fact that
\[
\int_{\cU(N;B)} z(\var_{p_1q_1})(\o) \dots z(\var_{p_mq_m})(\o)d\mu_N(\o) = 0,
\]
which follows immediately from the invariance of the integration under $\o \to e^{i\theta}\o$.\qed

\subsection{Totally bounded holomorphic functions on $\cD_0(\fM_k)$}
\label{sec14.4}

The norm $\|\ \ \|_{\i,1}$ of totally bounded functions in $\cA(\cD_0(\fM_k))$ is related to a compression of the full free product $C^*$-algebra $\fM_k *_{\bC} C(\bT)$ (i.e., $\fM_k$ and $C(\bT)$ have the same unit element).  More precisely if $e_{ij}$, $1 \le i,j \le k$ are the matrix units of $\fM_k \subset \fM_k *_{\bC} C(\bT)$ and $u \in C(\bT) \subset \fM_k *_{\bC} C(\bT)$ is the unitary element corresponding to the identical function in $C(\bT)$ we shall work with the $C^*$-algebra $B_k = e_{11}(\fM_k *_{\bC} C(\bT))e_{11}$ and use the elements $u_{ij} = e_{1i} ue_{j1} \in B_k$.  The $C^*$-algebra $\fM_k *_{\bC} C(\bT)$ has sufficiently many finite dimensional representations and representations are in bijection with pairs $(\rho,U)$ where $U$ is a unitary operator and $\rho$ a unital $*$-representation of $\fM_k$ on the same Hilbert space (these observations are certainly not new).

\bigskip
\noindent
{\bf Lemma.} {\em 
Let $\b: \cZ(\fM_k^d) \to B_k$ be the unital homomorphism so that $\b(\cZ(\var_{ij})) = u_{ij}$, $1 \le i,j \le k$.  We have
\[
\|\b(f)\| = \|f\|_{\i,1}.
\]
}

\bigskip
\noindent
{\bf {\em Proof.}} Each element of $\cU(N;\fM_k)$ is a unitary operator on $\bC^{Nk} \simeq \bC^N \otimes \bC^k$ on which we have the representation $T \to I_N \otimes T$ of $\fM_k$.  Up to unitary equivalence these are the pairs of unitary operator and representation of $\fM_k$ for all finite-dimensionaal representations of $\fM_k *_{\bC} C(\bT)$.  Remark also that if $f \in \cZ(\fM_k^d)$ and $\o \in \cU(N;\fM_k)$ then $f_N(\o)$ is just the image of $\b(f)$ via the representation of $B_k$ obtained from restricting the representation of $\fM_k * C(\bT)$ corresponding to $\o$ to the subalgebra $B_k$.  Thus
\[
\|\b(f)\| = \sup_N \sup_{\o \in \cU(N;\fM_k)} \|f_N(\o)\|.
\]
Since $f_N$ is holomorphic
\[
\sup_{\o \in \cU(N;\fM_k)} \|f_N(\o)\| = \sup_{\o \in (\cD_0(\fM_k))_N} \|f_N(\o)\|
\]
and hence $\|\b(f)\| = \|f\|_{\i,1}$.\qed

\bigskip
Using the proof of Proposition~\ref{sec14.2} it follows that $\b$ extends by continuity to an isometric homomorphism (which we shall still denote by $\b$) of $\bH^{\i}((1+\e)\cD_0(\fM_k))$ endowed with the $\|\ \ \|_{\i,1}$-norm into $B_k$.  If $f \in \bH^{\i}(\cD_0(\fM_k))$ we apply this result to $f(r\cdot)$ $0 < r < 1$, to obtain the following proposition.

\bigskip
\noindent
{\bf Proposition.} {\em 
Let $f \in \bH^{\i}(\cD_0(\fM_k))$ have the series expansion
\[
a_01\!\!1 + \sum_{m \ge 1} \sum_{\substack{1 \le p_i,q_i \le k \\ 1 \le i \le m}} a_m(p_1,q_1;\dots;p_m,q_m)z(\var_{p_1q_1}) \dots z(\var_{p_mq_m}).
\]
Then for each $0 < r < 1$ the series (summation over $m \ge 1$)
\[
\b(f(r\cdot)) = a_01 + \sum_{m \ge 1} \sum_{\substack{1 \le p_i,q_i \le k \\ 1 \le i \le m}} a_m(p_1,q_1;\dots;p_m,q_m)u_{p_1q_1} \dots u_{p_mq_m}
\]
is convergent in $B_k$ and we have
\[
\|\b(f(r\cdot))\| = \|f(r\cdot)\|_{\i,1}
\]
and
\[
\|f\|_{\i,1} = \sup_{0 < r < 1} \|\b(f(r\cdot))\|.
\]
}

\subsection{Extending the coefficient formula to $\bH^{\i}(\cD_0(\fM_k))$}
\label{sec14.5}

Some standard arguments about boundary values of bounded holomorphic functions can be used to extend the Corollary of section~\ref{sec14.3} to functions in $\bH^{\i}(\cD_0(\fM_k))$.  By basic facts about boundary values for holomorphic functions (see 6.10 in ch.~III of \cite{9}) given $f \in \bH^{\i}(\cD_0(\fM_k))$ for each $N \ge 1$, there is
\[
\b_N(f_N) \in L^{\i}(\cU(N;\fM_k);\mu_N) \otimes \fM_N,
\]
where $\mu_N$ is Haar measure on $\cU(N;\fM_k) \simeq U(Nk)$, so that $\b_N(f_N)(u) = \lim_{r \uparrow 1} f_N(ru)$ for $\mu_N$---almost all $u \in \cU(N;\fM_k)$.

If $f \in \cA(t\cD_0(\fM_k))$ for some $t > 0$, it will be convenient to denote by $T_m(f)$ and $P_m(f)$ the Taylor polynomial of order $m$ of $f$ and respectively its leading term (i.e., $(m!)^{-1}$ times the $m$-th order differential) viewed as elements of $\cZ(\fM_k^d)$, so that $T_m(f) = P_0(f) + P_1(f) + \dots + P_m(f)$.

\bigskip
\noindent
{\bf Theorem.} {\em 
Assume $f \in \bH^{\i}(\cD_0(\fM_k))$.  Then the Taylor expansion of $f$ at the origin is
\[
f = a_01\!\!1 + \sum_{m \ge 1} \sum_{\substack{1 \le p_i,q_i \le k \\ 1 \le i \le m}} a_m(p_1,q_1;\dots;p_m,q_m)z(\var_{p_1q_1}) \dots z(\var_{p_mq_m})
\]
where
\[
a_0 = \int_{\cU(N;\fM_k)} N^{-1}\Tr \b_N(f_N)(\o)d\mu_N(\o)
\]
for all $N \ge 1$, and where
\begin{align*}
&a_m(p_1,q_1;\dots;p_m,q_m) \\
&= k^m \lim_{N \to \i} \int_{\cU(N;\fM_k)} N^{-1} \Tr(\b_N(f_N)(\o)(z(\var_{p_1q_1})_N(\o)\dots z(\var_{p_mq_m})(\o))^*)d\mu_N(\o).
\end{align*}
}

\bigskip
\noindent
{\bf {\em Proof.}} It will be convenient to use the following notations in this proof:
\[
Q_N(\o) = z(\var_{p_1q_1})_N(\o) \dots z(\var_{p_mq_m})(\o)
\]
and if $F$ is an integrable $\fM_N$-valued function on $\cU(N;\fM_k)$
\[
\wt{\Phi}_N(F(\o)) = N^{-1} \int_{\cU(N;\fM_k)} \Tr F(\o)d\mu_N(\o).
\]

We shall apply the corollary in \ref{sec14.2} to $f(r\cdot) \in \bH^{\i}(r^{-1}\cD_0(\fM_k))$.  The Taylor expansion of $f(r\cdot)$ being
\[
a_01\!\!1 + \sum_{m \ge 1} \sum_{\substack{1 \le p_i,q_i \le k \\ 1 \le i \le m}} r^ma_m(p_1,q_1;\dots;p_m,q_m)z(\var_{p_1q_1})\dots z(\var_{p_mq_m})
\]
we have
\[
a_0 = \wt{\Phi}_N(f_N(r\o))
\]
for all $N \ge 1$ and which as $r \uparrow 1$ gives the formula for $a_0$ using $\b_N(f_N)$.  For $m \ge 1$ we get
\[
r^ma_m(p_1,q_1;\dots;p_mq_m) = k^m \lim_{N \to \i} \wt{\Phi}_N(f_N(r\o)(Q_N(\o))^*).
\]
To prove the theorem we must show that
\[
\lim_{r \uparrow 1} \lim_{N \to \i} \wt{\Phi}_N(f_N(r\o)(Q_N(\o))^*) = \lim_{N \to \i} \wt{\Phi}_N(\b_N(f_N)(\o)(Q_N(\o))^*).
\]
The integral in the left-hand side for fixed $0 < r < 1$, because of the uniform convergence of the Taylor series equals
\[
\sum_{n \ge 0} \wt{\Phi}_N(P_n(f)_N(r\o)(Q_N(\o))^*).
\]
The terms with $n - m \ne 0$ vanish since the integrand is homogeneous of degree $n-m$ with respect to $\o \to e^{i\theta}\o$.  Thus the left-hand side equals 
\begin{align*}
&\lim_{r \uparrow 1} \lim_{N \to \i} \wt{\Phi}_N(P_m(f)_N(r\o)(Q_N(\o))^*) \\
&= \lim_{r \uparrow 1} \lim_{N \to \i} r^m \wt{\Phi}_N(P_m(f)_N(\o)(Q_N(\o))^*) \\
&= \lim_{N \to \i} \wt{\Phi}_N(P_m(f)_N(\o)(Q_N(\o))^*).
\end{align*}
On the other hand using dominated convergence the right-hand side of the equality we want to prove equals
\[
\lim_{N \to \i} \lim_{r \uparrow 1} \wt{\Phi}_N(f_N(r\o)(Q_N(\o))^*)
\]
which by the same homogeneity argument used for the left-hand side equals
\begin{align*}
&\lim_{N \to \i} \lim_{r \uparrow 1} \wt{\Phi}_N(P_m(f)_N(r\o)(Q_N(\o))^*) \\
&= \lim_{N \to \i} \wt{\Phi}_N(P_m(f)_N(\o)(Q_N(\o))^*)
\end{align*}
which concludes the proof.\qed

\section{The Large $N$ Limit of a Totally Bounded Holomorphic Function when $B = \fM_k$}
\label{sec15}

\subsection{}
\label{sec15.1} We will show in this section that there is a large $N$ limit for the boundary values on $\cU(\fM_k)$ of a function in $\bH^{\i}(\cD_0(\fM_k))$, which is an operator in a certain subalgebra of a $II_1$-factor.  The asymptotic integral formulae for coefficients give rise to formulae involving the large $N$ limit for the coefficients of the function.

\subsection{The limit algebras $\cL_k^{\i}$ and $\cH_k^{\i}$}
\label{sec15.2}

Let $L^{\i}(\bT) * \fM_k$ be the von~Neumann algebra with the free-product trace-state $\mu * k^{-1}\Tr$, where $\mu$ is Haar measure and let $U$ be the Haar unitary arising from the identical function in $L^{\i}(\bT)$ and consider the matrix units $e_{ij}$ from $\fM_k$, $1 \le i,j \le k$, viewed as elements of the free product.  We define $\cL_k^{\i}$ to be the compression $e_{11}(L^{\i}(\bT) * \fM_k)e_{11}$ which is generated by the elements $v_{ij} = e_{1i}Ue_{j1}$, $1 \le i,j \le k$.  With respect to the trace-state $\psi$ on $\cL_k^{\i}$ (which is a $II_1$-factor), we have that the elements $\{v_{p_1q_1}\dots v_{p_mq_m} \mid m \ge 1, 1 \le p_j,q_j \le k, 1 \le j \le m\}$ form an orthogonal family in $L^2(\cL_k^{\i},\psi)$ and $|v_{p_1q_1}\dots v_{p_mq_m}|_2 = k^{-m/2}$ (an immediate consequence of the lemma in \ref{sec14.3}).

We define $\cH_k^{\i}$ to be the weakly closed subalgebra of $\cL_k^{\i}$ generated by the $v_{ij}$ $1 \le i,j \le k$ (the weak topology is with respect to the standard form of $\cL_k^{\i}$).

Since Haar measure $\mu$ on $\bT$ is invariant under rotations we get that these yield automorphisms of the free product $L^{\i}(\bT) * \fM_k$ which act as the identity on $\fM_k$.  We infer the existence of automorphisms $\a(e^{i\theta})$ of $\cL_k^{\i}$ so that $\a(e^{i\theta})(v_{pq}) = e^{i\theta}v_{pq}$ for $1 \le p,q \le k$.  Also $\cH_k^{\i}$ is invariant under the $\a(e^{i\theta})$.

For $0 \le r < 1$ let $P_r(e^{i\theta}) = (1-r^2) |e^{i\theta}-r|^{-2}$ be the Poisson kernel and if $x \in L^2(\cL_k^{\i},\psi)$ let
\[
\g(r)(x) = (2\pi)^{-1} \int_0^{2\pi} P_r(e^{i\theta})\a(e^{i\theta})(x)d\theta
\]
which is defined in the $L^2$-sense.  Then $\g(r)$ is an ultra-weakly continuous unital completely positive and completely contractive map $\g(r): \cL_k^{\i} \to \cL_k^{\i}$.  It is also easily seen that $\g(r)$ leaves $\cH_k^{\i}$ invariant and that $\g(r)(v_{p_1q_1} \dots v_{p_mq_m}) = r^mv_{p_1q_1} \dots v_{p_mq_m}$.  Moreover $\lim_{r \uparrow 1} |\g(r)(x)-x|_2 = 0$ for all $x \in L^2(\cL_k^{\i},\psi)$.

\subsection{The large $N$ limit map $\b_{\i}$}
\label{sec15.3}

The construction of the large $N$ limit homomorphism $\b_{\i}: \bH^{\i}(\cD_0(\fM_k)) \to \cH_k^{\i}$ will be done in several steps.

\bigskip
\noindent
{\bf Step 1.} On $\cZ(\fM_k^d)$ we define $\b_{\i}$ as the linear map such that
\[
\b_{\i}(z(\var_{p_1q_1}) \dots z(\var_{p_mq_m})) = v_{p_1q_1} \dots v_{p_mq_m}
\]
and $\b_{\i}(1\!\!1) = 1$.  From the large $N$ limit result in the proposition of \ref{sec14.3} it follows easily that $\|\b_{\i}(f)\| \le \|f\|_{\i,1}$ for $f \in \cZ(\fM_k^d)$.

\bigskip
\noindent
{\bf Step 2.} Since $\b_{\i}$ is contractive we can extend its definition to a contractive homomorphism of the closure of $\cZ(\fM_k^d)$ in $\|\cdot\|_{\i,1}$-norm.  This Banach subalgebra of $\bH^{\i}(\cD_0(\fM_k))$ contains in particular all $f(r\cdot)$ where $f \in \bH^{\i}(\cD_0(\fM_k))$ and $0 < r < 1$ in view of the proof of the proposition in \ref{sec14.2}.

\bigskip
\noindent
{\bf Step 3.} If $f \in \cZ(\fM_k^d)$ and if 
\[
a_01\!\!1 + \sum_{m \ge 1} \sum_{\substack{1 \le p_i,q_i \le k \\ 1 \le i \le m}} a_m(p_1,q_1;\dots;p_m,q_m)z(\var_{p_1q_1})\dots z(\var_{p_mq_m})
\]
is its (finite) series expansion then
\begin{align*}
&|a_0|^2 + \sum_{m \ge 1} \sum_{\substack{1 \le p_i,q_i \le k \\ 1 \le i \le m}} |a_m(p_1,q_1;\dots;p_m,q_m)|^2 k^{-m} \\
&= |\b_{\i}(f)|_2^2 \le \|\b_{\i}(f)\|^2 \le \|f\|_{\i,1}.
\end{align*}
Hence, more generally, if $f$ is in the $\|\ \ \|_{\i,1}$ closure of $\cZ(\fM_k^d)$ we have
\[
\b_{\i}(f) = a_01 + \sum_{m \ge 1} \sum_{\substack{1 \le p_i,q_i \le k \\ 1 \le i \le m}} a_m(p_1,q_1;\dots;p_m,q_m)v_{p_1,q_1} \dots v_{p_m,q_m}
\]
where the series converges in $|\ \ |_2$-norm and
\[
|a_0|^2 + \sum_{m \ge 1} \sum_{\substack{1 \le p_i,q_i \le k \\ 1 \le i \le m}} |a_m(p_1,q_1;\dots;p_m,q_m)|^2k^{-m} \le \|f\|_{\i,1}^2.
\]

\bigskip
\noindent
{\bf Step 4.} If $f \in \bH^{\i}(\fM_k^d)$ we can apply the result of Step~3 to the $f(r\cdot)$ for $0 < r < 1$ and we get that
\[
\|f\|_{1,\i}^2 \ge \|f(r\cdot)\|_{\i,1}^2 \ge |a_0|^2 + \sum_{m \ge 1} \sum_{\substack{1 \le p_i,q_i \le k \\ 1 \le i \le m}} |a_m(p_1,q_1;\dots;p_m,q_m)|^2r^{2m}k^{-m}.
\]
Hence, letting $r \uparrow 1$ we find that
\[
\|f\|_{1,\i}^2 \ge |a_0|^2 + \sum_{m \ge 1} \sum_{\substack{1 \le p_i,q_i \le k \\ 1 \le i \le m}} |a_m(p_1,q_1;\dots;p_m,q_m)|^2k^{-m}.
\]

\bigskip
\noindent
{\bf Step 5.} We define $\b_{\i}(f) \in L^2(\cL_k;\psi)$ for $f \in \bH^{\i}(\cD_0(\fM_k))$, using the result of Step~4, to be
\[
\b_{\i}(f) = a_01 + \sum_{m \ge 1} \sum_{\substack{1 \le p_i,q_i \le k \\ 1 \le i \le m}} a_m(p_1,q_1;\dots;p_m,q_m)v_{p_1q_1}\dots v_{p_mq_m}
\]
which coincides on the $\|\ \ \|_{\i,1}$ closure $\cZ(\fM_k^d)$ with the $\b_{\i}$ already constructed.  Moreover $\|f\|_{\i,1} \ge |\b_{\i}(f)|_2$ for general $f$.  Remark that with this definition
\[
|\b_{\i}(f) - \b_{\i}(f(r\cdot))|_2 \to 0
\]
as $r \uparrow 1$ and since
\[
\|\b_{\i}(f(r\cdot))\| \le \|f(r\cdot)\|_{1,\i} \le \|f\|_{\i,1}
\]
we infer that $\|\b_{\i}(f)\| \le \|f\|_{\i,1}$ and $\b_{\i}(f(r\cdot))$ converges to $\b_{\i}(f)$ $*$-strongly and hence is in $\cH_k^{\i}$.

We will record our conclusions as the following proposition.

\bigskip
\noindent
{\bf Proposition.} {\em 
If $f \in \bH^{\i}(\cD_0(\fM_k))$ then $\b_{\i}(f)$ is defined by the series convergent in $|\ \ |_2$-norm
\[
\b_{\i}(f) = a_01 + \sum_{m \ge 1} \sum_{\substack{1 \le p_i,q_i \le k \\ 1 \le i \le m}} a_m(p_1,q_1;\dots;p_m,q_m)v_{p_1q_1} \dots v_{p_mq_m}
\]
and we have $\b_{\i}(f) \in \cH_k^{\i}$ and $\|\b_{\i}(f)\| \le \|f\|_{\i,1}$.  Moreover
\[
|\b_{\i}(f) - \b_{\i}(f(r\cdot))\|_2 \to 0\
\]
as $r \to \i$ and $\b_{\i}$ is a homomorphism of Banach algebras.
}

\bigskip
\noindent
{\bf {\em Proof.}} The only assertion we still must prove is that $\b_{\i}$ is multiplicative.  This is clear on $\cZ(\fM_k^d)$ by the definition and hence also on its $\|\ \ \|_{\i,1}$-closure.  Thus we have
\[
\b_{\i}((fg)(r\cdot)) = \b_{\i}(f(r\cdot))\b_{\i}(g(r\cdot))
\]
and for $r \uparrow 1$ the equality converges strongly to
\[
\b_{\i}(fg) = \b_{\i}(f)\b_{\i}(g).
\]
\qed

\bigskip
\noindent
{\bf Corollary.} {\em 
If $f \in \bH^{\i}(\cD_0(\fM_k))$ then we have $a_0 = \psi(\b_{\i}(f))$ and
\[
k^{-m}a_m(p_1,q_1;\dots;p_m,q_m) = \psi(\b_{\i}(f)v_{p_mq_m}^* \dots v_{p_1q_1}^*).
\]
}

\subsection{Remark}
\label{sec15.4} It is a natural question whether the formal expression
\[
1 \otimes 1 + \sum_{m \ge 1} \sum_{\substack{1 \le p_i,q_i \le k \\ 1 \le i \le m}} k^mv_{p_1q_1} \dots v_{p_mq_m} \otimes v_{p_mq_m}^* \dots v_{p_2q_2}^*
\]
which plays the role of the Cauchy kernel of the large $N$ limit gives rise to an unbounded affiliated operator of the von~Neumann algebra $\cL_k^{\i} \otimes \cL_k^{\i}$ or perhaps $\cL_k^{\i} \otimes (\cL_k^{\i})^{op}$ (op denotes the opposite algebra).  It is also not clear whether there is a chance that the formal kernel
\[
1\!\!1 \otimes 1 + \sum_{m \ge 1} \sum_{\substack{1 \le p_i,q_i \le k \\ 1 \le i \le m}} z(\var_{p_1q_1}) \dots z(\var_{p_mq_m}) \otimes v_{p_mq_m}^* \dots v_{p_1q_1}^*
\]
exhibits some good analytic properties.  It is also natural to ask whether some group object is involved here, possibly a free quantum group (see \cite{3} and references therein).

\section{The Asymptotic Integral Formula and the Large $N$ Limit in the Unit Disk when $B = \bC^k$}
\label{sec16}

\subsection{}
\label{sec16.1} Here we will work out the analogue of the results of sections~\ref{sec14} and \ref{sec15} in the case of the commutative $C^*$-algebra $\bC^k$, which is at times somewhat simpler.  Often the arguments will be quite similar to those for $B = \fM_k$ and our presentation will be more compressed.  Note that when $B = \bC^k$ the components of the stably matricial unitary group are $\cU(N;\bC^k) \simeq U(N)^k$.

\subsection{The series expansion when $B = \bC^k$}
\label{sec16.2}

The natural basis in $B^d$ when $B = \bC^k$ is $\var_1,\dots,\var_k$, where $\var_j(w_1,\dots,w_k) = w_j$, $1 \le j \le k$ and the Taylor series for a fully matricial holomorphic function can be written
\[
f_1(0)1\!\!1 + \sum_{m \ge 1} \sum_{\substack{1 \le j_p \le k \\ 1 \le p \le m}} a(j_1,\dots,j_m)z(\var_{j_1})\dots z(\var_{j_n}).
\]
Note that an element $\b \in \fM_N(\bC^k)$ identifies with a $k$-tuple $\b = (\b_1,\dots,\b_k) \in (\fM_N(\bC))^k$ and $z(\var_{j_1})_N(\b) \dots z(\var_{j_m})_N(\b) = \b_{j_1}\b_{j_2} \dots \b_{j_m}$.

\subsection{The asymptotic integral formula for the coefficients}
\label{sec16.3}

The asymptotic integral formula for the coefficients is a consequence of our basic asymptotic freeness result for random Haar unitary matrices \cite{10}.  We record the result as the next lemma.

\bigskip
\noindent
{\bf Lemma} (\cite{10}).  {\em 
Let $U_{1,N},\dots,U_{k,N}$ be the $k$ projection functions on $U(N)^k \simeq \cU(N;\bC^k)$ with values in $\fM_N$ viewed as elements of the noncommutative probability space $(A_N,\Phi_N)$, where $A_N = \fM_N(L^{\i}(U(N)^k;d\mu_N))$ and where
\[
\Phi_N(T) = \int_{U(N)^k} N^{-1}\Tr Td\mu_N
\]
with $\mu_N$ denoting Haar measure.  Then $\{U_{1N},U_{1,N}^{-1}\},\dots,\{U_{kN},U_{kN}^{-1}\}$ are asymptotically free as $N \to \i$ in $(A_N,\Phi_N)$.
}

\bigskip
Remarking that $U_{j,N} = z(\var_j)_N$ this immediately implies the asymptotic integral formula for coefficients.  Indeed in the large $N$ limit the $U_{j,N}$ behave like the generating unitaries in a free group algebra with respect to the von~Neumann trace.

\bigskip
\noindent
{\bf Proposition.} {\em 
If $B = \bC^k$ and $i_1,\dots,i_m,j_1,\dots,j_n \in \{1,\dots,k\}$ then
\begin{align*}
&\lim_{N \to \i} \int_{\cU(N;B)} N^{-1} \Tr(z(\var_{i_1})_N(\o) \dots z(\var_{i_m})_N(\o)(z(\var_{j_1})_N(\o) \dots z(\var_{j_n})_N(\o))^*)d\mu_N(\o) \\
&= \begin{cases}
0 &\mbox{if $m \ne n$} \\
\d_{i_1j_1} \dots \d_{i_mj_m} &\mbox{if $m = n$.}
\end{cases}
\end{align*}
If $n = 0$ and $m > 0$ the above integral is zero.
}

\bigskip
The same argument as in the proof of Corollary~\ref{sec14.3} then yields the following corollary.

\bigskip
\noindent
{\bf Corollary.} {\em 
Assume $f \in \bH^{\i}((1+\e)\cD_0(\bC^k))$ for some $\e > 0$.  Then the Taylor expansion of $f$ at the origin is
\[
f = a_01\!\!1 + \sum_{m \ge 1} \sum_{\substack{1 \le j_p \le k \\ 1 \le p \le m}} a(j_1,\dots,j_m)z(\var_{j_1}) \dots z(\var_{j_m})
\]
where
\[
a_0 = \int_{\cU(N;\bC^k)} N^{-1}\Tr f_N(\o)d\mu_N(\o)
\]
for all $N \ge 1$, and where
\begin{align*}
&a(j_1,\dots,j_m) \\
&= \lim_{N \to \i} \int_{\cU(N;\bC^k)} N^{-1} \Tr f_N(\o)(z(\var_{j_1})_N(\o) \dots z(\var_{j_m})_N(\o))^*d\mu_N(\o).
\end{align*}
}

\subsection{Totally bounded holomorphic functions on $\cD_0(\bC^k)$}
\label{sec16.4}

In case $B = \bC^k$ the norm $\|\ \ \|_{\i,1}$ and totally bounded functions in $\cA(\cD_0(\bC^k))$ are connected with the full $C^*$-algebra of the free group on $k$ generators $F_k$.  We shall denote by $u_1,\dots,u_k$ the $k$ generating unitary elements of the full $C^*$-algebra of $F_k$ $C^*(F_k)$.  Unital representations of $C^*(F_k)$ are in bijection with unitary representations of $F_k$ which in turn are in bijection with $k$-tuples of unitary operators on Hilbert spaces and $C^*(F_k)$ has sufficiently many finite-dimensional representations.  This immediately implies the following lemma.

\bigskip
\noindent
{\bf Lemma.} {\em 
Let $\b: \cZ((\bC^k)^d) \to C^*(F_k)$ be the unital homomorphism so that $\b(z(\var_j)) = u_j$, $1 \le j \le k$.  We have
\[
\|\b(f)\| = \|f\|_{\i,1}.
\]
}

\bigskip
Using the proof of Proposition~\ref{sec14.2} it follows that $\b$ extend by continuity to an isometric homomorphism (which we shall still denote by $\b$)of $\bH^{\i}((1+\e)\cD_0(\bC^k))$ endowed with the $\|\ \ \|_{\i,1}$-norm into $C^*(F_k)$.  If $f \in \bH^{\i}(\cD_0(\bC^k))$ we apply this result to $f(r\cdot)$, $0 < r < 1$, to obtain the following proposition.

\bigskip
\noindent
{\bf Proposition.} {\em 
Let $f \in \bH^{\i}(\cD_0(\bC^k))$ have the series expansion
\[
a_01\!\!1 + \sum_{m \ge 1} \sum_{\substack{1 \le j_p \le k \\ 1 \le p \le m}} a(j_1,\dots,j_m)z(\var_{j_1})\dots z(\var_{j_m}).
\]
Then for each $0 < r < 1$ the series (summation over $m \ge 1$)
\[
\b(f(r\cdot)) = a_01 + \sum_{m \ge 1} \left(\sum_{\substack{1 \le j_p \le k \\ 1 \le p \le m}} r^ma(j_1,\dots,j_m)u_{j_1}\dots u_{j_m}\right)
\]
is convergent in $C^*(F_k)$ and we have
\[
\|\b(f(r\cdot))\| = \|f(r\cdot)\|_{1,\i}
\]
and
\[
\|f\|_{\i,1} = \sup_{0 < r < 1} \|\b(f(r\cdot))\|.
\]
}

\bigskip
It is also possible to strengthen the asymptotic integral formula to $\bH^{\i}(\cD_0(\bC^k))$.  Like in the case of $\fM_k$ also for $\bC^k$ there is a boundary value map $\b_N(f_N)(u) = \lim_{r \uparrow 1} f_N(ru)$ for $\mu_N$---almost all $u \in \cU(N;\bC^k) \simeq (U(N))^k$ and $\b_N(f_N) \in \fM_N(L^{\i}(\cU(N;\bC^k);d\mu_N))$ (see ch.~III of \cite{9}).

\bigskip
\noindent
{\bf Theorem.} {\em 
Assume $f \in \bH^{\i}(\cD_0(\bC^k))$.  Then the Taylor expansion of $f$ at the origin is
\[
f = a_01\!\!1 + \sum_{m \ge 1} \sum_{\substack{1 \le j_p \le k \\ 1 \le p \le m}} a(j_1,\dots,j_m)z(\var_{j_1}) \dots z(\var_{j_m})
\]
where
\[
a_0 = \int_{\cU(N;\bC^k)} N^{-1} \Tr \b_N(f_N)(\o)d\mu_N(\o)
\]
for all $N \ge 1$, and where
\begin{align*}
&a(j_1,\dots,j_m) \\
&= \lim_{N \to \i} \int_{\cU(N;\bC^k)} N^{-1} \Tr(\b_N(f_N)(z(\var_{j_1})_N(\o)\dots z(\var_{j_m})_N(\o))^*d\mu_N(\o).
\end{align*}
}

\bigskip
The proof is along the same lines as in the case of $\fM_k$ and will be omitted.

\subsection{The large $N$ limit homomorphism $\b_{\i}$ when $B = \bC^k$}
\label{sec16.5}

In case $B = \bC^k$ the role of $\cL_k^{\i},\cH_k^{\i}$ will be played by the free group $II_1$-factor $L(F_k)$ and by $H(F_k)$ the weakly closed nonselfadjoint subalgebra of $L(F_k)$ generated by $1$ and the unitary operators $\l(g_1),\dots,\l(g_k)$ corresponding to the generators of $F_k$.  Also here there is an action of $\bT$ by automorphisms of $\a(e^{i\theta})$ of $L(F_k)$ where $\a(e^{i\theta})\l(g_j) = e^{i\theta}\l(g_j)$ and this gives rise to the semigroup of completely positive maps $\g(r)$, $0 \le r \le 1$, where
\[
\g(r)(x) = (2\pi)^{-1} \int_0^{2\pi} P_r(e^{i\theta})\a(e^{i\theta})(x)d\theta
\]
with $P_r(e^{i\theta}) = (1-r^2)|e^{i\theta}-r|^{-2}$.  Then $\g(r)(\l(g_{i_1})\dots \l(g_{i_m})) = r^m\l(g_{i_1})\dots \l(g_{i_m})$.

The construction of the homomorphism $\b_{\i}: \bH^{\i}(\cD_0(\bC^k)) \to H(F_k)$ is also similar to the construction in case $B = \fM_k$.

\bigskip
\noindent
{\bf Step 1.} On $\cZ((\bC^k)^d)$ we define $\b_{\i}$ as the linear map so that $\b_{\i}(1\!\!1) = 1$ and $\b_{\i}(z(\var_{j_1}) \dots z(\var_{j_m})) = \l(g_{j_1}) \dots \l(g_{j_m})$.  Comparing with the lemma and proposition in \ref{sec16.3}, we see that $\b_{\i}$ arises from $\b$ via the homomorphism $C^*(F_k) \to L(F_k)$ which passes via the reduced $C^*$-algebra and as such $\|\b_{\i}(f)\| \le \|f\|_{\i,1}$ for $f \in \cZ((\bC^k)^d)$ and this extends by continuity to the norm-closure of $\cZ((\bC^k)^d)$ in $\|\ \ \|_{\i,1}$-norm.  This shows in particular that if $f \in \bH^{\i}(\cD_0(\bC^k))$ then $\b_{\i}(f(r\cdot))$ is defined for $0 < r < 1$ and $\|\b_{\i}(f(r\cdot))\| \le \|f(r\cdot)\|_{\i,1}$.

\bigskip
\noindent
{\bf Step 2.} Like in the case of $\fM_k$ also in the case of $\bC^k$ if $f \in \bH^{\i}(\cD_0(\bC^k))$ has series expansion
\[
a_01\!\!1 + \sum_{m \ge 1} \sum_{\substack{1 \le j_p \le k \\ 1 \le p \le m}} a(j_1,\dots,j_m)z(\var_{j_1}) \dots z(j_m)
\]
we find that
\[
|a_0|^2 + \sum_{m \ge 1} \sum_{\substack{1 \le j_p \le k \\ 1 \le p \le m}} |a(j_1,\dots,j_m)|^2 \le \|f\|_{\i,1}^2.
\]
This is first shown when the sum is finite, then for $f(r\cdot)$ and then in full generality.

\bigskip
\noindent
{\bf Step 3.} If $f \in \bH^{\i}(\cD_0(\bC^k))$ we define $\b_{\i}(f)$ as an element of $L^2(L(F_k),\tau)$ ($\tau$ the unique trace-state) by
\[
\b_{\i}(f) = a_01 + \sum_{m \ge 1} \sum_{\substack{1 \le j_p \le k \\ 1 \le p \le m}} a(j_1,\dots,j_m)\l(g_{j_1})\dots \l(g_{j_m}).
\]
We have $\lim_{r \uparrow 1} |\b_{\i}(f) - \b_{\i}(f(r\cdot))|_2 = 0$ and $\sup_{0 < r < 1} \|\b_{\i}(f(r\cdot))\| \le \|f\|_{\i,1}$ from which we infer that $\b_{\i}(f)$ is bounded and $\|\b_{\i}(f)\| \le \|f\|_{\i,1}$ and $\b_{\i}(f)$ is the $*$-strong limit of $\b_{\i}(f(r\cdot))$ as $r \uparrow 1$.

\bigskip
We conclude that $\b_{\i}$ has properties similar to those in the case of $\fM_k$.

\bigskip
\noindent
{\bf Proposition.} {\em 
If $f \in \bH^{\i}(\cD_0(\bC^k))$ then $\b_{\i}(f)$ is defined by the series convergent in $|\ \ |_2$-norm
\[
\b_{\i}(f) = a_01 + \sum_{m \ge 1} \sum_{\substack{1 \le j_p \le k \\ 1 \le p \le m}} a(j_1,\dots,j_p)\l(g_{j_1})\dots \l(g_{j_m})
\]
and we have $\b_{\i}(f) \in H(F_k)$ and $|\b_{\i}(f) - \b_{\i}(f(r\cdot))|_2 \to 0$ as $r \to \i$ and $\b_{\i}$ is a homomorphism of Banach algebras.  Moreover $a_0 = \tau(\b_{\i}(f))$ and $a(j_1,\dots,j_m) = \tau(\b_{\i}(f)\l(g_{j_m}^{-1}) \dots \l(g_{j_1}^{-1}))$.
}

\section{Unbounded Fully Matricial Analytic Functions}
\label{sec17}

\subsection{}
\label{sec17.1} We saw in the preceding sections that totally bounded fully matricial analytic functions behave well with respect to the large $N$ limit.  In this section we show that without the totally boundedness assumption the components of a fully matricial analytic function may increase wildly.

\bigskip
\noindent
{\bf Theorem.} {\em 
Let $B$ be a finite-dimensional $C^*$-algebra with $\dim B > 1$.  Then there is $f = (f_k)_{k \ge 1} \in \cA(\fM(B))$ so that each $f_k$ is a polynomial function and $\|f\|_{\i,R} = \i$ for all $R > 0$.
}

\bigskip
Underlying the proof of the theorem is the following lemma.

\bigskip
\noindent
{\bf Lemma.} {\em 
Let $B$ be a finite-dimensional $C^*$-algebra, with $\dim B > 1$ and let $\var_1,\dots,\var_n$ be a basis of $B^d$.  Let $m_j = z(\var_{\a(1,j)}) \dots z(\var_{\a(N,j)}) \in \cZ(B^d)$, $1 \le j \le p$ be monomials of the same degree $N$ and which are algebraically free in the free semigroup of all monomials in $z(\var_1),\dots,z(\var_n)$ and let $g = \sum_{\s \in \fS_p} \mbox{sign}(\s)m_{\s(1)} \dots m_{\s(p)}$ be the total antisymmetrization of $m_1\dots m_p$.  If $g = (g_k)_{k \ge 1}$, then each $g_k$ is a polynomial function and $g_k \equiv 0$ if $k^2 < p$, but there is some $N$ such that $g_N$ is not identically zero.
}

\bigskip
\noindent
{\bf {\em Proof of the lemma.}} If $\b \in \fM_k(B)$ then $g_k(\b)$ is the total antisymmetrization of the product of $p$ matrices of size $k \x k$.  Since the vector space of $k \x k$ matrices has dimension $k^2$, there are no non-zero antisymmetric $p$-tensors in $\underset{\mbox{$p$-times}}{\underbrace{\fM_k \otimes \dots \otimes \fM_k}}$ if $p > k^2$.  On the other hand $\cZ(B^d)$ is a free algebra with generators $z(\var_1),\dots,z(\var_d)$ so $g$ is non-zero in $\cZ(B^d) \subset \cA(\fM(B))$ and hence for some $N$, $g_N$ is not identically zero.\qed

\bigskip
\noindent
{\bf {\em Proof of the theorem.}} Using the lemma, we can find a sequence $1 < N_1 < N_2 < \dots$ so that $N_{j+1} > N_{j+1}^2$ for all $j \in \bN$ and $g_j = (g_{jk})_{k \ge 1} \in \cZ(B^d) \subset \cA(\fM(B))$ so that $k \le N_j \Rightarrow g_{jk} \equiv 0$ and $g_jN_{j+1}$ is not identically zero ($g_j$ is an $N_{j+1}^2$-antisymmetrization).  Recursively, replacing $g_j$ by some multiple $\l_jg_j$, $\l_j\in \bC$ we can make sure that
\[
\sup_{\substack{\b \in \fM_{N_{j+1}(B)} \\ \|\b\| < 1/j}} \|(g_{1N_{j+1}} + \dots + g_{jN_{j+1}})(\b)\| > j.
\]
Note that $g_{pk} \equiv 0$ if $p > j$ and $k \le N_{j+1}$.  Hence we may define $g = \sum_{1 \le j \le \i} g_j$, which will have the desired properties.\qed

\section*{Appendix I: Duality for the difference quotient bialgebras on the Riemann sphere}
\label{appI}

We present here the duality result for the difference quotient bialgebras on the Riemann sphere, which is perhaps the simplest case of such duality.

Let $K \subset \bC$ be a nonempty compact set and let $G = \bC{\backslash}K$ and $\wt{G} = G \cup \{\i\}= \bP^1(\bC){\backslash}K$ be its complements in $\bC$ and in the Riemann sphere.

Then $\cO(K)$, the algebra of germs of holomorphic functions around $K$ becomes a topological infinitesimal bialgebra when endowed with the comultiplication-derivation given by the difference-quotient $\p f(z_1,z_2) = (z_1-z_2)^{-1}(f(z_1)-f(z_2))$.  Moreover there is a coderivation $Lf = (zf)'$ with respect to $\p$, so that $Lf - f$ is a derivation of the algebra $\cO(K)$.

On the other hand $\cO_{\i}(\wt{G})$, the algebra of holomorphic functions $g: \wt{G} \to \bC$ which vanish at infinity, $g(\i) = 0$, is also a topological infinitesimal bialgebra with a coderivation.  Here, the comultiplication-derivation is given by $-\p$, the negative of the difference-quotient
\[
-\p g(\z_1,\z_2) = -(\z_1-\z_2)^{-1}(g(\z_1)-g(\z_2))
\]
(if one or both $\z_1,\z_2$ is the point at infinity $\p g$ vanishes).  We also define $\L g = (\z g)'$ which is a coderivation and $\L$-id is a derivation.

The well-known duality pairing between $\cO(K)$ and $\cO_{\i}(\wt{G})$ can be described by
\[
\<f \mid g\> = \frac {1}{2\pi i} \int_{\g} f(z)g(z)dz
\]
where $f$ is defined in an  open neighborhood $\o$ of $K$, $g \in \cO_{\i}(\wt{G})$ and $\g$ for each $\o$ is a finite collection of smooth oriented simple nonintersecting curves in $\o$, such that the total winding number of $\g$ is $0$ or $1$ around any point outside its support and for the points of $K$ the winding number is $1$ (for such $\g$ the Cauchy integral formula for points in $K$ holds with integration over $\g$, see for instance \cite{2}).

Then we have the kind of duality relations like in Thm.~5.3 of \cite{13}.

\bigskip
\noindent
{\bf Fact.} If $f,f_1,f_2 \in \cO(K)$ and $g,g_1,g_2 \in \cO_{\i}(\wt{G})$, then
\begin{align*}
\<\p f \mid g_1 \otimes g_2\> &= \<f \mid g_1g_2\>, \\
\<f_1f_2 \mid g\> &= -\<f_1\otimes f_2 \mid \p g\>, \\
\<Lf \mid g\> + \<f \mid \L g\> &= \<f \mid g\>.
\end{align*}

\bigskip
\noindent
{\bf {\em Proof.}}  By density results for holomorphic functions, it suffices to check the first duality relation for the functions $f(z) = (\z - z)^{-1}$ for some $\z \in \bC{\backslash}K$.  After computing $\p f = (\z- z_1)^{-1}(\z - z_2)^{-1}$ we can apply three times the Cauchy integral for the outer region of $\g$ (containing $\i$) to $g_1,g_2$ and $g_1g_2$ and show that the left-hand side is $g_1(\z)g_2(\z)$ and the right-hand side $(g_1g_2)(\z)$.

One can deal along similar lines with the second duality relation.  One first chooses an open neighborhood $\o$ of $K$ where $f_1$ and $f_2$ are defined and a suitable $\g$.  The second duality relation, written as an equality of integrals over $\g$ and $\g \x \g$, respectively, can be reduced by a density argument to the case when $g(\z) = (\z - w)^{-1}$ where $w \in \o$ is a point around which $\g$ has winding number $1$.  After computing $-\p g = (\z_1 - w)^{-1}(\z_2 - w)^{-1}$ we can apply $3$ times the Cauchy integral formula and reduce the equality to the obvious one
\[
(f_1f_2)(w) = f_1(w)f_2(w).
\]

The last duality relation is just that
\begin{align*}
\<Lf \mid g\> + \<f \mid \L g\> &= (2\pi i)^{-1} \int_{\g} ((zf)'g + f(zg)')dz \\
&= (2\pi i)^{-1} \int_{\g} ((zfg)' + fg)dz \\
&= (2\pi i)^{-1} \int_{\g} fgdz.
\end{align*}
\qed

\bigskip
\noindent
{\bf Remark.} In case $K$ is invariant under complex conjugation $z \in K \Rightarrow {\bar z} \in K$, then the same holds for $\wt{G}$ and this yields involutions $f^*(z) = \overline{f({\bar z})}$ in $\cO(K)$ and $g^*(\z) = \overline{g({\bar \z})}$ in $\cO_{\i}(\wt{G})$.  We leave it as an easy exercise to check that with these definitions the bialgebra structures, coderivations and duality are compatible with the involutions along the lines of Thm.~5.3 in \cite{13}.

\section*{Appendix II: The Fully Matricial $R$-transform}
\label{appII}

It is a rather straightforward fact that the operator-valued $R$-transform \cite{11} is actually part of a fully matricial $R$-transform.  We explain this here in the case of a bounded variable, but it is clear that the Grassmannian completion we developed is the appropriate context to deal also with unbounded variables.  We also explain how taking the series expansion of the fully matricial $R$-transform one arrives at the unsymmetrized $R$-transform of \cite{4a}.  In particular the additivity of the unsymmetrized $R$-transform can then be inferred from that of the (symmetric) operator-valued $R$-transform \cite{11}.

We shall assume $1 \in B \subset E$ is a unital inclusion of $C^*$-algebras and $\Phi: E \to B$ is a conditional expectation which is a projection of norm $1$.  Let further $a \in E$ be an element which will be viewed as a $B$-valued noncommutative random variable in $(E,\Phi)$.

Let $G_a(b) = \Phi((1-ba)^{-1}b)$ which is defined for $\|b\| < \|a\|^{-1}$ and which corresponds to $\Phi$ applied to the resolvent of $a$ at $b^{-1}$, when $b$ is invertible.  Since $G_a(0) = 0$ and the differential of $G_a$ at $0$ is the identity, using the inverse function theorem there is a local inverse $L_a$ of $G_a$.  Since we can make $\|\id_B - DG_a(b)\|$ arbitrarily small choosing $\|b\|$ small, it is an exercise for the reader, that there is $\e >0$, depending only on $C > 0$ so that if $\|a\| < C$, then $L_a(b)$ is defined for $\{b \in B \mid \|b\| < \e\}$ and we have $\|L_a(b)\| < (3C)^{-1}$.

Let $H_a(b) = \sum_{k \ge 1} (ab)^{k-1}a$ which is convergent and holomorphic on $\{b \in B \mid \|b\| < C^{-1}\}$.  Then if $\|b\| < C^{-1}$ we have $G_a(b) = \Phi(b+bH(a)b)$.  If $\|b\| < \e$, then $\|L_a(b)\| < C^{-1}$ and we have
\[
b = G_a(L_a(b)) = L_a(b) + L_a(b)\Phi(H_a(L_a(b)))L_a(b).
\]
Then $R_a(b)$ which corresponds to ``$(L_a(b))^{-1} - b^{-1}$'' is defined as
\begin{align*}
R_a(b) &= (1 + \Phi(H_a(L_a(b)))L_a(b))^{-1}\Phi(H_a(L_a(b))) \\
&= \sum_{k \ge 1} (-1)^{k-1}(\Phi(H_a(L_a(b)))L_a(b))^{k-1}\Phi(H_a(L_a(b))).
\end{align*}
Remark that if $\|b\| < \e$ then $\|L_a(b)\|< (3C)^{-1}$ and $\|H_a(L_a(b))\| < C/2$ which gives $\|\Phi(H_a(L_a(b)))L_a(b)\| < 1/6$.  This insures that the series for $R_a(b)$ is convergent.

The same considerations apply to $a \otimes I_n \in \fM_n(E)$ as an element of $(\fM_n(E),\Phi \otimes \id_{\fM_n})$.  Note that the constants $C$ and $\e$ stay the same for all $n \in \bN$.

Note further that the $G_{a \otimes I_n}$ are the components of a stably matricial function on $C^{-1}\cD_0(B)$.  It follows then using \ref{sec11.5} that the inverse functions $L_{a \otimes I_n}$ are the components of a stably matricial function and also the $H_{a \otimes I_n}(L_{a \otimes I_n})$ are the components of a stably matricial function.  This then leads to the same conclusion for the $R_{a \otimes I_n}$ and moreover this function is also totally bounded.  Thus we obtain the following result.

\bigskip
\noindent
{\bf Fact.} If $a \in E$ is a noncommutative $B$-valued random variable in $(E,\Phi)$, then there is $\e > 0$, which depends only on $\|a\|$, such that the $R$-transforms $(R_{a \otimes I_n})_{n \in \bN}$ form a $B$-valued totally bounded stably matricial holomorphic function on $\e\cD_0(B)$.  In particular $(R_{a \otimes I_n})_{n \in \bN}$ has a fully matricial extension to $\wt{\e\cD_0(B)}$.

\bigskip
\noindent
{\bf Remark.} Since $(R_{a \otimes I_n})_{n \in \bN}$ gives rise to a fully matricial function the results about the series expansion apply and the Taylor series has the special form we found in section \ref{sec13}.  If $a_1,a_2 \in E$ are $B$-free in $(E,\Phi)$, then $a_1 \otimes I_n$, $a_2 \otimes I_n$ being $\fM_n(B)$-free in $(\fM_n(E),\fM_n(\Phi)) = (E \otimes \fM_n,\Phi \otimes \id_{\fM_n})$ we infer \cite{11} that $(\cR_{(a_1+a_2) \otimes I_n})_{n \in \bN} = (\cR_{a_1 \otimes I_n})_{n \in \bN} + (\cR_{a_2 \otimes I_n})_{n \in \bN}$ is some $\wt{\a \cD_0(B)}$, $\a > 0$.  Clearly the same additivity property holds for the Taylor series expansions at the origin.  This gives rise to series-versions of the fully matricial $\cR$-transform.  As noted in \cite{4a} the equations defining the operator-valued $R$-transform of \cite{11} and the unsymmetrized $R$-transform formally coincide and hence also formally coincide with the equations defining the fully matricial $R$-transform.  Hence we may invoke \ref{sec13.10} to derive that the series expansion of the fully matricial $R$-transform is the unsymmetrized $R$-transform.

\section*{Acknowledgment}

The author's research was supported in part by NSF Grant DMS-0501178.  The author benefitted from attending and giving lectures at conferences related to the subject of this work, in particular the Free Analysis Workshop at the American Institute of Mathematics in 2006, the Takagi Lectures at Tokyo University in 2007 and the free probability workshops in March 2007 at U.C.\ Berkeley (supported by NSF funds) and at the Fields Institute and BIRS during the academic year 2007--2008.  The author also thanks the referee for useful suggestions concerning Appendix~II.


\begin{thebibliography}{14}

\bibitem{1} M. Aguiar, Infinitesimal Hopf algebras, New Trends in Hopf Algebra Theory (La Falda, 1999) Contemp. Math., vol.~267, American Mathematical Society, Rhode Island, 2000, pp.~1--29.

\bibitem{2} L. V. Ahlfors, Complex Analysis, McGraw--Hill, 1979.

\bibitem{3} T. Banica and B. Collins, Integration over compact quantum groups, Publ. Res. Inst. Math. Sci. {\bf 43} (2007), 277--302.

\bibitem{4} B. Collins, Moments and cumulants of polynomial random variables on unitary groups, Itzykson--Zuber integral and free probability, International Math. Res. Notices (17):953--982, 2003.

\bibitem{4a} K. J. Dykema, Multilinear function series and transforms in free probability theory, Adv. Math. {\bf 208} (2007), 351--407.

\bibitem{5} J. W. Helton, J. A. Ball, C. R. Johnson and J. N. Palmer, Operator Theory, Analytic Functions, Matrices and Electrical Engineering, CBMS Regional Conference Series in Mathematics, vol.~68, American Mathematical Society, Rhode Island, 1987.

\bibitem{6} E. Hille and R. Phillips, Functional Analysis and Semigroups, Colloquium Publications, American Mathematical Society, 1982.

\bibitem{7} S. A. Joni and G.-C. Rota, Coalgebras and bialgebras in combinatorics, Stud. Appl. Math. {\bf 61} (1979), no.~2, 93--139.

\bibitem{8} P. Muhly and B.Solel, Hardy algebras associated with $W^*$-correspondences (point evaluation and Schur class functions), Operator theory, system theory and scattering theory: multidimensional generalizations, Operator Theory: Adv. and Appl. {\bf 157} (2005), 221--241.

\bibitem{8a} R. Speicher, Combinatorial theory of the free product with amalgamation and operator-valued free probability theory, Mem. AMS {\bf 132} (1998), no.~627.

\bibitem{9} E. M. Stein and G. Weiss, Introduction to Fourier Analysis on Euclidean Spaces, Princeton University Press, Princeton, NJ, 1971.

\bibitem{10} D. V. Voiculescu, Limit laws for random matrices and free products, Invent. Math. {\bf 104} (1991), 201--220.

\bibitem{11} D. V. Voiculescu, Operations on certain non-commutative operator-valued random variables, Asterisque (1995), no.~232, 243--275.

\bibitem{12} D. V. Voiculescu, A strengthened asymptotic freeness result for random matrices with applications to free entropy, International Math. Res. Notices (1):41--63, 1998.

\bibitem{13} D. V. Voiculescu, The coalgebra of the free difference quotient and free probability, International Math. Res. Notices {\bf 2000} (2000), no.~2, 79--106.

\bibitem{14} D. V. Voiculescu, Free Analysis Questions I: Duality Transform for the Coalgebra of $\p_{X:B}$, International Math. Res. Notices {\bf 16} (2004), 793--822.

\bibitem{15} D. V. Voiculescu, K. Dykema and A. Nica, Free Random Variables, CRM Monograph Series 1, American Mathematical Society, Providence, Rhode Island, 1992.

\bibitem{16} F. Xu, A random matrix model from two-dimensional Yang--Mills theory, Communications Math. Phys. {\bf 190} (2):287--307, 1997.

\end{thebibliography}
\end{document}